\documentclass [leqno, 10pt]{amsart}

\usepackage[all]{xy}
\usepackage{amssymb, amsmath, amsfonts, amsthm, amsbsy, amscd, amsxtra}
\usepackage[bbgreekl]{mathbbol} 
\usepackage[mathscr]{euscript}

\usepackage{url}
\usepackage[linktocpage]{hyperref}

\newtheorem{theorem}{Theorem}[section]

\newtheorem{corollary}{Corollary}[section]
\newtheorem{lemma}{Lemma}[section]

\newtheorem{definition}{Definition}[section]

\theoremstyle{definition}

\newtheorem{remark}{Remark}[section]

\numberwithin{equation}{section}

\newcommand{\bG}{\mathbb{G}}
\newcommand{\bOmega}{\mathbb{\Omega}}
\newcommand{\sphereast}{\mathbb{S}^{2\,\ast}}
\newcommand{\sphere}{\mathbb{S}^{2}}
\newcommand{\torus}{\mathbb{T}^{2}}
\newcommand{\sing}{\mathscr{S}}

\newcommand{\dlef}{\mathsf{\Lambda}}
\newcommand{\Om}{\Omega}
\newcommand{\lne}{\mathcal{E}}
\newcommand{\tcm}{T_{\com}M}
\newcommand{\cscal}{\mathsf{S}}
\newcommand{\imt}{\iota}
\newcommand{\elp}{\mathscr{C}_{\mathsf{ell}}}

\def\blt#1{\ensuremath{\mkern+0mu{\mkern-0mu ^{\nv}\mkern-3mu#1}}}
\newcommand{\hu}{\blt{h}}
\newcommand{\gu}{\blt{g}}

\newcommand{\hi}{\vphantom{}^{\infty}{h}}

\newcommand{\hsR}{\hat{\sR}}
\newcommand{\hbt}{\hat{\bt}}
\newcommand{\tractor}{\mathcal{T}}
\newcommand{\oper}{\mathscr{O}}
\newcommand{\prebij}{\mathsf{\Psi}}
\newcommand{\bij}{\mathbb{\Psi}}
\newcommand{\bj}{\mathbb{J}}
\newcommand{\bJ}{\mathbb{J}}

\newcommand{\nv}{\nu}

\renewcommand{\H}{\mathsf{H}}

\newcommand{\cano}{\mathscr{K}}
\newcommand{\tf}{\operatorname{\mathsf{tf}}}

\renewcommand{\div}{\mathsf{div}\,}

\newcommand{\symk}{S^{k}(\ctm)}

\newcommand{\symkt}{S^{k}_{0}(\ctm)}

\newcommand{\symkp}{S^{k+1}(\ctm)}

\newcommand{\symkv}{S^{k}(TM)}
\newcommand{\symktv}{S^{k}_{0}(TM)}

\newcommand{\symkptv}{S^{k+1}_{0}(TM)}

\newcommand{\om}{\omega}

\newcommand{\hess}{\operatorname{\mathsf{Hess}}}

\newcommand{\monge}{\mathcal{M}}

\newcommand{\ten}{[\tilde{\nabla}]}
\newcommand{\sth}{{\vphantom{}}^{\ast}\!{h}}
\def\sbrt{\ensuremath{\mkern+6mu{\mkern-6mu ^{\ast}\mkern-3mu}}}
\newcommand{\bsth}{\bar{\sbrt{h}}}

\newcommand{\ahop}{\mathcal{A}}
\newcommand{\map}{\mathsf{Map}}
\newcommand{\prekhodge}{\mathsf{Q}^{k}}
\newcommand{\khodge}{\mathscr{Q}^{k}}
\newcommand{\precubic}{\mathsf{Q}^{3}}
\newcommand{\cubic}{\mathscr{Q}^{3}}
\newcommand{\dproj}{\mathscr{P}}
\newcommand{\jpremod}{\mathsf{J}}
\newcommand{\epremod}{\mathsf{E}}
\newcommand{\emod}{\mathscr{M}}
\newcommand{\tmod}{\mathscr{E}}
\newcommand{\teich}{\mathscr{T}}

\newcommand{\hs}{\star\,}
\newcommand{\hsh}{\star_{h}\,}

\newcommand{\hdelbar}{\hat{\del}}
\newcommand{\tdelbar}{\tilde{\del}}
\newcommand{\hodge}{\square}
\newcommand{\dad}{d^{\ast}}

\newcommand{\vol}{\mathsf{vol}}

\newcommand{\tD}{\tilde{D}}
\newcommand{\delbar}{\bar{\del}}

\newcommand{\ka}{\kappa}

\newcommand{\sR}{\mathscr{R}}

\renewcommand{\part}{\vdash}

\newcommand{\rad}{\mathbb{E}}

\newcommand{\Id}{\operatorname{Id}}

\newcommand{\dum}{\,\cdot\,\,}
\newcommand{\Ga}{\Gamma}

\newcommand{\nm}{\mathsf{W}}

\newcommand{\lap}{\Delta}

\renewcommand{\j}{\mathsf{i}}

\newcommand{\la}{\lambda}
\newcommand{\ep}{\epsilon}

\newcommand{\reat}{\mathbb{R}^{\times}}
\newcommand{\comt}{\mathbb{C}^{\times}}
\newcommand{\reap}{\mathbb{R}^{+}}

\newcommand{\cinf}{C^{\infty}}

\newcommand{\Det}{\operatorname{Det}}

\newcommand{\ben}{[\Bar{\nabla}]}

\newcommand{\hD}{\Hat{D}}

\newcommand{\eno}{\operatorname{End}}

\newcommand{\lin}{\mathsf{L}}

\newcommand{\si}{\sigma}
\newcommand{\pr}{\partial}
\newcommand{\fd}{\mathsf{F}}
\newcommand{\ctm}{T^{\ast}M}

\newcommand{\sign}{\text{sgn}}
\newcommand{\bnabla}{\bar{\nabla}}
\newcommand{\integer}{\mathbb{Z}}

\def\tbf#1{\textbf{#1}\index{#1}}

\newcommand{\en}{[\nabla]}

\newcommand{\lie}{\mathfrak{L}}
\newcommand{\clie}{\mathscr{L}}

\newcommand{\re}{\text{Re\,}}
\newcommand{\im}{\text{Im\,}}

\newcommand{\B}{\mathcal{B}}

\newcommand{\lb}{\langle}
\newcommand{\ra}{\rangle}

\newcommand{\ste}{\mathbb{V}}
\newcommand{\stec}{\mathbb{V}\tensor_{\rea}\com}
\newcommand{\sted}{\mathbb{V}^{\ast}}

\newcommand{\projp}{\mathbb{P}^{+}}

\newcommand{\al}{\alpha}
\newcommand{\be}{\beta}
\newcommand{\ga}{\gamma}

\newcommand{\emf}{\mathcal{V}}

\newcommand{\uR}{\mathsf{R}}
\newcommand{\hnabla}{\widehat{\nabla}}
\newcommand{\tnabla}{\tilde{\nabla}}

\newcommand{\gl}{\mathfrak{gl}}
\newcommand{\sll}{\mathfrak{sl}}

\newcommand{\proj}{\mathbb{P}}

\DeclareMathOperator{\diff}{Diff}
\DeclareMathOperator{\Aut}{Aut}

\newcommand{\g}{\mathfrak{g}}

\newcommand{\del}{\partial}
\newcommand{\tensor}{\otimes}

\newcommand{\rea}{\mathbb R}
\newcommand{\com}{\mathbb C}
\newcommand{\tr}{\operatorname{\mathsf{tr}}}

\newcommand{\bt}{\mathcal{L}}
\newcommand{\nbt}{|\mathcal{L}|^{2}_{H}}
\newcommand{\tnbt}{|\tilde{\mathcal{L}}|^{2}_{H}}

\newcommand{\defeq}{:=}
\def\pwr#1#2{{{^{#1}\mkern-3mu} #2}}

\makeindex

\begin{document}
\title{Einstein-like geometric structures on surfaces}
\author{Daniel J. F. Fox} 
\address{Departamento de Matem\'atica Aplicada\\ EUIT Industrial\\ Universidad Polit\'ecnica de Madrid\\Ronda de Valencia 3\\ 28012 Madrid Espa\~na}
\email{daniel.fox@upm.es}

\begin{abstract}
An AH (affine hypersurface) structure is a pair comprising a projective equivalence class of torsion-free connections and a conformal structure satisfying a compatibility condition which is automatic in two dimensions. They generalize Weyl structures, and a pair of AH structures is induced on a co-oriented non-degenerate immersed hypersurface in flat affine space. The author has defined for AH structures Einstein equations, which specialize on the one hand to the usual Einstein Weyl equations and, on the other hand, to the equations for affine hyperspheres. Here these equations are solved for Riemannian signature AH structures on compact orientable surfaces, the deformation spaces of solutions are described, and some aspects of the geometry of these structures are related. Every such structure is either Einstein Weyl (in the sense defined for surfaces by Calderbank) or is determined by a pair comprising a conformal structure and a cubic holomorphic differential, and so by a convex flat real projective structure. In the latter case it can be identified with a solution of the Abelian vortex equations on an appropriate power of the canonical bundle. On the cone over a surface of genus at least two carrying an Einstein AH structure there are Monge-Amp\`ere metrics of Lorentzian and Riemannian signature and a Riemannian Einstein K\"ahler affine metric. A mean curvature zero spacelike immersed Lagrangian submanifold of a para-K\"ahler four-manifold with constant para-holomorphic sectional curvature inherits an Einstein AH structure, and this is used to deduce some restrictions on such immersions.
\end{abstract}

\maketitle

\setcounter{tocdepth}{1}  
{\footnotesize }
\tableofcontents

\section{Introduction}\label{introduction}
In \cite{Fox-ahs} the author defined a class of geometric structures, called \textbf{AH (\textit{affine hypersurface}) structures}, which generalize both Weyl structures and the structures induced on a co-oriented non-degenerate immersed hypersurface in flat affine space by its second fundamental form, affine normal and co-normal Gau\ss\, map, and defined for these AH structures equations, called \textbf{Einstein}, specializing on the one hand to the usual Einstein Weyl equations and, on the other hand, to the equations for an affine hypersphere. In the present paper these equations are solved on compact orientable surfaces and some aspects of their geometry in this case are described. It is hoped that, aside from their interest as such, the results provide motivation for studying higher dimensional AH structures.

\subsection{}\label{introsection}
An \textbf{AH structure} on an $n$-manifold is a pair $(\en, [h])$ comprising a projective equivalence class $\en$ of torsion-free affine connections and a conformal structure $[h]$ such that for each $\nabla \in \en$ and each $h \in [h]$ there is a one-form $\si_{i}$ such that $\nabla_{[i}h_{j]k} = 2\si_{[i}h_{j]k}$, or, what is the same, the completely trace-free part of $\nabla_{[i}h_{j]k}$ vanishes (most notational and terminological conventions can be found in sections \ref{backgroundsection} and \ref{ahsection}). When $n = 2$ this compatibility condition is automatic; any pair $(\en, [h])$ is AH. On the one hand, an AH structure for which $\nabla_{i}h_{jk}$ is pure trace for any $\nabla \in \en$ and any $h \in [h]$ is simply a Weyl structure (what is usually called the Weyl connection is the \textbf{aligned} representative $\nabla \in \en$ distinguished by the requirement that $h^{pq}\nabla_{p}h_{iq} = nh^{pq}\nabla_{i}h_{pq}$ for any $h \in [h]$). On the other hand, there is induced on any non-degenerate co-orientable immersed hypersurface in flat affine space a pair of AH structures. Namely, for both, $[h]$ is generated by the second fundamental form, while the projective structures are those induced via the affine normal bundle and the conormal Gau\ss\, map. There is a canonical duality associating to each AH structure $(\en, [h])$ a conjugate AH structure $(\ben, [h])$ having the same underlying conformal structure. The Weyl structures are exactly the self-conjugate AH structures, and the two AH structure induced on an affine hypersurface are conjugate in this sense. It may be helpful to think of the formalism of AH structures as giving an instrinsically formulated generalization of the geometry of hypersurfaces in flat affine space in a manner similar to how CR structures abstract the geometry of a pseudoconvex real hypersurface in flat complex Euclidean space. As for CR structures, the generalization is genuine; there are local obstructions to realizing a given AH structure as that induced on a hypersurface in flat affine space. 

The specialization to surfaces of the notion of Einstein AH structure defined in \cite{Fox-ahs} says that an AH structure on a surface is \textbf{Einstein} if there vanishes the trace-free symmetric part of the Ricci curvature of its aligned representative $\nabla$, while the scalar trace of its Ricci curvature satisfies the condition \eqref{conservationcondition}, generalizing the constancy of the scalar curvature of an Einstein metric, and taken by D. Calderbank in \cite{Calderbank-mobius, Calderbank-faraday, Calderbank-twod} as the definition of a two-dimensional Einstein Weyl structure.  Calderbank's point of view in the two-dimensional case was motivating for the definition of Einstein AH equations in general. In all dimensions the notion of Einstein AH structure has the following properties: for Weyl structures it specializes to the usual Einstein Weyl equations; the conjugate of an Einstein AH structure is again Einstein AH; and the AH structures induced on a hypersurface in affine space are Einstein if and only if the hypersurface is an affine hypersphere. 

There seem to be two principal reasons for interest in these equations. On the one hand, in section $7$ of \cite{Fox-ahs} or \cite{Fox-calg} there are given examples of Einstein AH structures in dimension $4$ and higher which are neither Weyl nor locally immersable as non-degenerate hypersurfaces in flat affine space, though otherwise as nice as possible (in the terminology of \cite{Fox-ahs}, they are exact, with self-conjugate curvature). This means the Einstein AH equations are a genuine generalization of the Einstein Weyl and affine hypersphere equations. On the other hand, via the theorem of Cheng-Yau associating an affine hypersphere to the universal cover of a manifold with strictly convex flat projective structure, methods for solving the Einstein AH equations should lead to analytic methods for producing convex flat projective structures. The two-dimensional case studied here illustrates that class of AH structures is rich yet sufficiently limited as to be amenable to characterization. 

The primary purpose of the present paper is to describe the classification up to the action of the group of orientation-preserving diffeomorphisms isotopic to the identity of the Riemannian signature Einstein AH structures on a compact orientable surface. It turns out that these all are either Einstein Weyl structures or have underlying projective structure which is flat and convex. Such a structure is determined in the former case by a conformal structure and a holomorphic vector field the real part of which is Killing for some metric in the conformal class, and in the latter case by a conformal structure and cubic holomorphic differential, and, consequently, by results going back to C.~P. Wang's \cite{Wang} and E. Calabi's \cite{Calabi-affinemaximal} and completed by F. Labourie's \cite{Labourie-flatprojective} and J. Loftin's \cite{Loftin-affinespheres, Loftin-riemannian, Loftin-compactification}, by a convex flat real projective structure. In the latter case it can be viewed as a solution of the Abelian vortex equations. Precise statements are given later in the introduction. 

For either Einstein Weyl structure or convex flat projective structures on surfaces the relevant classifications were already understood. The description of the deformation space of strictly convex flat projective structures is due independently to F. Labourie and J. Loftin, while the description of Einstein-Weyl structures on surfaces is mostly completed in the papers \cite{Calderbank-mobius, Calderbank-twod} of Calderbank. They are recounted here in part to show concretely how they fit into the formalism used here, and in part to highlight the relation with the Abelian vortex equations, which appears not to have been noted previously. The main novelty is the point of view, that there is a common framework encompassing both kinds of structures. The picture that emerges, and is suggestive of what might be true in higher dimensions, is roughly that for Einstein AH structures, positive curvature implies Weyl while negative curvature implies flatness and convexity of the underlying projective structure. The situation recalls that for extremal K\"ahler metrics. If there are no holomorphic vector fields then the scalar curvature of an extremal K\"ahler metric must be constant; the analogous statement here is that an Einstein AH structure on a compact Riemann surface admitting no holomorphic vector fields must be exact, and, as a consequence, have underlying flat projective structure which is strictly convex. On the other hand, if there are holomorphic vector fields, then there can be extremal K\"ahler metrics with nonconstant scalar curvature, as was shown by E. Calabi in \cite{Calabi-extremal}; the analogous examples here are just the Einstein Weyl structures on spheres and tori.

The final sections show that Einstein AH structures arise naturally in at least two other contexts, namely in the construction of Hessian metrics, and on mean curvature zero Lagrangian submanifolds of para-K\"ahler manifolds. In the penultimate section it is shown that the trivial real line bundle over a surface admitting an Einstein AH structure with parallel negative scalar curvature admits several interesting Hessian metrics of Riemannian and Lorentzian signatures which satisfy various Einstein type conditions. In the final section it is shown that an Einstein AH structure is induced on a mean curvature zero Lagrangian submanifold of a para-K\"ahler space form.

The remainder of the introduction describes the contents in more detail.

\subsection{}
Section \ref{backgroundsection} describes background needed in the remainder of the paper. The reader is advised to read it enough to be familiar with the notational and terminological conventions employed throughout. Sections \ref{ahsection} and \ref{ahcurvaturesection} describes the basic properties and local curvature invariants of AH structures on surfaces. As for a Riemann surface, the geometric structures considered admit equivalent descriptions in one-dimensional holomorphic terms or two-dimensional smooth real terms. The complex description generally leads to a more efficient and more transparent description, while the real description is more natural for comparing the two-dimensional results to the higher-dimensional case. Here both, and the transition from one to the other, are described, relying on material recounted in section \ref{holomorphicdifferentialsection} relating holomorphic differentials with conformal Killing and Codazzi tensors. 

\subsection{}
In section \ref{einsteinsection} the Einstein AH equations are defined and their most basic properties are noted. The defining conditions are given several reformulations. In particular, Lemma \ref{einsteinhololemma} shows that a Riemannian AH structure on an oriented surface is Einstein if and only if its Ricci curvature has type $(1,1)$ and its complex weighted scalar curvature (defined in section \ref{fdhsection}) is holomorphic.

\subsection{}
In order to state the main result of section \ref{scalarsection} it is necessary to recall some definitions. An open subset of the projective sphere is \textbf{convex} if its intersection with any projective line is connected. It is \textbf{properly convex} if its closure contains no pair of antipodal points. A properly convex domain is \textbf{strictly convex} if its boundary contains no open segment of a projective line. For example the positive orthant in $\rea^{n}$ (which is projectively equivalent to a standard simplex) is properly convex but not strictly convex. A flat real projective structure is \textbf{(strictly) convex} if its developing map is a diffeomorphism onto a (strictly) properly convex subset of the projective sphere.

Perhaps the principal technical result in the paper is Theorem \ref{scalarexacttheorem} which describes what are the possible types of Einstein AH structure on a compact orientable surface in terms of the genus of the surface, the sign of the weighted scalar curvature, the exactness or not of the AH structure, and whether the AH structure is Weyl or not. The precise statement is long, so is not repeated here, but roughly the result is that a Riemannian Einstein AH structure $(\en, [h])$, with conjugate $(\ben, [h])$, on a compact orientable surface $M$ of genus $g$ satisfies one of the following mutually exclusive possibilities:
\begin{enumerate}
\item $(\en, [h])$ is exact and Weyl, or, what is the same, $\en$ is generated by the Levi-Civita connection of a constant curvature metric representative of $[h]$.
\item The genus is $g \geq 2$. $(\en, [h])$ is exact and not Weyl with parallel negative scalar curvature, $\en$ and $\ben$ are strictly convex flat real projective structures, the cubic torsion is the real part of a holomorphic cubic differential, and a distinguished metric has negative scalar curvature.
\item $M$ is a torus. $(\en, [h])$ is exact and not Weyl with parallel negative scalar curvature, $\en$ and $\ben$ are flat real projective structures which are convex but not strictly convex, the cubic torsion is the real part of a holomorphic cubic differential, and a distinguished metric is flat.
\item $M$ is a torus, $(\en, [g])$ is Weyl and closed but not exact, the scalar curvature is zero, and the $(1,0)$ part of the aligned representative $\nabla \in \en$ is a holomorphic affine connection.
\item $M$ is a torus, $(\en, [h])$ is Weyl and not closed, and the scalar curvature changes signs. 
\item $M$ is a sphere, $(\en, [h])$ is Weyl and not closed, and the scalar curvature is somewhere positive.
\end{enumerate}
In particular an Einstein AH structure on $M$ which is not simply that generated by a conformal structure is either Weyl or exact, but not both. The proof of Theorem \ref{scalarexacttheorem} is based on Theorem \ref{classtheorem} which shows that for a distinguished metric representing the conformal structure the Einstein equations have a technically convenient form. Although the existence of this distinguished metric is a consequence of the Hodge decomposition, it is said to be \textbf{Gauduchon} because in the higher-dimensional setting the corresponding construction in the context of Einstein Weyl structures is due to P. Gauduchon; see \cite{Gauduchon, Gauduchon-circlebundles}.

Theorem \ref{magnetictheorem} shows that for an Einstein Weyl structure on a sphere or torus the integral curves of the metric dual of the Faraday primitive of a Gauduchon metric $h$ are magnetic geodesics for the magnetic flow generated by $h$ and a scalar multiple of the Faraday two-form. This interpretation of Einstein Weyl structures seems to be new.

\subsection{}
In section \ref{vortexsection} it is shown that an exact Einstein AH structure on a compact orientable surface is equivalent to a special sort of solution of the Abelian vortex equations. Recall that a triple $(\nabla, g, s)$ comprising a Hermitian structure $g$ on a complex line bundle $\emf$, a Hermitian connection $\nabla$ on $\emf$, and a smooth section $s$ of $\emf$ solves the abelian vortex equations if $\nabla$ induces a holomorphic structure on $\emf$ with respect to which $s$ is holomorphic, and there is satisfied a third equation relating the curvature of $\nabla$ and the Hermitian norm of $s$ (see \eqref{vortex}). This extension to compact K\"ahler manifolds of equations proposed by Landau and Ginzburg to model superconductivity was first studied by M. Noguchi in \cite{Noguchi} and S. Bradlow in \cite{Bradlow}. A special class of solutions, the \textbf{$p$-canonical} solutions, is distinguished by the requirement that the resulting $s$ be a section of the $p$th power of the canonical bundle with respect to the underlying complex structure, and that $\nabla$ be the Hermitian connection induced by the underlying K\"ahler structure. The Abelian vortices corresponding to exact Einstein AH structures are restricted in the sense that they correspond to $3$-canonical solutions. This gives a geometric interpretation to $3$-canonical solutions of the Abelian vortex equations which appears not to have been previously observed. It should be interesting to understand exactly how this plays out at the level of moduli spaces, though note that diffeomorphism inequivalent Einstein AH structures can determine gauge equivalent Abelian vortices, so at this level the correspondence is neither injective nor surjective. See the final paragraphs of section \ref{vortexsection} for a brief discussion.

\subsection{}
Section \ref{constructionsection} is devoted mainly to showing constructively that all of the possibilities identified in Theorem \ref{scalarexacttheorem} are realized and to describing the moduli/deformation spaces of solutions.

Combining Theorems \ref{scalarexacttheorem} and \ref{2dmodulitheorem} leads to the following theorem, the proof of which is completed in section \ref{convexsection}.

\begin{theorem}\label{summarytheorem}
On a compact orientable surface $M$ of genus at least $2$, the following spaces are in canonical bijection:
\begin{list}{(\arabic{enumi}).}{\usecounter{enumi}}
\renewcommand{\theenumi}{Step \arabic{enumi}}
\renewcommand{\theenumi}{(\arabic{enumi})}
\renewcommand{\labelenumi}{\textbf{Level \theenumi}.-}
\item \label{ds1} The fiber bundle over Teichm\"uller space of $M$ the fibers of which comprise the cubic holomorphic differentials.
\item \label{ds2} The deformation space of convex flat real projective structures.
\item \label{ds3} The deformation space of Einstein AH structures. 
\end{list}
The same equivalences are true for $M$ of genus $1$ provided that \ref{ds3} is replaced by
\begin{list}{(\arabic{enumi}).}{\usecounter{enumi}}
\setcounter{enumi}{3}
\item \label{ds3b} The deformation space of exact Einstein AH structures.
\end{list}
\end{theorem}
\noindent
Aside from its interest as such, Theorem \ref{summarytheorem} is suggestive of what can be expected to be true about the corresponding structures in higher dimensions. It plays for Einstein AH structures a role something like the uniformization theorem plays in the theory of higher-dimensional Einstein metrics.

The main point of section \ref{constructionsection} is to prove directly the implication \ref{ds1}$\Rightarrow$\ref{ds3} of Theorem \ref{summarytheorem}, which is Theorem \ref{2dmodulitheorem}. The equivalence \ref{ds1}$\iff$\ref{ds2} (and also, implicitly, the implication \ref{ds2}$\Rightarrow$\ref{ds3}) of Theorem \ref{summarytheorem} was known already, having been proved independently by F. Labourie (see \cite{Labourie-flatprojective}) and J. Loftin (see \cite{Loftin-affinespheres}). The proof given here of Theorem \ref{2dmodulitheorem} is not in essential points different from Loftin's proof of the implication \ref{ds1}$\Rightarrow$\ref{ds2}, but is set in a slightly more general context so as to yield Corollary \ref{vortexcorollary}, which shows how to construct certain $p$-canonical solutions of the Abelian vortex equations. Also there are extracted some bounds on the solutions which yield bounds on the volume and curvature of distinguished metrics of Einstein AH structures. These represent very preliminary steps in the direction of understanding the analogues in the present context of Teichm\"uller curves. As is explained in section \ref{singularmetricsection} there is a natural diffeomorphism equivariant action of $GL^{+}(2, \rea)$ on the space of cubic holomorphic differentials coming from its action on the singular flat Euclidean structure determined by such a differential on the complement of its zero locus. Theorem \ref{summarytheorem} means that the orbits of this action on the space of cubic holomorphic differentials determine disks in the deformation space of Einstein AH structures. An analysis of the structure of these disks awaits, but in section \ref{constructionsection} there are proved some results about the path in the deformation space corresponding to a ray in the space of cubic holomorphic differentials. Theorem \ref{liptheorem} shows that with a suitable parameterization of the ray the conformal factor relating the Gauduchon metric at time $t$ to the Gauduchon metric at time $0$ is pointwise non-decreasing and Lipschitz in $t$. This makes possible some statements about the limiting behavior along the ray of the volume and curvature of suitably scaled Gauduchon metrics.

In order to say a bit more about the background to Theorem \ref{summarytheorem}, some context is recalled. In \cite{Choi-Goldman}, S. Choi and W.~M. Goldman showed that the deformation space of convex flat real projective structures on a two-dimensional orbifold $M$ of negative orbifold characteristic $\chi(M)$ is homeomorphic to a cell of a dimension equal to $-8\chi(M) + k$, where $k$ is a number expressible in terms of the orders of the stabilizers of the singular points; this generalizes an earlier theorem of Goldman in \cite{Goldman-convex} for the manifold case. It follows from a theorem of Thurston that a compact $2$-orbifold of non-positive orbifold Euler characteristic is a quotient of a manifold by a finite group; see the end of section $1.2$ of \cite{Choi-Goldman}. In \cite{Wang}, C.~P. Wang shows how a hyperbolic affine hypersphere gives rise to a conformal structure and a cubic holomorphic differential, and conversely, how given such data on a compact oriented surface there is associated to its universal cover a hyperbolic affine hypersphere. On the affine hypersphere over the universal cover of $M$ the difference tensor of the affine connection induced via the affine normal and the Levi-Civita connection of an equiaffine metric is the real part of a holomorphic cubic differential (it is often called the \textit{Pick form}). This observation underlies the Labourie-Loftin theorem.

The content of the implication \ref{ds2}$\Rightarrow$\ref{ds3} of Theorem \ref{summarytheorem} is the claim that the Einstein AH structure is determined by its underlying projective structure (which is necessarily flat by Lemma \ref{2deinsteinlemma}). As is described briefly now, and in more detail in section \ref{convexsection}, this can be deduced from various results of S.~Y. Cheng and S.~T. Yau. By a theorem of Cheng and Yau (see \cite{Loftin-survey} for the full history) resolving a conjecture of E. Calabi, the interior of the cone over a properly convex domain admits a unique foliation by hyperbolic affine hyperspheres asymptotic to the cone. In particular, this can be applied to the cone over the universal cover of a $2$-orbifold $M$ carrying a convex flat real projective structure, and in the manifold case the AH structure induced on the affine hypersphere descends to $M$. Since these AH structures are always exact, this means $M$ carries a canonical homothety class of metrics (those induced by the equiaffine metrics on the affine hyperspheres foliating the interior of the cone) and so a distinguished connection, the Levi-Civita connection of any one of these metrics. Technically, this has two aspects. One is that an immersed hyperbolic affine hypersphere is properly embedded if and only if the induced affine metric is complete; for references to a proof see section $5$ of \cite{Trudinger-Wang-survey}. The other is that a convex flat real projective structure determines a distinguished conformal structure. The latter can be obtained from either of two theorems of Cheng and Yau solving certain Monge-Amp\`ere equations. In \cite{Cheng-Yau-mongeampere}, it is shown that on a bounded convex domain $\Omega \subset \rea^{n}$ there is a unique smooth, convex solution of the equation $u^{n+2}\det \hess u = (-1)^{n}$ vanishing on the boundary of $\Omega$. The radial graph of $u$ is the desired affine hypersphere (see Theorem $3$ of \cite{Loftin-riemannian}). Alternatively, in \cite{Cheng-Yau-realmongeampere} it is shown that on a convex cone $\Omega \subset \rea^{n+1}$ containing no complete affine line there is a unique smooth function $F$ solving $\det \hess F = e^{2F}$, tending to $+\infty$ at the boundary of the cone, and such that $\hess F$ is a complete Riemannian metric on the interior of the cone; the level sets of $F$ are the desired affine hyperspheres asymptotic to the cone (the uniqueness claim follows by passing to the corresponding K\"ahler metric on the tube domain over $\Omega$ and appealing to the Schwarz lemma for volume forms in \cite{Mok-Yau}). Because this approach does not appear to have been written anywhere, it is sketched in section \ref{hessianmetricsection}. These theorems of Cheng-Yau should be viewed as real analogues of the theorems of Cheng, Mok and Yau producing complete K\"ahler Einstein metrics, e.g. \cite{Cheng-Yau-completekahler}, \cite{Mok-Yau}. (In fact the latter theorem follows from the specialization to a tube domain over a pointed convex cone of the theorem of N. Mok and Yau in \cite{Mok-Yau} producing a K\"ahler Einstein metric on a bounded domain of holomorphy). For further background the reader can consult, in addition to papers already mentioned, \cite{Loftin-survey, Loftin-compactification}, \cite{Cheng-Yau-mongeampere, Cheng-Yau-affinehyperspheresI}, and \cite{Calabi-improper, Calabi-bernsteinproblems, Calabi-completeaffine, Calabi-affinemaximal}. In \cite{Labourie-flatprojective}, Labourie has shown that in two dimensions these results admit more direct and simpler proofs, and has given a variety of other ways of understanding them.

\subsection{}
The remaining Einstein AH structures on compact orientable surfaces are all Weyl and occur on either the torus or the sphere. The existence of such structures in the case of zero scalar curvature (on the torus) is straightforward, and they correspond in a natural way to holomorphic affine connections. This is explained in section \ref{torusconstantsolutionssection}. The remaining cases are Einstein Weyl structures which are not closed. Such a structure determines a vector field $X$ which is Killing for the Gauduchon metric. It follows that $X$ is the real part of a holomorphic vector field, though it is not an arbitrary holomorphic vector field, for $X$ preserves some metric and so generates a circle action; in particular on the torus it must generate a rational flow. Using this circle action the equations that need to be solved to construct an Einstein Weyl structure are reduced to an ODE, which after further reductions admits explicit solutions in terms of elementary functions. In the case of the sphere, the moduli space of Einstein Weyl structures is essentially parameterized by conjugacy classes of elliptic elements of $PSL(2, \com)$, that is by a half-open interval. The precise statement is Theorem \ref{spheremodulitheorem}. In the torus case, because of the aforementioned rationality condition, it is not clear how to describe nicely the deformation space.

Much of the description of Einstein-Weyl structures given in section \ref{spheretorussection} can be found in some similar form in \cite{Calderbank-mobius} and \cite{Calderbank-twod}, although the presentation here is perhaps more elementary, and is made to illustrate the realization of the possibilities stated in Theorem \ref{scalarexacttheorem}; also the description of when the solutions found are equivalent is made more explicitly than it is in \cite{Calderbank-mobius} and \cite{Calderbank-twod}, although it seems that the results were understood by Calderbank.

\subsection{}\label{maintrosection}
In section \ref{hessianmetricsection} it is shown that an Einstein AH structure $(\en, [h])$ on a compact orientable surface $M$ of genus $g > 1$ gives rise to two metrics $f_{IJ}$ and $g_{IJ}$ on $M \times \reat$ which together with the flat affine structure induced on $M \times \reat$ by $\en$ constitute Hessian metrics. The metric $g_{IJ}$ is a Lorentz signature Monge-Amp\`ere metric, while $f_{IJ}$ is a Riemannian signature Hessian metric which is Einstein when viewed as a K\"ahler affine metric in the sense of \cite{Cheng-Yau-realmongeampere}. Precise statements are given in Theorem \ref{hessiantheorem}. In particular, it is shown that, with respect to $g_{IJ}$, $M$ is a smoothly immersed, spacelike, umbilic hypersurface of constant mean curvature. In section \ref{dustsection} is explained that the metric $g_{IJ}$ can also be viewed as a solution of $2+1$ gravitational equations with stress energy tensor of the form corresponding to a pressureless perfect fluid.

In section \ref{mongeamperemetricsection} it is shown that for each $C > 0$ there is a Riemannian signature Monge-Amp\`ere metric on $M \times [-\log C, \infty)$. Its potential has the form $\Psi(F)$ where $F$ is the potential for the metric $f_{IJ}$ described in the previous paragraph, and $\Psi$ is the function given in \eqref{psit}. The precise statement is Theorem \ref{mongeamperetheorem}.

\subsection{}
In \cite{Aledo-Espinar-Galvez}, J.~A. Aledo, J.~M. Espinar, and J.~A. G\'alvez have studied a surface equipped with what they call a \textbf{Codazzi pair}, which is a pair comprising a Riemannian metric and a second symmetric covariant tensor satisfying the Codazzi equations with respect to the Levi-Civita connection of the metric. They view this as an abstraction of the geometric structure induced on a hypersurface in Euclidean three space, and have shown that many classical results can be strengthened in this context. Though the setting is different from that considered here, the perspective is similar. One of the motivations of \cite{Aledo-Espinar-Galvez} is that such Codazzi pairs arise in other ways, e.g. the real part of the Hopf differential of a harmonic mapping yields such a pair. Similarly AH structures on surfaces arise naturally in the study of submanifolds of para-K\"ahler manifolds. In section \ref{parakahlersection} it is explained that a mean curvature zero spacelike Lagrangian immersed submanifold of a para-K\"ahler manifold of constant para-holomorphic sectional curvature inherits an Einstein AH structure; see Theorem \ref{parakahlertheorem} and the final paragraph of section \ref{pksection}. It is indicated also how to associate to certain Einstein AH structures such an immersion. I thank an anonymous referee for pointing out that essentially equivalent ideas have been worked out independently in R. Hildebrand's papers \cite{Hildebrand-crossratio} and \cite{Hildebrand-parakahler}. In particular Theorem \ref{parakahlertheorem} is equivalent to theorems in \cite{Hildebrand-crossratio} modulo changes of terminology. Related constructions appeared already in L. Vrancken's \cite{Vrancken-centroaffine}. Corollary \ref{parakahlercorollary} states some restrictions on mean curvature zero Lagrangian immersions in four dimensional para-K\"ahler space forms resulting from applying Theorem \ref{scalarexacttheorem} to the induced Einstein AH structure. Closely related results about minimal Lagrangian immersions in complex hyperbolic space are the focus of the preprint \cite{Loftin-Mcintosh} of J. Loftin and I. McIntosh.

\section{Notation and terminology}\label{backgroundsection}
This section records the notational and terminological conventions in use throughout the paper.

\subsection{}
Throughout $M$ is a connected smooth ($\cinf$) manifold. For a vector bundle $E$, $\Ga(E)$ denotes the space of its smooth sections (even if $E$ has a holomorphic structure). Its $k$th symmetric power and top exterior power are written respectively as $S^{k}(E)$ and $\Det E$. For a line bundle $E$, $|E|$ is the tensor product of $E$ with its orientation bundle. Sections of $|\Det \ctm|^{\la}$ are called \textbf{$\la$-densities}. A tensor taking values in $|\Det \ctm|^{\la}$ for some $\la \in \reat$ is said to be \textbf{weighted}. If $M$ is compact the integral of a $1$-density has sense, and so there is a bilinear pairing $\lb\dum, \dum \ra:\Ga(|\Det \ctm|^{\la}) \times \Ga(|\Det \ctm|^{1-\la}) \to \rea$ between densities of weight $\la$ and the complementary weight $(1-\la)$. 

\subsection{}\label{holobackgroundsection}
Given a complex structure $J \in \eno(\ste)$ on the real vector space $\ste$, the complexification $\stec$ decomposes as the direct sum $\stec = \ste^{1,0} \oplus \ste^{0, 1}$ of the $\pm \j$ eigensubbundles of the extension of $J$ to $\stec$ by complex linearity. The induced action of $J$ on $\sted$ is defined by $J(\al) \defeq  \al \circ J$, so that $(\sted, J)$ is an a complex vector space, and the $(1,0)$ part $\ste^{\ast\, 1,0}$ of $\sted$ is the $\com$-dual of $\ste^{1, 0}$ and annihilates $\ste^{0,1}$. A completely symmetric or anti-symmetric tensor decomposes by type. If $B$ is in $S^{p+q}(\sted)$ or $S^{p+q}(\ste)$ denote by $B^{(p, q)}$ the $(p,q)$ part of its complex linear extension as an element of $S^{p+q}(\sted\tensor_{\rea}\com)$. For example, for $\al \in \sted$, $2\al^{(1, 0)} = \al - \j J(\al) = \al - \j \al \circ J$.

The preceeding makes sense fiberwise on a complex vector bundle. If $(M, J)$ is a complex manifold, there is written $TM \tensor_{\rea}\com = T^{1, 0} \oplus T^{0, 1}$, while the complex tangent bundle $\tcm$ is $TM$ viewed as a complex vector bundle; it is identified as a complex vector bundle with $T^{1,0}$. A \textbf{holomorphic structure} on a complex vector bundle $E$ over a complex manifold is a linear differential operator $\hat{D}$ sending $E$-valued $(p,q)$-forms to $E$-valued $(p, q+1)$-forms, satisfying the Leibniz rule $\hat{D}(fs) = \delbar f \wedge s + f\hat{D}s$ for any $f \in \cinf(M)$ and any smooth section $s$ of $E$, and such that $\hD^{2} = 0$. A smooth local section of $E$ is \textbf{$\hat{D}$-holomorphic} if it is in $\ker \hat{D}$. If $(M, J)$ is a complex manifold then any bundle of complex tensors (a tensor product of powers of $\tcm$ and its dual) over $M$ is naturally a holomorphic vector bundle, and the corresponding holomorphic structure is denoted $\delbar$. The exterior differential decomposes by type as $d = \del + \delbar$. 

For a one-complex dimensional manifold and $p \in \integer$, define $\cano^{p}$ to be the $p$th power of the complex cotangent bundle $\tcm^{\ast}$ viewed as a holomorphic line bundle, and viewed also as the $(-p)$th symmetric power of $T^{1,0}$ (a negative power means the power of the dual).  A smooth (holomorphic) section of $\cano^{p}$ is called a complex (holomorphic) \textbf{$p$-differential}. In a local holomorphic coordinate $z$, a holomorphic section $p$-differential $\si$ has the form $\phi(z)dz^{p}$ for a holomorphic funtion $\phi(z)$.

\subsection{}
Tensors are usually indicated using the abstract index notation. Ordinary tensors are indicated by lowercase Latin abstract indices, so that, for instance, $a_{ij}$ indicates a covariant two tensor. If a complex structure is given, lowercase Greek indices, e.g. $\al, \be, \ga$, etc., decorate sections of the tensor powers of $T^{1,0}$ and its complex dual, while barred lowercase Greek indices, e.g. $\bar{\al}, \bar{\be}, \bar{\ga}$, etc., decorate sections of the tensor powers of $T^{0, 1}$ and its complex dual. Enclosure of indices in square brackets (resp. parentheses) indicates complete skew-symmetrization (resp. complete symmetrization), so that for example $a^{ij}= a^{(ij)} + a^{[ij]}$ indicates the decomposition of a contravariant two-tensor into its symmetric and skew-symmetric parts, and $(X \wedge Y)^{ij} = 2X^{[i}Y^{j]}$ for vector fields $X$ and $Y$. The summation convention is always in effect in the following form: indices are in either \textit{up} position or \textit{down} and a label appearing as both an up index and a down index indicates the trace pairing. Since polynomials on the vector space $\ste$ are tautologically identified with symmetric tensors on the dual vector space $\ste^{\ast}$, the index $i$ in $\tfrac{\pr}{\pr y_{i}}$ has to be regarded as an \textit{up} index. The horizontal position of indices is always maintained when they are raised or lowered. The interior multiplication of a vector field $X^{i}$ in a covariant tensor $B_{i_{1}\dots i_{k}}$ is defined by $\imt(X)B_{i_{1}\dots i_{k-1}} \defeq X^{p}B_{pi_{1}\dots i_{k-1}}$.

\subsection{}
The \textbf{curvature} $R_{ijk}\,^{l}$ of a torsion-free affine connection $\nabla$ is defined by $2\nabla_{[i}\nabla_{j]}X^{k} = R_{ijp}\,^{k}X^{p}$. The \textbf{Ricci curvature} is the trace $R_{ij} \defeq R_{pij}\,^{p}$. 

\subsection{}\label{riemannsurfacehodgestarsection}
A non-degenerate weighted covariant two-tensor $h_{ij}$ determines a contravariant two-tensor $h^{ij}$ of complementary weight defined by $h^{ip}h_{jp} = \delta_{j}\,^{i}$, in which here, as always, $\delta_{i}\,^{j}$ is the tautological $\binom{1}{1}$-tensor determined by the pairing of vectors with covectors. By $\det h$ is meant the $2$-density which satisfies $\lb \det h, E_{1}\wedge \dots \wedge E_{n}\ra = \det h(E_{i}, E_{j})$ for any frame $E_{1}, \dots, E_{n}$.

A \textbf{pseudo-Riemannian metric} or, simply, a \textbf{metric} means a non-degenerate covariant two-tensor $h_{ij}$. The metric is \textbf{Riemannian} if it is positive definite. A \textbf{conformal structure} $[h]$ means a pseudo-Riemannian metric determined up to multiplication by a positive function. A conformal structure is identified with its \textbf{normalized representative} $H_{ij}\defeq |\det h|^{-1/\dim M}h_{ij}$ which takes values in the bundle of $-(2/\dim M)$-densities.

For conformal metrics $\tilde{h}_{ij} = fh_{ij}$ the Levi-Civita connections are written $\tD$ and $D$, and their difference tensor is written $\tD - D =   2\si_{(i}\delta_{j)}\,^{k} - h_{ij}h^{kp}\si_{p}$, for $2\si_{i} = d\log{f}_{i}$. The curvature of $D$ is written $\sR_{ijk}\,^{l}$. Objects corresponding to $\tD$ are written with the same notations as those corresponding to $D$, but decorated with a $\tilde{\,}$, although for the scalar curvatures it will be convenient to write $\sR_{\tilde{h}}$ and $\sR_{h}$ rather than $\tilde{\sR}$ and $\sR$. For example, the scalar curvature changes under conformal rescaling by 
\begin{align}
\label{conformalscalardiff}&f\sR_{\tilde{h}} = \sR_{h} - 2h^{pq}D_{p}\si_{q}  = \sR_{h} - \lap_{h}\log f .
\end{align}
Here $\lap_{h}$ is the rough Laplacian $h^{pq}D_{p}D_{q}$, which acts on tensors as well as on functions. 
 
Given a metric $h$, for any $X_{i_{1}\dots i_{k}}\,^{j_{1}\dots j_{l}}$ the notations $X^{\sharp}$ and $X^{\flat}$ indicate the tensors obtained by raising and lowering all indices using $h$. That is $X^{\flat}_{i_{1}\dots i_{k+l}} = X_{i_{1}\dots i_{k}}\,^{j_{1}\dots j_{l}}h_{j_{1}i_{k+1}}\dots h_{j_{l}i_{k+l}}$, and similarly for $X^{\sharp}$. The $h$-norm of a tensor $X$ is defined by complete contraction, e.g. for $A_{ij}\,^{k}$, $|A|_{h}^{2} = A^{\flat}_{ijk}A^{\sharp\,ijk} = A_{ij}\,^{k}A_{ab}\,^{c}h^{ia}h^{jb}h_{kc}$.

\subsection{}\label{riemannsurfacesection}
A \textbf{Riemann surface} is an oriented surface $M$ equipped with a Riemannian signature conformal structure $[h]$. (Note that \textit{Riemannian surface} and \textit{Riemann surface} are not synonyms; the former indicates there is given a distinguished metric, while the underlying smooth structure need not be orientable). On a Riemann surface there is a unique almost complex structure $J$ defined in terms of any $h \in [h]$ by the requirements that $J_{i}\,^{p}J_{j}\,^{q}h_{pq} = h_{ij}$ and that the two-form $\om_{ij} \defeq J_{i}\,^{p}h_{pj}$ determine the given orientation. As these conditions do not depend on the choice of $h \in [h]$, $J$ is determined by the conformal structure $[h]$. On the other hand, given a complex structure $J$ on a surface $M$, the orientation determined by $X \wedge JX$ does not depend on the choice of $X$, and there is a unique Riemannian signature conformal structure $[h]$ such that $J_{i}\,^{p}J_{j}\,^{q}h_{pq} = h_{ij}$. Any almost complex structure on a surface is integrable, or, equivalently $DJ = 0$, where $D$ is the Levi-Civita connection of $h$, and so $(h, J, D)$ is a K\"ahler structure on $M$ with K\"ahler form $\om_{ij}$ and associated \textbf{Hermitian metric} $h^{(1, 1)}$. On any bundle $E$ of complex tensors, $h$ determines a Hermitian structure, and $D$ induces the unique Hermitian connection such that $D^{0, 1} = \delbar$.

\subsection{}
For a Hermitian metric $h^{(1,1)}_{\al\bar{\be}}$ it is convenient to omit the superscript $(1,1)$, and to write instead simply $h_{\al\bar{\be}}$, although it should be kept in mind that $h_{ab}$ and $h_{\al\bar{\be}}$ refer to different objects. The dual bivector $h^{\al\bar{\be}}$ is defined by $h^{\al\bar{\be}}h_{\ga\bar{\be}} = \delta_{\al}\,^{\ga}$. The conventions are such that $\delta_{\al}\,^{\be}$ indicates the tautological endomorphism of $T^{1,0}$, while $\delta_{i}\,^{j}$ indicates the tautological endomorphism of $TM$. The convention is that for a section $X^{a_{1}\dots a_{k}}$ of $\tensor^{k}TM$, the Hermitian norm of its $(k, 0)$ part $X^{\al_{1}\dots \al_{k}}$ is defined by complete contraction with its complex conjugate using $h_{\al\bar{\be}}$; e.g. for $X \in \Ga(TM)$, $|X^{(1,0)}|^{2}_{h} = X^{\al}X^{\bar{\be}}h_{\al\bar{\be}} = \tfrac{1}{2}|X|_{h}^{2}$. 

\subsection{}
The formal adjoint $\dad_{h}$ of the exterior differential with respect to the pairing of forms given by $h$ is given on $k$-forms by $\dad_{h}\al_{i_{1}\dots i_{k-1}} = -D^{p}\al_{p i_{1}\dots i_{k-1}}$. Since, for a $k$-form $\al_{i_{1}\dots i_{k}}$ there holds $f\dad_{\tilde{h}}\al_{i_{1}\dots i_{k-1}} = \dad_{h}\al_{i_{1}\dots i_{k-1}} + 2(k-1)\si^{\sharp\, p}\al_{pi_{1}\dots i_{k-1}}$, that a one-form be co-closed (in $\ker \dad_{h}$) is a conformally invariant condition. For a one-form whether there is written $D^{p}\ga_{p}$ or $-\dad_{h}\ga$ will depend on context. On a Riemann surface the action of the Hodge star operator $\hs = \hsh$ on one-forms depends only on the conformal class $[h]$, and is given in terms of the complex structure $J$ determined by $[h]$ by $(\hs \al)_{i} = -\al_{p}J_{i}\,^{p}$. On an oriented surface $\dad_{h} = -\star d \star$. 

\section{Holomorphic differentials and conformal Killing and Codazzi tensors}\label{holomorphicdifferentialsection}
This section \ref{holomorphicdifferentialsection} records some basic facts about holomorphic differentials on Riemann surfaces. The purely real equations for symmetric tensors characterizing conformal Killing and Codazzi tensors make sense in higher dimensions (see section $6$ of \cite{Fox-ahs}), but on surfaces their complex counterparts are easier to work with. Lemma \ref{kdifferentialslemma} states the identification on Riemann surfaces of Codazzi and conformal Killing tensors as the real parts of holomorphic differentials, and Lemma \ref{flatmetriclemma} shows that such a differential determines a singular flat metric. For quadratic differentials these statements are well known and widely utilized. Although the general statements are surely also well known, there does not seem to be any convenient reference, and so it seems useful to include them. The correspondence between holomorphic differentials and singular flat metrics is discussed in more detail in sections \ref{khodgesection} and \ref{singularmetricsection}, where it is explained how it yields a diffeomorphism equivariant action of $GL^{+}(2, \rea)$ on the space of such differentials. While not much use is later made of this material, it motivates some estimates proved in section \ref{constructionsection} and when coupled with Theorem \ref{summarytheorem} and compared with usual Teichm\"uller theory, suggests many questions.

\subsection{}
To a Young diagram the boxes of which are labeled with distinct indices corresponds the irreducible $GL(n, \rea)$ module comprising tensors skew-symmetric in the indices in a given column of the Young diagram and vanishing when skew-symmetrized over the indices in a given column and any index in any box to the right of the given column. The irreducible representations of the subgroup $CO(h)$ of $GL(n, \rea)$ acting conformally with respect to a fixed metric $h$ on $\rea^{n}$ are described in \cite{Weyl}. The subspace of an irreducible $GL(n, \rea)$ representation comprising tensors completely trace-free with respect to $h$ is a representation of $CO(h)$. Lemma \ref{weylcriterion} will be invoked repeatedly.

\begin{lemma}[\cite{Weyl}, Theorem $5.7.$A]\label{weylcriterion}
The $CO(h)$-module of covariant trace-free tensors on $\rea^{n}$ having symmetries corresponding to a Young diagram is trivial if the sum of the lengths of the first two columns of the Young diagram is greater than $n$. 
\end{lemma}

For instance, Lemma \ref{weylcriterion} implies that the usual conformal Weyl tensor of a Riemannian metric vanishes identically on a manifold of dimension at most $3$.

If there is given a fiberwise metric on the vector bundle $E$ then $S^{k}_{0}(E)$ denotes the subbundle of $S^{k}(E)$ comprising elements trace-free with respect to the given metric. The convention is that $S^{1}_{0}(E) = E$; this corresponds to regarding the trace as the zero map on vectors. 
\begin{lemma}\label{twodablemma}
Let $h_{ij}$ be a constant metric on $\rea^{2}$, and for $k \geq 1$ let $A, B \in S^{k}_{0}(\rea^{2})$. Then 
\begin{align}\label{twodab}
2A_{a_{1}\dots a_{k}(i}B_{j)}\,^{a_{1}\dots a_{k}} = A_{a_{1}\dots a_{k}}B^{a_{1}\dots a_{k}}h_{ij}, 
\end{align}
in which indices are raised and lowered with $h_{ij}$ and its inverse $h^{ij}$. Let $X^{i}$ be a vector. Then $2|\imt(X)B|^{2}_{h} = |X|_{h}^{2}|B|_{h}^{2}$. In particular if $h_{ij}$ has definite signature the equations $|B|_{h}^{2}|X|_{h}^{2} = 0$ and $X^{p}B_{pi_{1}\dots i_{k-1}} = 0$ are equivalent.
\end{lemma}
\begin{proof}
The tensor
\begin{align*}
\om_{ija_{1}\dots a_{k+1}} = h_{i(a_{1}}A_{a_{2}\dots a_{k})j} + h_{j(a_{1}}A_{a_{2}\dots a_{k})i} - h_{ij}A_{a_{1}\dots a_{k}} - h_{(a_{1}a_{2}}A_{a_{3}\dots a_{k})ij}
\end{align*}
is completely trace-free. As $\om_{ija_{1}\dots a_{k}} = \om_{(ij)(a_{i}\dots a_{k})}$ and $\om_{i(ja_{1}\dots a_{k})} = 0$, it follows from Lemma \ref{weylcriterion} that $\om = 0$. Hence $2A_{a_{1}\dots a_{k}(i}B_{j)}\,^{a_{1}\dots a_{k}} - A_{a_{1}\dots a_{k}}B^{a_{1}\dots a_{k}}h_{ij} = B^{a_{1}\dots a_{k}}\om_{ija_{1}\dots a_{k}} = 0$. By \eqref{twodab} there holds $2|\imt(X)B|^{2}_{h} = 2X^{p}B_{pi_{1}\dots i_{k-1}}X^{q}B_{q}\,^{i_{1}\dots i_{k-1}} = |X|_{h}^{2}|B|_{h}^{2}$.
\end{proof}

\subsection{}
On a Riemannian manifold $(M, h, D)$, a symmetric tensor $\si \in \Ga(S^{k}(\ctm))$ is \textbf{Codazzi} if $D\si \in \Ga(S^{k+1}(\ctm))$. From
\begin{align}\label{dsym2}
&(k+1)\left(D_{i}\si_{i_{1}\dots i_{k}} - D_{(i}\si_{i_{1}\dots i_{k})}\right) =  2\sum_{s = 1}^{k}D_{[i}\si_{i_{s}]i_{1}\dots \hat{i}_{s}\dots i_{k}},
\end{align}
(in which a $\hat{\,}$ denotes the omission of the index), it is evident that $\si$ is Codazzi if and only if $D_{[i_{1}}\si_{i_{2}]\dots i_{k+1}} = 0$. Write $\div_{h}$ (or simply $\div$) for the divergence operator $\div_{h}(\si)_{i_{1}\dots i_{k}} \defeq D^{p}\si_{p i_{1}\dots i_{k}}$ determined by $h$. For any tensor $\si$ let $\tf_{h}(\si)$ be its trace-free part. Let $\clie_{h}:\Ga(\symk) \to \Ga(\symkp)$ be the formal adjoint of the composition $\div_{h} \circ \tf_{h}$ with respect to the pairing of sections of $\symkp$ determined by integration. Explicitly, for $\si \in \Ga(\symkt)$ and $M$ a surface,
\begin{align}\label{cliedefined}
\clie_{h}(\si)_{i_{1}\dots i_{k+1}} = D_{(i_{1}}\si_{i_{2}\dots i_{k+1})} - \tfrac{1}{2}h_{(i_{1}i_{2}}\div_{h}(\si)_{i_{3}\dots i_{k+1})}.
\end{align}

\begin{lemma}\label{codazzilemma}
On a Riemannian surface $(M, h)$, for $\si \in \Ga(S^{k}_{0}(\ctm))$ there hold
\begin{align}\label{dsi}
 D_{i}\si_{i_{1}\dots i_{k}} &= \clie_{h}(\si)_{ii_{1}\dots i_{k}} + h_{i(i_{1}}\div_{h}(\si)_{i_{2}\dots i_{k})} - \tfrac{1}{2}h_{(i_{1}i_{2}}\div_{h}(\si)_{i_{3}\dots i_{k})i}, \qquad&&  \text{if}\,\, k > 1,\\
\label{dsi2}D_{i}\si_{j} & = \clie_{h}(\si)_{ij} + \tfrac{1}{2}d\si_{ij} + \tfrac{1}{2}\div_{h}(\si)h_{ij}, \qquad&& \text{if}\,\, k = 1.
\end{align} 
For $k > 1$, the following are equivalent for $\si \in \Ga(\symkt)$: 1. $\si$ is Codazzi; 2. $\si$ is divergence free; 3. $D\si = \clie_{h}(\si)$.
\end{lemma}
\begin{proof}
The tensor on the righthand side of \eqref{dsi} is trace-free and its complete symmetrization vanishes, so it lies in the irreducible $O(h)$-module corresponding to a Young diagram with $k-1$ boxes in its first row and $1$ box in its second row. By Lemma \ref{weylcriterion} this module is trivial if $k > 1$. When $k > 1$, for $\si \in \Ga(S^{k}_{0}(\ctm))$ it follows from \eqref{dsi} that $D\si \in \Ga(S^{k+1}_{0})(\ctm))$ if and only if $\div_{h}(\si) = 0$. This also follows by tracing the identity
\begin{align}\label{2dci}
\begin{split}
2D_{[a}\si_{b]i_{1}\dots i_{k-1}}& 
= h_{a(i_{1}}\div_{h}(\si)_{i_{2}\dots i_{k-1})b} -h_{b(i_{1}}\div_{h}(\si)_{i_{2}\dots i_{k-1})a},
\end{split}
\end{align}
which also follows from Lemma \ref{weylcriterion}. From \eqref{dsi} there follows $\div_{h}(\om) =0$ if and only if $D\si =\clie_{h}(\si)$.
\end{proof}

On a surface $f\div_{\tilde{h}}(\si) = \div_{h}(\si)$, so that the space $\ker \div_{h} \cap \Ga(\symkt)$ of trace-free Codazzi tensors depends only on the conformal class of $h$. As $\clie_{e^{\phi}h}(e^{k\phi}\si)_{i_{1}\dots i_{k+1}} = e^{k\phi}\clie_{h}(\si)_{i_{1}\dots i_{k+1}}$, the subspace $\Ga(\symkt)\cap \ker \clie$ is conformally invariant. The operator $\clie_{h}^{\sharp}:\Ga(\symktv) \to \Ga(\symkptv)$ defined by $\clie_{h}^{\sharp}(X) = \clie_{h}(X^{\flat})^{\sharp}$ satisfies $e^{\phi}\clie^{\sharp}_{e^{\phi}h}(X) = \clie^{\sharp}_{h}(X)$. For a vector field $X$ there holds $2\clie_{h}(X^{\flat})_{ij} = \tf(\lie_{X}h)_{ij}$, which motivates the notation resembling that for the Lie derivative. Moreover, this shows that $\Ga(TM) \cap \ker \clie_{h}^{\sharp}$ comprises the conformal Killing fields. For this reason the sections of $\ker \clie^{\sharp} \cap \Ga(\symktv)$ are called \textbf{conformal Killing tensors}. 

Lemma \ref{kdifferentialslemma} explains the conformal invariance on surfaces of trace-free Codazzi and conformal Killing tensors by showing that they are exactly the real parts of holomorphic differentials.

\subsection{}
An almost complex structure $J$ on a real vector space $\ste$ determines an almost complex structure $\bJ$ on $\tensor^{k}(\sted)$ defined by $\bJ(B)_{i_{1}\dots i_{k}} = J_{i_{1}}\,^{p}B_{pi_{2}\dots i_{k}}$.

\begin{lemma}\label{apolaritylemma}
Fix a two-dimensional real vector space $\ste$ with an almost complex structure $J$. Let $k > 1$. The map $B \to \bJ(B)$ is a complex structure on the $2$-dimensional real vector space $S^{k}_{0}(\sted)$.  For $B \in S^{k}_{0}(\sted)$ there holds $2B^{(k, 0)} = B - i\bJ(B)$, so that the $(1,0)$ part of $B$ \textnormal{qua} element of $(S^{k}_{0}(\sted), \bJ)$ equals the $(k, 0)$ part of $B$ relative to $J$. There results a complex linear isomorphism between the complexification $S^{k}_{0}(\sted) \tensor_{\rea}\com$ and $\cano^{k} \oplus \bar{\cano}^{k}$ such that the $\pm \j$ eigenspaces of $\bJ$ on $S^{k}_{0}(\sted)\tensor_{\rea}\com$ are identified respectively with $\cano^{k}$ and $\bar{\cano}^{k}$. 
\end{lemma}

\begin{proof}
Let $h$ be a definite signature metric on $\ste$ compatible with $J$. Let $B \in S^{k}(\sted)$. If $X \in \ste$ is non-zero, then $\{X, JX\}$ is an $h$-orthogonal basis of $\ste$. From the evidently equivalent identities
\begin{align}\label{jb0}
&\bJ(B)_{[i_{1}i_{2}]i_{3}\dots i_{k}} = J_{[i_{1}}\,^{p}B_{i_{2}]i_{3}\dots i_{k}p} = 0, &  &J_{i_{1}}\,^{p}J_{i_{2}}\,^{q}B_{pqi_{3}\dots i_{k}} = -B_{i_{1}\dots i_{k}}.&
\end{align}
it follows that $B \in S^{k}_{0}(\sted)$ if and only if $\bj(B)$ is completely symmetric. Since $B = -\bj(\bj(B))$, the same statement with the roles of $B$ and $\bj(B)$ interchanged shows that in this case $\bj(B) \in S^{k}_{0}(\sted)$. Thus $B \in S^{k}_{0}(\sted)$ if and only if $\bj(B) \in S^{k}_{0}(\sted)$. These conditions are obviously equivalent to the vanishing of $B^{(p, k-p)}$ whenever $0< p< k$.
For $\be \in \cano^{k}$, $\bj(\re \be) + \j\bj(\im \be) = \bj(\be) = \j \be = -\im \be + \j \re \be$, so $\be = \re \be - \j\bj(\re \be)$. Since $\be$ is symmetric so are $\re \be$ and $\bj(\re \be)$, and hence, by the preceeding, $\re \be$ and $\im \be$ are in $S^{k}_{0}(\sted)$.

Since $\bj$ preserves $S^{k}_{0}(\sted)$, it is a complex structure on $S^{k}_{0}(\sted)$. This means that $S^{k}_{0}(\sted)\tensor_{\rea}\com$ decomposes into $(1,0)$ and $(0,1)$ parts with respect to the action of $\bJ$ and so for $B \in S^{k}_{0}(\sted)$ it makes sense to speak of its $(1,0)$ part $\tfrac{1}{2}(B - \j \bJ(B))$. For $Z \in \ste^{0, 1}$ there holds $\imt(Z)\bj(B) = \imt(J(Z))B = -\j\imt(Z)B$, so $\imt(Z)(B - \j\bj(B)) = 0$. Hence $B - \j \bJ(B) \in \cano^{k}$ if $B\in S^{k}_{0}(\sted)$. The map $B \to B^{(k, 0)}$ also sends $S^{k}_{0}(\sted)$ to $\cano^{k}$, and it is claimed that $2B^{(k, 0)} = B - \j \bj(B)$. By the complete symmetry of $B$ there holds
\begin{align}\label{rs1}
\begin{split}
&2^{k}B^{(k,0)}(X, \dots, X) = B(X-iJX, \dots, X - iJX) = \sum_{s = 0}^{k}(-i)^{s}\tbinom{k}{s}B(X, \dots, X, \underbrace{JX, \dots, JX}_{\text{$s$ times}})\\
&= \sum_{s = 0}^{\llcorner k/2 \lrcorner}(-1)^{s}\tbinom{k}{2s}B(X, \dots, X, \underbrace{JX, \dots, JX}_{\text{$2s$ times}}) +i \sum_{s = 1}^{\ulcorner k/2 \urcorner}(-1)^{s}\tbinom{k}{2s-1}B(X, \dots, X, \underbrace{JX, \dots, JX}_{\text{$2s-1$ times}}).
\end{split}
\end{align}
Using the second equation of \eqref{jb0} in \eqref{rs1} yields
\begin{align}\label{rs3}
\begin{split}
&2^{k}B^{(k, 0)}(X, \dots, X) = \left(\sum_{s = 0}^{\llcorner k/2 \lrcorner}\tbinom{k}{2s}\right)B(X, \dots, X) -i \left(\sum_{s = 1}^{\ulcorner k/2 \urcorner}\tbinom{k}{2s-1}\right)\bJ(B)(X, \dots, X)\\
& = 2^{k-1}B(X, \dots, X) - i2^{k-1}\bJ(B)(X, \dots, X).
\end{split}
\end{align}
Polarizing \eqref{rs3} shows that for $B \in S^{k}_{0}(\sted)$ there holds $2B^{(k,0)} = B - i\bJ(B)$, so that $B$ is the real part of an element of $\cano^{k}$, and, similarly, for $\be \in \cano^{k}$, there holds $\be = 2(\re \be)^{(k,0)}$.
\end{proof}

It follows immediately from Lemma \ref{apolaritylemma} that on a Riemann surface $(M, [h], J)$, a smooth section $B$ of $\symk$ (resp. $\symkv$) is the real part of a smooth section of $\cano^{k}$ (resp. $\cano^{-k}$) if and only if it is completely $[h]$-trace-free, in which case $2B^{(k, 0)} = B - \j \bJ(B)$ and $\bJ(B)$ is also in $\Ga(S^{k}_{0}(\ctm))$.

\subsection{}
It follows from the identities $(\lie_{X}J)_{\bar{\al}}\,^{\be} =  2\j \delbar_{\bar{\al}} X^{\be}$ and $(\lie_{X}J)_{\al\be} = 2\j D_{\al}X^{\flat}_{\be}$ that on a Riemann surface $(M, J)$ the following conditions on $X \in \Ga(TM)$ are equivalent: $X$ is conformal Killing; $X$ is an infinitesimal automorphism of $J$; and $X^{(1,0)}$ is holomorphic. Since by the Riemann-Roch theorem a compact Riemann surface of genus greater than $1$ admits no non-zero holomorphic vector fields, on such a surface every conformal Killing vector field is identically $0$. These observations are generalized to sections of $S^{k}_{0}(\ctm)$ and $S^{k}_{0}(TM)$ by Lemma \ref{kdifferentialslemma}.

\begin{lemma}\label{kdifferentialslemma}
Let $(M, [h], J)$ be a Riemann surface and $k > 0$.
\begin{list}{(\arabic{enumi}).}{\usecounter{enumi}}
\renewcommand{\theenumi}{Step \arabic{enumi}}
\renewcommand{\theenumi}{(\arabic{enumi})}
\renewcommand{\labelenumi}{\textbf{Level \theenumi}.-}
\item\label{cd1} If $k > 1$, a section $B \in \Ga(S^{k}(\ctm))$ is the real part of a holomorphic section of $\cano^{k}$ if and only if it is a trace-free Codazzi tensor. In this case $(DB)^{(k+1, 0)} = DB^{(k, 0)}$.
\item\label{cd2} A one-form is the real part of a holomorphic section of $\cano^{1}$ if and only if it is closed and co-closed. In this case $(DB)^{(2, 0)} = DB^{(1, 0)}$. 
\item\label{cd3} A section $B \in \Ga(S^{k}(TM))$ is the real part of a holomorphic section of $\cano^{-k}$ if and only if it is a conformal Killing tensor. 
\end{list}
\end{lemma}

\begin{proof}[Proof of Lemma \ref{kdifferentialslemma}]
By Lemma \ref{apolaritylemma}, $B \in \Ga(S^{k}(\ctm))$ is the real part of a smooth section of $\cano^{k}$ if and only if $B = 2\re B^{(k, 0)}$, in which case $B \in \Ga(S^{k}_{0}(\ctm))$. For $B \in \Ga(S^{k}_{0}(\ctm))$, since $\clie_{h}(B)$ and $\div_{h}(B)$ are $[h]$-trace-free, there follow from \eqref{dsi} and \eqref{dsi2},
\begin{align}
\label{delb} & D^{1,0} B^{(k,0)} = \clie_{h}(B)^{(k+1, 0)}, & & \delbar B^{(k, 0)} = \overline{h^{(1,1)}}\tensor \div_{h}(B)^{(k-1,0)},& &\text{if}\,\, k > 1,\\
\label{db2}  & D^{1,0} B^{(1, 0)} = \clie_{h}(B)^{(2, 0)},  & & 2\delbar B^{(1,0)} = dB + \div_{h}(B)\overline{h^{(1,1)}},& &\text{if}\,\, k = 1 .
\end{align}
Comparing \eqref{delb} and \eqref{db2} with \eqref{dsi} and \eqref{dsi2} of Lemma \ref{codazzilemma} shows \ref{cd1} and \ref{cd2}. In the case of \ref{cd1} or \ref{cd2} then $DB \in \Ga(S^{k}_{0}(\ctm))$ and it follows from the definition of $\bj$ that $\bj DB = D\bj B$. By Lemma \ref{apolaritylemma}, $2B^{(k, 0)} = B - \j \bj B$, so $2DB^{(k, 0)} = DB - \j D\bj B =  DB - \j \bj D B = (DB)^{(k+1, 0)}$, the last equality also by Lemma \ref{apolaritylemma}.

Again by Lemma \ref{apolaritylemma} if $B\in \Ga(S^{k}(TM))$ is to be the real part of a section of $\cano^{-k}$ it must be trace-free, in which case $2B^{(k, 0)} = B - i\bJ(B)$, so that $B$ is the real part of a holomorphic section of $\cano^{-k}$ if and only if $\delbar B^{(k, 0)} = 0$. Since raising and lowering indices interchanges type, it follows from \eqref{delb} that $\delbar B^{(k, 0)} = 0$ if and only if $0 = \delbar B^{\flat\,(0, k)} = \clie_{h}(B^{\flat})^{(0, k+1)}$. Hence, $B$ is the real part of a holomorphic section if and only if $\clie_{h}(B^{\flat}) = 0$, or, what is by definition the same, $B$ is a conformal Killing tensor.
\end{proof}

\begin{lemma}\label{flatmetriclemma}
Let $(M, [h], D)$ be a Riemann surface and $\si$ the real part of a holomorphic section of $\cano^{k}$ which is not identically zero. View $\si$ as a tensor of rank $|k|$, covariant or contravariant according to whether $k$ is positive or negative. For any $h \in [h]$ there hold
\begin{align}\label{keromlap}
&2\lap_{h}\si = k\sR_{h} \si,& &\lap_{h}|\si|_{h}^{2} = 2|D\si|^{2}_{h}+ k\sR_{h} |\si|_{h}^{2}.
\end{align}
On the subset $M^{\ast} = \{|\si|_{h}^{2} \neq 0\}$, which is the complement of a discrete set of points, there hold 
\begin{align}
\label{kato}
&2|d|\si||_{h}^{2} = |D\si|_{h}^{2},& &\lap_{h}\log|\si|^{2}_{h} = k \sR_{h}.
\end{align} 
When $k \neq 0$, the metric $\sth_{ij} \defeq |\si|_{h}^{2/k}h_{ij}$ on $M^{\ast}$ is flat. If $M$ is a torus then $M^{\ast} = M$.
\end{lemma}
\begin{proof}
Tracing the (K\"ahler) identity $\sR_{ijp}\,^{l}J_{k}\,^{p} = \sR_{ijk}\,^{p}J_{p}\,^{l}$ and using the algebraic Bianchi identity yields $\omega^{pq}\sR_{pqij} = -2J_{i}\,^{p}\sR_{pj}$, from which follows $\sR_{\al\bar{\be}\ga}\,^{\ga} = \sR_{\al\bar{\be}} = (\sR_{h}/2)h_{\al\bar{\be}}$. For $s \in \Ga(\cano^{k})$ there results $2D_{[\al}D_{\bar{\be}]}s = -kR_{\al\bar{\be}\ga}\,^{\ga}s = - (k/2)\sR_{h} h_{\al\bar{\be}}s$. View $\si$ as the real part of a holomorphic section $s$ of $\cano^{k}$. That $s$ be holomorphic implies the second equality of  
\begin{align*}
\begin{split}
2\lap_{h}s &= 2h^{ij}D_{i}D_{j}s 
= 2h^{\bar{\al}\be}D_{\bar{\al}}D_{\be}s  = h^{\bar{\al}\be}\left(2D_{\be}D_{\bar{\al}}s + k \sR_{h}h_{\al\bar{\be}}\right)s = k\sR_{h} s.
\end{split}
\end{align*}
Taking the real part shows the first equation of \eqref{keromlap}, from which the second equation of \eqref{keromlap} follows. By \ref{cd1} and \ref{cd2} of Lemma \ref{kdifferentialslemma}, $(D\si)^{(k+1, 0)} = D^{1,0}s = Ds$, and so  $|D\si|_{h}^{2} = 2|Ds|_{h}^{2}$. From $d|\si|_{h}^{2} = 2D(s\bar{s}) = 2(\bar{s}Ds + sD\bar{s})$ there results $|d|\si|^{2}|^{2} = 8|s|^{2}|Ds|^{2} = 2|\si|^{2}|D\si|^{2}$, which proves the first equality of \eqref{kato}. On $M^{\ast}$ there holds by \eqref{keromlap} and the first equality of \eqref{kato},
\begin{align*}
\lap_{h}\log|\si|^{2}_{h} = 2|\si|^{-2}_{h}\left(\lap_{h}|\si|_{h}^{2} - 4|d|\si||_{h}^{2}\right) =  2|\si|^{-2}_{h}\left(\lap_{h}|\si|_{h}^{2} - 2|D\si|_{h}^{2}\right)  = k\sR_{h}. 
\end{align*}
That $\sth$ is flat follows from \eqref{conformalscalardiff} and the second equality of \eqref{kato}. If $M$ is a torus, then $\si$ is parallel for a flat representative of $[h]$, so has constant norm, which is not zero, as $\si$ is not identically zero, and so $M^{\ast} = M$.
\end{proof}

Note that $\sth$ does not depend on the choice of $h \in [h]$, and is determined by the requirement that $|\si|_{\sth}^{2} = 1$. Corollary \ref{rrcorollary} is the specialization to Riemann surfaces of the results of \cite{Kobayashi-holomorphicsymmetric}.

\begin{corollary}\label{rrcorollary}
If $M$ is a sphere then there is no non-trivial trace-free Codazzi tensor nor any non-trivial harmonic one-form. If $M$ is a torus then any trace-free Codazzi tensor, harmonic one-form, or conformal Killing tensor is parallel with respect to a flat metric conformal to $h$. If a Riemann surface $(M, [h])$ is compact with genus $g > 1$ then any conformal Killing tensor is identically $0$. 
\end{corollary}
\begin{proof}
These claims follow either from Riemann Roch together with Lemma \ref{kdifferentialslemma}, or from the maximum principle applied to \eqref{keromlap} for a constant scalar curvature representative $h \in [h]$. 
\end{proof}

\subsection{}\label{khodgesection}
For an oriented smooth surface $M$, let $\diff(M)$ be its group of diffeomorphisms viewed as a topological group in the $\cinf$ compact-open topology, and let $\diff^{+}(M)$ be the subgroup of orientation-preserving diffeomorphisms. The connected component $\diff_{0}(M)$ of the identity of $\diff(M)$ is evidently contained in $\diff^{+}(M)$, and comprises the diffeomorphisms of $M$ smoothly isotopic to the identity (see Corollary $1.2.2$ of \cite{Banyaga}). Let $\jpremod(M)$ be the space of complex structures on $M$ inducing the given orientation, with the topology of $\cinf$-convergence. The group $\diff^{+}(M)$ acts on $\jpremod(M)$ by pullback and the quotient $\jpremod(M)/\diff_{0}(M)$ is the \textbf{Teichm\"uller space} $\teich(M)$. The oriented mapping class group $\map^{+}(M) \defeq \diff^{+}(M)/\diff_{0}(M)$ acts on $\teich(M)$ with quotient $\jpremod(M)/\diff_{0}(M)$, the \textbf{moduli space of complex structures} on $M$. For background on these spaces from a point of view compatible with that here see \cite{Earle-Eells} and \cite{Wolf-teichmuller}.

Let $\prekhodge(M)$ be the space comprising pairs $(J, B) \in \jpremod(M)\times \Ga(S^{k}(\ctm))$ such that with respect to the decomposition of tensors by type determined by $J$, $B^{(k,0)}$ is $\delbar_{J}$-holomorphic. The group $\diff^{+}(M)$ acts by pullback on $\prekhodge(M)$. Denote the quotient by $\khodge(M) = \prekhodge(M)/\diff^{+}(M)$. The projection $\prekhodge(M) \to \jpremod(M)$ commutes with the $\diff^{+}(M)$ action so descends to a projection $\rho:\khodge(M) \to \teich(M)$ sending the equivalence class $[J, B]$ to the equivalence class $[J]$. If a representative $J \in [J]$ is chosen then the fiber $\rho^{-1}([J])$ is identified with the space $H^{0}(M, \cano_{J}^{k})$ of $J$-holomorphic $k$-differentials. The space $\khodge(M)$ also will be referred to as \textit{the vector bundle over the Teichm\"uller space of $M$ the fiber of which over $[h]$ comprises the $k$-holomorphic differentials with respect to the complex structure induced by $[h]$}. (This makes sense for negative $k$ if $S^{k}(\ctm)$ is replaced by $S^{|k|}(TM)$). 

\subsection{}\label{singularmetricsection}
Here are recalled some aspects of the correspondence between holomorphic differentials and singular flat metrics which are motivating for the discussion in section \ref{modulisection} of the action of $GL^{+}(2, \rea)$ on the deformation space of strictly convex flat real projective structures. Related background can be found in many places, e.g. \cite{Masur-Tabachnikov}, \cite{Troyanov-coniques}, \cite{Troyanov-moduli}, or \cite{Viana-ergodic}. 

For $\tau > 0$, the metric $dr^{2} + r^{2}dt^{2}$ on the cone $V_{\tau} = \{(r, t): r \geq 0, t \in [0, \tau)\}$ is Euclidean away from the vertex (where $r = 0$), where it is said to have a \textbf{conical singularity} of \textbf{cone angle} $\tau$. 
The change of variables $z = ((b+1)r)^{1/(b + 1)}e^{ \j t/(b+1)}$ identifies the cone $V_{2\pi(b + 1)}$ isometrically with the singular metric $|z|^{2b}|dz|^{2}$ on $\com$. If there are integers $\be$ and $k$ such that $b = \be/k > -1$, this metric is that determined as in Lemma \ref{flatmetriclemma} by the holomorphic $k$-differential $z^{\be}dz^{k}$. Conversely, this shows how to associate to a singularity with cone angle $2\pi(\be/k + 1)$ a holomorphic $k$-differential with a zero of order $\be/k$ at the singularity.

By a \textbf{flat Euclidean structure} is meant an atlas of charts in which the transition functions are restrictions of Euclidean isometries. Such a structure determines a positive homothety class of flat Riemannian metrics, and an underlying flat real affine structure. If the transition functions can be chosen to be orientation preserving, then it determines also a corresponding complex structure, and flat complex affine structure. For a positive integer $k$, a \textbf{$1/k$-translation surface} is a compact Riemann surface $M$ equipped with a flat Euclidean structure on the complement $M^{\ast}$ in $M$ of a finite subset $\sing(M) \subset M$ generated by an atlas $\{z_{i}:U_{i} \to \com\}$ for which the transition functions have the form $z_{i} = e^{2\pi\j m_{ij}/k} z_{j} + u_{ij}$ for some integers $0 \leq m_{ij} \leq k-1$ and some complex numbers $u_{ij}$. In particular, the linear part of the affine holonomy around each $p \in \sing(M)$ is contained in the finite subgroup $\integer/k\integer \subset SO(2)$. In the cases $k$ is $1$ or $2$, such a surface is called a \textbf{translation} or \textbf{half-translation} surface, which explains the terminology. The $k$-differentials $dz_{i}^{k}$ and $dz_{j}^{k}$ agree on the overlaps $U_{i}\cap U_{j}$, so patch together to give a holomorphic $k$-differential on $M^{\ast}$. A point of $M^{\ast}$ will be called a \textbf{regular} point of $\si$, while points in $\sing(M)$ will be called \textbf{singular}. From the form of the transition funtions it follows that in a neighborhood of $p \in \sing(M)$ there is a chart in which a flat metric $\sth$ representing the flat Euclidean structure is isometric to a conical singularity of cone angle $2\pi(\be/k + 1) > 0$ for some integer $\be$. Then the holomorphic $k$-differential constructed on $M^{\ast}$ extends to the singular points via the local model described in the preceeding paragraph.

Let $(M, [h], J)$ be a compact Riemann surface of genus $g > 1$ and, for $k > 0$, let $\si = B^{(k,0)}$ be a non-trivial holomorphic $k$-differential. By Lemma \ref{flatmetriclemma}, $\sth = |B|^{2/k}_{h}h$, which does not depend on the choice of $h \in [h]$, is a flat metric on the complement $M^{\ast}$ of the zero set of $M$.  Around a regular point $p_{0}$ choose a local holomorphic coordinate $w$ (such that $w = 0$ corresponds to $p_{0}$) and write $\si = \phi(w) dw^{k}$ for a holomorphic function $\phi(w)$. Choose a branch of $\phi^{1/k}(w)$ near $w = 0$ and define a new coordinate, said to be \textbf{adapted} to $\si$, by $z(p) = \int_{p_{0}}^{p}\phi(w)^{1/k}\,dw$. In the $z$ coordinate, $\si = dz^{k}$. If $\tilde{z}$ is another coordinate constructed in this way in a neighborhood of $p_{0}$ then, since $dz^{k} = d\tilde{z}^{k}$ on the overlap, there are an integer $0 \leq m \leq k-1$ and a complex number $\be$ such that $\tilde{z} = e^{2\pi\j m/k} z + \be$. Consequently the local charts constructed in this way determine on $M^{\ast}$ a flat Euclidean structure, which makes $M$ a $1/k$-translation surface, and for which the underlying homothety class of flat Riemannian metrics is generated by the metric $\sth$. In a neighborhood of a singular point $p_{0}$ of $\si$ choose as before a local holomorphic coordinate $w$ and write $\si(w) = \phi(w)dw^{k}$, where now it is supposed that $\phi(w)$ is holomorphic with a zero of order $\be$ at $0$. Define a coordinate $z(p)$ adapted to $\si$ by 
\begin{align*}
z(p) = \left(\tfrac{\be + k}{k}\int_{p_{0}}^{p}\phi(w)^{1/k}\,dw\right)^{\tfrac{k}{\be + k}}.
\end{align*}
The coordinate $z$ is determined up to multiplication by a $(\be + k)$th root of unity, and $\si = \phi(w)dw^{k} = z^{\be}dz^{k} = (\tfrac{k}{\be+k}d(z^{\be/k + 1}))^{k}$. In an adapted coordinate around a singular point $p$ of $\si$ of order $\be$ the flat metric $\sth$ has the form $|z|^{2\be/k}|dz|^{2}$, and so $p$ is a cone point of angle $2\pi(\be/k + 1)$, and the linear holonomy of $\sth$ around $p$ is a rotation of angle $2\pi \be/k$. By Riemann Roch the sum of the orders $\be(p)$ of the zeroes $p$ of $\si$ is $2k(g-1)$, and so the cone angles $\tau(p) = 2\pi(\be(p)/k + 1)$ of $\sth$ satisfy the relation $\sum_{p \in M}(\tau(p)/2\pi - 1) = 2(g-1)$. The same conclusion follows from a version of the Gau\ss-Bonnet theorem for metrics with conic singularities (see \cite{Troyanov-moduli}).

The preceeding establishes a bijective correspondence between $1/k$-translation surfaces and Riemann surfaces with a holomorphic $k$-differential. 

A pair $(J, \si)$ for which $\si$ is a $J$-holomorphic $k$-differential generates a complex curve in Teichm\"uller space, and, moreover, in $\khodge(M)$, as follows. The curve in $\teich(M)$ comprises those conformal structures generated by singular flat metrics real affinely equivalent to the flat metric $\sth$ determined by $(J, \si)$. The evolution of the $k$-differential $\si$ along this curve in $\teich(M)$ is determined tautologically by the requirement that the singular flat metric associated to the evolved $k$-differential be a member of a conformal structure representing the point in the curve over which it lies. These curves should be relevant to understanding compactifications of the space of strictly convex flat real projective structures on a compact surface of genus $g > 1$ (see section \ref{modulisection}).

These curves can be described analytically in terms of an action of the group $GL^{+}(2, \rea)$ on $\prekhodge(M)$, with respect to which they are images of the hyperbolic disk $GL^{+}(2, \rea)/C^{+}O(2)$. The structure of this $GL^{+}(2, \rea)$ action on $\prekhodge(M)$ has been intensively studied in the $k = 1, 2$ cases; see \cite{Kerckhoff-Masur-Smillie}, \cite{Masur-Tabachnikov}, and \cite{Herrlich-Schmithusen} for background and references. An identification of real vector spaces $\com \simeq \rea^{2}$ is fixed so that $g \in GL^{+}(2, \rea)$ acts on $z = x + \j y\in \com$ real linearly by $g \cdot z = (ax + by) + \j(cx + dy)$. The complex field $\comt$ acting by multiplication is identified with the oriented conformal subgroup $C^{+}O(2)$. That is $z = re^{\j\theta}$ corresponds to $a = r\cos\theta = d$, $-b = r\sin \theta = c$. Given $(J, \si) \in \prekhodge(M)$, let $\{z_{i}\}$ be an atlas on $M^{\ast}$ of $J$-holomorphic charts adapted to $\si$. This atlas makes $M^{\ast}$ a flat Euclidean manifold with linear holonomy in $\integer/k\integer$. For $g \in GL^{+}(2, \rea)$ the collection $\{\tilde{z}_{i} = g z_{i}\}$ determines a flat Euclidean structure which has also holonomy in $\integer/k\integer$, so a structure of $1/k$-translation surface on $M$, the cone angles of which are same as those of the original structure determined by $(J, \si)$. This $1/k$-translation surface corresponds to a pair $g \cdot (J, \si) = (g\cdot J, g\cdot \si) \in \prekhodge(M)$, and this defines the desired $GL^{+}(2, \rea)$ action on $\prekhodge(M)$. In the case $k = 1$, $g\cdot \si$ is given in $z$-coordinates by the real linear action of $g$ on $\si$; that is $g \cdot \si = (a \re \si + b \im \si) + \j (c \re \si + d\im \si)$, but in general, with respect to a coordinate $z$ adapted to $\si$, the new $k$-differential $g \cdot \si$ is not given by a linear action. The underlying complex structure $g\cdot J$ will equal the original one if $g \in C^{+}O(2) = \comt$. In this case $g \cdot \si$ is expressible in adapted coordinates in terms of a linear represnetation of $C^{+}O(2)$; namely, it is given by the product by $z^{k}\si= z^{k}B^{(k, 0)} = r^{k}(\cos (k \theta) B - \sin (k \theta) \bj(B))^{(k, 0)}$, where $z$ is the complex number corresponding to $g$. This action of $C^{+}O(2)$ is considered in section \ref{gaugesection}. The more interesting actions of hyperbolic and parabolic elements of $GL^{+}(2, \rea)$ are more difficult to describe; their study is left for the future. 

If $\Phi \in \diff(M)$ then $z_{i}\circ \Phi$ are $\Phi^{\ast}(\si)$ adapted coordinates with respect to the conformal structure $\Phi^{\ast}([h])$, and since $g(z_{i}\circ F) = (gz_{i})\circ F$, it follows that the action of $GL^{+}(2,\rea)$ on $\prekhodge(M)$ is diffeormorphism equivariant in the sense that $\Phi^{\ast}(g \cdot (J, \si)) = g \cdot \Phi^{\ast}(J, \si)$, so descends to an action on $\khodge(M)$ which commutes with $\map^{+}(M)$.

\section{AH structures on surfaces}\label{ahsection}
In this section are given the basic definitions related to AH structures. Though adapted to the peculiarities of the two-dimensional case, the exposition is consistent with the conventions of \cite{Fox-ahs}.
\subsection{}
Two affine connections are \textbf{projectively equivalent} if they have the same unparameterized geodesics in the sense that the image of any geodesic of one connection is the image of a geodesic of the other connection. This is the case if and only if the symmetric part of their difference tensor is pure trace. A \textbf{projective structure} $\en$ is an equivalence class of projectively equivalent affine connections. For a torsion-free affine connection $\nabla$ on a surface there vanishes the usual projective Weyl tensor, or, what is the same,
\begin{align}\label{rijklpij}
\begin{split}
R_{ijk}\,^{l}  &= \delta_{i}\,^{l}R_{(jk)} - \delta_{j}\,^{l}R_{(ik)} - R_{[ij]}\delta_{k}\,^{l}.
\end{split}
\end{align}
The \textbf{projective Cotton} tensor $C_{ijk} \defeq -\nabla_{i}R_{(jk)} + \nabla_{j}R_{(ik)} + \tfrac{1}{3}\nabla_{k}R_{[ij]}$ does not depend on the choice of $\nabla \in \en$. On a surface, since a $3$-form vanishes, there hold $C_{[ijk]} = 0$ and $\nabla_{[i}C_{jk]l} = 0$.

\subsection{}
The definition of an AH structure $(\en, [h])$ on an $n$-manifold was given in the first paragraph of section \ref{introsection} of the introduction. Let $H_{ij}$ be the normalized representative of $[h]$. An equivalent definition is that for each $\nabla \in \en$ there be a one-form $\si_{i}$ such that $\nabla_{[i}H_{j]k} = 2\si_{[i}H_{j]k}$. In two dimensions, by Lemma \ref{weylcriterion} the tensor $\nabla_{[i}H_{j]k}$ is pure trace, so there is always such a one-form. Hence in two dimensions an AH structure could simply be defined to be a pair $(\en, [h])$, the necessary compatibility being automatic. If $M$ is oriented, an alternative way to see this is the following. On a surface with K\"ahler structure $(h, J, \om)$, the Hodge star operator on one-forms is $(\star \si)_{i} = -\si_{p}J_{i}\,^{p}$, while on two-forms it is $\star \be = \tfrac{1}{2}\om^{\sharp\,ij}\be_{ij}$. Because any two-form is a multiple of any other, given $\nabla \in \en$ there is a one-form $\si_{i}$ such that $\nabla_{[i}h_{j]k} = \om_{ij}\si_{k}$. An easy computation shows 
\begin{align}\label{omsi}
\om_{ij}\si_{k} = -2(\star \si)_{[i}h_{j]k},
\end{align} 
so $\nabla_{[i}h_{j]k} =-2(\star \si)_{[i}h_{j]k}$.

\subsection{}
Given an AH structure $(\en, [h])$, there is a unique torsion-free $\nabla \in \en$ such that $\nabla_{i}H_{jk} = 0$; given any torsion-free $\bnabla \in \en$ with  $\bnabla_{[i}H_{j]k} = 2\si_{[i}H_{j]k}$ it is given by $\nabla = \tnabla - 2\si_{(i}\delta_{j)}\,^{k}$. This distinguished representative of $\en$ is called the \textbf{aligned} representative of the AH structure.

From now on, except where stated otherwise, indices and raised and lowered using $H_{ij}$ and the dual bivector $H^{ij}$. Because $\det H_{ij} = 1$ there holds $H^{pq}\bnabla_{i}H_{pq} = 0$ for any $\bnabla \in \en$.  Alternative characterizations of the alignment condition are given in Lemma \ref{special}, the proof of which is straightfoward, using the identity
\begin{align}\label{lcidentity}
A_{ijk} = A_{i(jk)} + A_{j(ik)} - A_{k(ij)} + A_{[ij]k} + A_{[ki]j} - A_{[jk]i},
\end{align}
valid for any covariant $3$-tensor. 
\begin{lemma}\label{special}
Let $\en$ be a projective structure and $[h]$ a conformal structure on a surface $M$. There is a unique torsion-free representative $\nabla \in \en$ satisfying any one of the following equivalent conditions
\begin{enumerate}
\item $\nabla_{[i}H_{j]k} = 0$.
\item $\nabla_{i}H_{jk} = \nabla_{(i}H_{jk)}$.
\item $\nabla_{p}H^{ip} = 0$.
\item $H^{pq}\nabla_{p}H_{qi} = 0$.  That is, $\nabla_{i}H_{jk}$ is completely trace-free.
\item For any $h \in [h]$ there holds $2h^{pq}\nabla_{p}h_{qi} = h^{pq}\nabla_{i}h_{pq}$.
\item For any $h \in [h]$ there holds $\nabla_{[i}h_{j]k}  = 2\ga_{[i}h_{j]k}$ with $4\ga_{i} = h^{pq}\nabla_{i}h_{pq}$.
\end{enumerate}
\end{lemma}
\noindent
Henceforth, except where stated otherwise, $\nabla$ denotes the aligned representative of an AH structure. While it may seem perverse speak of the projective structure $\en$ if one works only with a distinguished representative $\nabla \in \en$, later developments will show the utility of the perspective.

\subsection{}
The \textbf{cubic torsion} of an AH structure is the tensor $\bt_{ij}\,^{k}$ defined in terms of an arbitrary representative $\nabla \in \en$ by setting $\bt_{ij}\,^{p}H_{pk}$ equal to the completely trace-free part of $\nabla_{(i}H_{jk)}$. For the aligned representative $\nabla \in \en$ the cubic torsion is just $\nabla_{i}H_{jk} = \nabla_{(i}H_{jk)} = \bt_{ijk}$. An AH structure for which $\bt_{ij}\,^{k} \equiv 0$ is a \textbf{Weyl structure}. The aligned representative of a Weyl structure is what is usually called a Weyl connection.

\subsection{}
In the split signature case, an appropriate finite cover $\hat{M}$ of $M$ is orientable and such that the null cone of the lifted metric on $\hat{M}$ is orientable; in this case the null cone of the split signature metric is a pair of transverse nowhere vanishing line fields, and so $\hat{M}$ has Euler characteristic $0$, and hence $M$ has as well, and so if $M$ is compact it is a torus or a Klein bottle. The study of Riemannian Einstein AH structures on surfaces makes use of Hodge theory, the associated complex structure, and so forth. While the study of split signature Einstein AH structures on surfaces is also interesting, it requires a different set of techniques and will be mostly ignored here, although it will always be indicated when the Riemannian hypothesis is necessary.

\subsection{}
The basic example of an AH structure is the following. A hypersurface immersion in flat affine space is \textbf{non-degenerate} if its second fundamental form (which takes values in the normal bundle) is non-degenerate. If the immersion is also co-oriented the second fundamental form determines a conformal structure on the hypersurface. A choice of subbundle transverse to the immersion induces on the hypersurface a torsion-free affine connection, and there is a unique choice of transverse subbundle such that the induced connection is aligned with respect to the conformal structure determined by the second fundamental form and the co-orientation. This choice of transverse subbundle is the \textbf{affine normal subbundle}. That this definition coincides with the customary one is proved in section $4.2$ of \cite{Fox-ahs}.

The second fundamental form with respect to a particular vector field spanning the affine normal subbundle is a pseudo-Riemannian metric on the hypersurface. Usually an \textbf{(equi)affine normal vector field} $\nm$ is distinguished by fixing a volume form on the ambient affine space and requiring that the volume density induced by the metric $h_{ij}$ determined by the vector field agree with the volume density induced via interior multiplication of the vector field with the chosen volume form. Though it is often omitted, the prefix \textit{equi} ought to be included because this construction is only invariant under equiaffine transformations of the ambient space (those preserving the given volume form). The metric $h_{ij}$ determined by the affine normal is called the \textbf{(equi)affine} or \textbf{Blaschke} metric. The affine normal field admits the following description (that this is equiaffinely invariant is not self-evident; rather it follows from the definition in the previous paragraph, which is manifestly so). Fix the standard flat Euclidean metric $\delta_{IJ}$ on the ambient $\rea^{3}$ and let $g_{ij} = i^{\ast}(\delta)_{ij}$ be the induced metric on the immersed, co-oriented, non-degenerate hypersurface $i:\Sigma \to \rea^{3}$. Let $N$ be the unit Euclidean normal consistent with the given co-orientation, and let $\Pi_{ij}$ be the second fundamental form defined with respect to $N$. The Gaussian curvature $K$ is the function $\det \Pi \tensor (\det g)^{-1}$, and the equiaffine metric $h_{ij}$ has the form $h_{ij} = |K|^{-1/4}\Pi_{ij}$. Let $\rad$ be the radial (position) vector field on $\rea^{3}$. Along $i(\Sigma)$ the equiaffine vector field $\nm$ satisfies $\lap_{h}\rad = 2\nm$. 

For the convenience of a reader not familiar with affine geometry, the translation of these definitions into local coordinates is recalled briefly. For details, consult \cite{Calabi-affinelyinvariant} or the textbook \cite{Nomizu-Sasaki}. Locally a non-degenerate hypersurface $\Sigma$ is given as a graph $z = f(x)$ where $x^{i} \in \Omega \subset \rea^{2}$. Let $\pr$ denote the flat connection on $\Omega$ with respect to which the $dx^{i}$ are parallel, and write $f_{i_{1}\dots i_{k}} = \pr_{i_{1}}\dots \pr_{i_{k}}f$ and $\hess f = f_{ij}$. Also let $f^{ij}$ be the tensor inverse to $f_{ij}$. Define $\H(f)$ by $\det \hess f= \H(f)(dx^{1}\wedge dx^{2})^{2}$. The normalized representative $H_{ij}$ of the AH structure $(\en, [h])$ induced on $\Sigma$ is given by $H_{ij} = |\H(f)|^{-1/4}f_{ij}$, while the aligned representative $\nabla \in \en$ is given by $\nabla = \pr + \tfrac{1}{4}f_{ij}f^{pq}f_{pqr}f^{rk}$. The affine normal is $\nm = |\H(f)|^{1/4}\left(Z - \tfrac{1}{4}f^{pq}f_{pqr}f^{ri}X_{i}\right)$ in which $X_{i} = \tfrac{\pr}{\pr x^{i}} - f_{i}\tfrac{\pr}{\pr z}$ and $Z = \tfrac{\pr}{\pr z}$. The Euclidean unit normal is $N = (1 + f^{p}f_{p})^{-1/2}\left(Z - (1 + f^{p}f_{p})^{-1}f^{i}X_{i}\right)$, in which $f^{i} = g^{ip}f_{p}$, and the second fundamental form $\Pi_{ij}$ induced by $N$ is $\Pi_{ij}= (1 + f^{p}f_{p})^{-1/2}f_{ij}$, while the Gaussian curvature is $K = (1 + f^{p}f_{p})^{-2}\H(f)$. Along $i(\Sigma)$ the radial field $\rad$ is equal to $x^{p}X_{p} - f^{\ast}Z$, in which $f^{\ast} = x^{p}f_{p} - f$ is the Legendre transform of $f$. The hypersurface is an \textbf{affine hypersphere} if its affine normal subbundles meet in a point, which may be at infinity. This holds if and only if either $\nm$ is parallel or there is a constant $S$ such that $\nm = -S\rad$. For non-zero $S$ this holds if and only if $f^{\ast}$ solves the equation $(f^{\ast})^{4}\H(f^{\ast}) = S^{4}$, in which the Hessian is taken with respect to the Legendre transformed variables $y_{i} \defeq f_{i}$ on the domain $\Omega^{\ast} = df(\Omega)$. Most of the deeper results in the study of affine hyperspheres are obtained by studying this Monge-Ampere equation. In the two-dimensional case under consideration, there is also a Weierstrass-like representation of affine hypersurfaces due to Calabi, \cite{Calabi-affinemaximal}, and Wang, \cite{Wang}, which allows the use of complex analytic methods.

The simplest affine hyperspheres which are not quadrics are the hypersurfaces 
\begin{align*}
\Sigma_{\al} = \{(X, Y, Z) \in \rea^{3}: (Y^{2}Z\cos \al - XYZ \sin \al) = 1\},
\end{align*}
for $\al \in (0, \pi)$. Writing $x = X/Z$ and $y = Y/Z$, $\Sigma_{\al}$ is the \textbf{radial graph} 
\begin{align*}
\{(x/u, y/u, -1/u) \in \rea^{3}: u(x, y) < 0\} 
\end{align*}
of the function $u(x, y) = (y^{2}\cos \al - xy \sin \al)^{1/3}$, which solves $27 u^{4}\det \hess u = \sin^{2}\al$. This $\Sigma_{\al}$ is asymptotic to the cone over the base triangle $\Omega^{\ast} = \{(x, y): u(x,y) < 0\}$; it can be written as the graph of the Legendre transform of $u$ (that is, $u$ corresponds to $f^{\ast}$ in the previous paragraph). The parameter $\al$ is not important, as these examples are all equivalent by an affine transformation, but is included to illustrate that as $\al \to \pi$ the function $u$ becomes degenerate; as the upper halfspace contains a complete affine line it supports no strictly convex negative function vanishing along its boundary.

The function $v(x, y) = -(x^{2} + y^{2})^{1/3}$ solves $v^{4}\H(v) = -(4/27)$ on $\rea^{2} \setminus\{0\}$. Its radial graph is the hypersurface $\{(X, Y, Z) \in \rea^{3}: (X^{2} + Y^{2})Z  = 1\}$, which is an affine hypersphere. This does not contradict the preceeding paragraph because the Hessian of $v$ has mixed signature (its eigenvalues are $-(2/3)v^{-2}$ and $(2/9)v^{-2}$), so in this case the equiaffine metric has split signature.

\subsection{}
Given an AH structure $(\en, [h])$ the torsion-free connection $\bnabla \defeq \nabla + \bt_{ij}\,^{k}$ satisfies $\bnabla_{i}H_{jk} = - \bt_{ijk}$, so is the aligned representative of the AH structure $(\ben, [h])$ formed by the projective structure $\ben$ generated by $\bnabla$ and the given $[h]$. This $(\ben, [h])$ is said to be \textbf{conjugate} to $(\en, [h])$. As its cubic torsion is $\bar{\bt}_{ij}\,^{k} = -\bt_{ij}\,^{k}$, the conjugate of the conjugate is the original AH structure.  

The conormal Gau\ss\, map of a non-degenerate co-oriented hypersurface immersion in flat affine space sends a point of the hypersurface to the annihilator of the space tangent at the image of the point to the hypersurface. The pullback of the flat projective structure on the projectivization of the dual to the flat affine space via this conormal Gau\ss\, map forms with the conformal structure determined by the second fundamental form and the co-orientation the AH structure conjugate to that determined by the affine normal subbundle. 

\subsection{}
The \textbf{curvature} and \textbf{Ricci curvature} of an AH structure are defined to be the curvature $R_{ijk}\,^{l}$ and Ricci curvature $R_{ij} \defeq R_{pij}\,^{p}$ of the aligned representative $\nabla$. The \textbf{scalar curvature} $R$ is the density $R \defeq R_{p}\,^{p} =  H^{pq}R_{pq}$. Sometimes, for emphasis, the qualifier \textit{weighted} will be added, and $R$ will be called the \textit{weighted scalar curvature}. It does not make sense to speak of the numerical value of $R$ because $R$ takes values in the line bundle $|\det \ctm|$; however it does make sense to speak of the vanishing of $R$ and because $|\det \ctm|$ is oriented, to speak of the positivity or negativity of $R$. An AH structure is \textbf{proper} if its weighted scalar curvature is non-vanishing. When a representative $h \in [h]$ is given there will be written $\uR_{h} = |\det h|^{-1/2}R = h^{ij}R_{ij}$.

An AH structure is \textbf{projectively flat} if $\en$ is projectively flat. If the conjugate AH structure $(\ben, [h])$ is projectively flat, then $(\en, [h])$ is \textbf{conjugate projectively flat}. 

\subsection{}
An AH structure is \textbf{exact} if there is a representative $h \in [h]$ such that $\nabla_{i}\det h = 0$ for the aligned representative $\nabla \in \en$. If there is such an $h$ it is determined uniquely up to positive homothety (on each connected component of $M$). Such an $h$ will be called a \textbf{distinguished representative} of the AH structure. For example, the AH structure induced on a hypersurface in flat affine space is always exact, and the equiaffine metric is a distinguished representative. An AH structure is exact if and only if there is a global $\nabla$-parallel non-vanishing density of non-trivial weight, for if there is such a density, then some power of it is a non-vanishing density $\mu$ such that $h_{ij} = \mu \tensor H_{ij} \in [h]$ verifies $\nabla |\det h| = \nabla \mu^{2} = 0$ (the converse is obvious). 

\subsection{}\label{faradaysection}
The \textbf{Faraday form} $F_{ij}$ of an AH structure $(\en, [h])$ is the curvature of the covariant derivative induced on the line bundle of $-1/2$-densities by the aligned representative $\nabla \in \en$. If $R_{ijk}\,^{l}$ is the curvature of $\nabla$, then by definition and the traced algebraic Bianchi identity there hold
\begin{align}\label{faradayexplicit}
2F_{ij} = R_{ijp}\,^{p} = -2R_{[ij]}. 
\end{align} 
Since $F_{ij}$ is the curvature of a connection on a trivial line bundle it is always exact. The AH structure $(\en, [h])$ is \textbf{closed} if $F_{ij} = 0$. The terminology directly extends that introduced for Weyl structures in \cite{Calderbank-faraday}.

If an AH structure $(\en, [h])$ has parallel weighted scalar curvature then either $R$ vanishes identically or $R$ vanishes nowhere and $(\en, [h])$ is exact, for if $R$ vanishes nowhere, then $2\si_{i} = -R^{-1}\nabla_{i}R$ satisfies $d\si_{ij} = 2\nabla_{[i}\si_{j]} = -R^{-1}\nabla_{[i}\nabla_{j]}R = F_{ij}$.

The \textbf{Faraday primitive} $\ga_{i}$ associated to $h \in [h]$ is the one-form $\ga_{i}$ defined by $4\ga_{i} = h^{pq}\nabla_{i}h_{pq}$. Note that the Faraday primitive associated to $h$ depends only on the positive homothety class of $h$ and not on $h$. From the Ricci identity follows $d\ga_{ij} = 2\nabla_{[i}\ga_{j]} = -F_{ij}$, so that $\ga_{i}$ is a primitive for $-F_{ij}$. If $\tilde{h}_{ij} = e^{2f}h_{ij} \in [h]$ then the corresponding one-form $\tilde{\ga}_{i}$ differs from $\ga_{i}$ by an exact one-form, $\tilde{\ga}_{i} = \ga_{i} + df_{i}$. The equivalence class $\{\ga\}$ of one-forms so determined is the \tbf{equivalence class of Faraday primitives} \textbf{induced by $(\en, [h])$}. The Faraday primitives associated to $h \in [h]$ of the AH structure $(\en, [h])$ and its conjugate $(\ben, [h])$ are the same, and so these AH structures determine the same equivalence class of Faraday primitives. In particular, an AH structure is closed (resp. exact) if and only if the conjugate AH structure is closed (resp. exact).

Most properties of the Faraday curvature of Weyl structures hold also for AH structures. For example, Lemma \ref{fcoclosedlemma} generalizes (trivially) Theorem $2.5$ of \cite{Calderbank-faraday}.
\begin{lemma}\label{fcoclosedlemma}
A definite signature AH structure on a surface is closed if and only if $\nabla^{p}F_{ip} = 0$.
\end{lemma}
\begin{proof}
By \eqref{faradayexplicit}, there holds $\nabla_{[p}\nabla_{q]}F^{pq} =  -R_{[pq]}F^{pq}- R_{pqa}\,^{a}F^{pq} = -F_{pq}F^{pq}$. Hence that $\nabla^{p}F_{ip} = 0$ implies $F_{ij}$ is $[h]$-null; in definite signature this holds if and only if $F_{ij} \equiv 0$.
\end{proof}

\subsection{}
Because  $f\hodge_{\tilde{h}}\al_{i} = \hodge_{h}\al_{i} + 2\si^{\sharp\,p}d\al_{pi} - 2\dad_{h}\al\, \si_{i}$, the \textbf{Hodge Laplacian} $\hodge_{h} \defeq d\dad_{h} + \dad_{h}d$ is not conformally invariant on one-forms. However, because $d$ is independent of the metric and $f\dad_{\tilde{h}}\al = \dad_{h}\al$, and because a form is harmonic if and only if it is closed and co-closed, the Hodge decomposition of one-forms is conformally invariant. On a compact orientable Riemannian surface the Hodge decomposition implies there is a unique representative of $\{\ga\}= \{\ga + df: f\in \cinf(M)\}$ which is co-closed. Consequently, on such a surface there is associated to an AH structure a unique positive homothety class of representative metrics $h \in [h]$ for which the associated Faraday primitive $\ga_{i}$ is co-closed with respect to $h$. Such a representative metric will be called a \textbf{Gauduchon metric}. In higher dimensions the existence for an AH structure $(\en, [h])$ of a representative of $[h]$ distinguished up to positive homothety follows from arguments of P. Gauduchon (e.g. \cite{Gauduchon} and \cite{Gauduchon-circlebundles}), and, although in two dimensions the existence of such representatives follows from the Hodge decomposition, in this case the terminology \textit{Gauduchon metric} is used also, for consistency. Note that without imposing some further normalization, such as setting the volume equal to a fixed constant, there is no naturally preferred Gauduchon metric. A distinguished metric on an exact AH structure is trivially also a Gauduchon metric.

\subsection{}\label{ahfrommetricsection}
If $(\en, [h])$ is an AH structure and $h \in [h]$ has corresponding Faraday primitive $\ga_{i}$ then the Levi-Civita connection $D$ of $h$ is related to the aligned representative $\nabla$ by
\begin{align}\label{dnabladiff}
D  = \nabla + \tfrac{1}{2}\bt_{ij}\,^{k} + 2\ga_{(i}\delta_{j)}\,^{k} - h_{ij}\ga^{\sharp\,k} = \nabla + \tfrac{1}{2}\bt_{ij}\,^{k} + 2\ga_{(i}\delta_{j)}\,^{k} - H_{ij}\ga^{k},
\end{align}
On the other hand, equation \eqref{dnabladiff} shows how to build from a metric and a one-form an AH structure pair with a given cubic torsion. Given a pseudo-Riemannian metric $h_{ij}$ with Levi-Civita connection $D$, a completely symmetric, completely trace-free covariant $3$-form $B_{ijk} = B_{(ijk)}$, and a one-form $\ga_{i}$, defining $\bt_{ij}\,^{k} = B_{ijp}h^{pk}$ and defining $\nabla$ by \eqref{dnabladiff} determines an AH structure with aligned representative $\nabla$, cubic torsion $\bt_{ij}\,^{k}$, and such that $\ga_{i}$ is the Faraday primitive associated to $h$.

\subsection{}
On an oriented surface the conformal structure underlying a Riemannian signature AH structure determines a complex structure. Lemma \ref{2dcomplexlemma} describes the compatibility between them.
\begin{lemma}\label{2dcomplexlemma}
Let $(\en, [h])$ be a Riemannian signature AH structure on an oriented surface and let $J_{i}\,^{j}$ be the complex structure determined by $[h]$. Then
\begin{align}\label{btj}
&\bt_{ij}\,^{k} = -J_{p}\,^{k}\nabla_{i}J_{j}\,^{p},& &\nabla_{[i}J_{j]}\,^{k} = 0,& &J_{[i}\,^{p}\bt_{j]p}\,^{k} = 0.
\end{align}
\end{lemma}
\begin{proof}
Fix $h \in [h]$ with Levi-Civita connection $D$ and associated Faraday primitive $\ga_{i}$. Since the tensor $B_{ijk} \defeq \bt_{ij}\,^{p}h_{pk}$ is completely trace-free, Lemma \ref{apolaritylemma} implies $J_{[i}\,^{p}B_{j]kp} = 0$, which shows the third equation of \eqref{btj}. By \eqref{omsi} there holds
\begin{align}\label{btj0}
\begin{split}
h_{kp}&\left( - \ga_{j}J_{i}\,^{p} + h_{ij}(\star \ga)^{\sharp\,p} - (\star \ga)_{j}\delta_{i}\,^{p} + \om_{ij}\ga^{\sharp\, p} \right) = \ga_{j}\om_{ki} + \ga_{k}\om_{ij} - 2(\star\ga)_{[j}h_{k]i} = 3\ga_{[i}\om_{jk]} = 0.
\end{split}
\end{align}
By the last equality of \eqref{btj}, $\bt_{ip}\,^{k}J_{j}\,^{p}  = -\bt_{ij}\,^{p}J_{p}\,^{k}$. Using \eqref{dnabladiff} and \eqref{btj0} there results
\begin{align*}
\begin{split}
0 = D_{i}J_{j}\,^{k} & = \nabla_{i}J_{j}\,^{k} - \bt_{ij}\,^{p}J_{p}\,^{k}  - \ga_{j}J_{i}\,^{k} + h_{ij}(\star \ga)^{\sharp\,k} - (\star \ga)_{j}\delta_{i}\,^{k} + \om_{ij}\ga^{\sharp\, k} = \nabla_{i}J_{j}\,^{k} - \bt_{ij}\,^{p}J_{p}\,^{k}.
\end{split}
\end{align*}
This gives the first equation of \eqref{btj}, from which the second equation of \eqref{btj} is immediate. 
\end{proof}

Thus the cubic torsion of an AH structure can be seen as measuring the failure of the aligned representative to preserve the associated complex structure. By Lemma \ref{apolaritylemma} the third equation of \eqref{btj} shows that for any $h \in [h]$ the $(3,0)$ part $B^{(3,0)}$ is a smooth section of the bundle of cubic holomorphic differentials.

From \eqref{btj} it follows also that $\nabla$ preserves $J$ if and only if the AH structure is Weyl. In \cite{Nomizu-Podesta} (see also the more easily obtained \cite{Nomizu-Sasaki}) a torsion-free affine connection on a complex manifold preserving the complex structure is said to be \textbf{affine K\"ahler} if the $(2,0)$ part of its curvature vanishes. (Though apt, the terminology is problematic because \textit{K\"ahler affine} has been used in \cite{Cheng-Yau-realmongeampere} to mean something else, in a somewhat related context; see section \ref{kahleraffinesection}). This curvature condition is automatic on a Riemann surface, so the aligned representative of a Weyl structure is an affine K\"ahler connection in this sense.

\section{Curvature of an AH structure}\label{ahcurvaturesection}
In this section there are described the basic local curvature invariants of an AH structure $(\en, [h])$ on a surface. The fundamental invariants are the weighted scalar curvature $R$, a trace-free symmetric tensor $E_{ij}$ which is a multiple of the trace-free Ricci tensor, as well as the cubic torsion $\bt_{ij}\,^{k}$ and Faraday curvature $F_{ij}$. Conceptually important are Lemmas \ref{eijvanishlemma} and \ref{selfconjugatecottonlemma} which show when $E_{ij}$ and $\bt_{ij}\,^{k}$ can be viewed as the real parts of holomorphic differentials. In \cite{Calabi-affinemaximal} Calabi develops the geometry of hypersurfaces in flat affine three space; the fundamental geometric invariants are a quadratic differential $B$ and a cubic differential $A$ which is holomorphic when $B$ vanishes; for the AH structure induced on such a hypersurface, constant multiples of the real parts of $A$ and $B$ are identified, respectively, with $\bt_{ij}\,^{k}$ and $E_{ij}$.

\subsection{}
Let $(\en, [h])$ be an AH structure on a surface $M$. In what follows indices are raised and lowered using the normalized representative $H$, $\nabla \in \en$ is the aligned representative, and $h \in [h]$ is a representative with Faraday primitive $\ga_{i}$ and Levi-Civita connection $D$ related to $\nabla$ as in \eqref{dnabladiff}.

On a manifold of dimension $n > 2$ it is necessary to consider, in addition to the Ricci trace, the trace $R_{ip}\,^{p}\,_{j}$, because in general the trace-free symmetric parts of $R_{ip}\,^{p}\,_{j}$ and $R_{ij}$ are independent, but in two dimensions tracing \eqref{rijklpij} shows that 
\begin{align}\label{qrrelated}
R_{ip}\,^{p}\,_{j} + R_{ij} = RH_{ij},
\end{align}
so in this case there is no need to speak of $R_{ip}\,^{p}\,_{j}$. It will be convenient to work instead with the trace-free symmetric tensor $E_{ij} = E_{(ij)}$ defined by 
\begin{align}\label{ricdecompose}
R_{ij} = -2E_{ij} + \tfrac{1}{2}RH_{ij} - F_{ij}.
\end{align}
The apparently unnatural coefficient $-2$ is chosen for consistency with the conventions of \cite{Fox-ahs}. Substituting $E_{ij}$ into \eqref{rijklpij} gives $R_{ijkl} = -4H_{l[i}E_{j]k} + RH_{l[i}H_{j]k} + F_{ij}H_{kl}$. However, because $H_{k[i}E_{j]l} - H_{l[i}E_{j]k}$ is trace-free and has symmetries corresponding to the Young diagram of the partition $(22)$ it vanishes by Lemma \ref{weylcriterion}, and so 
\begin{align}\label{rijkl}
& R_{ijkl} = -2H_{k[i}E_{j]l} -2H_{l[i}E_{j]k} + RH_{l[i}H_{j]k} + F_{ij}H_{kl}.
\end{align}
The Ricci identity implies
\begin{align}\label{skewnablah}
2H_{k[i}E_{j]l} + 2H_{l[i}E_{j]k} =  -R_{ij(kl)} + F_{ij}H_{kl} = \nabla_{[i}\nabla_{j]}H_{kl} = \nabla_{[i}\bt_{j]kl}.
\end{align}
There hold
\begin{align}
\label{gd1}
\begin{split}
&D_{i}\ga_{j} = \nabla_{i}\ga_{j} - \tfrac{1}{2}\bt_{ij}\,^{p}\ga_{p} - 2\ga_{i}\ga_{j} + H_{ij}\ga_{p}\ga^{p}, \qquad  
D^{p}\ga_{p}  = \nabla^{p}\ga_{p} = \ga_{p}\,^{p},
\end{split}
\end{align}
\begin{align}
  \label{dga3}
\begin{split}&D_{i}\bt_{jk}\,^{l} = \\ &\nabla_{i}\bt_{jk}\,^{l} - \tfrac{1}{2}\bt_{ij}\,^{p}\bt_{kp}\,^{l} -\tfrac{1}{2}\bt \delta_{[i}\,^{l}H_{j]k}+ \delta_{i}\,^{l}\ga_{p}\bt_{jk}\,^{p} + 2H_{i(j}\bt_{k)p}\,^{l}\ga^{p} - 3\ga_{(i}\bt_{jk)}\,^{l} - \bt_{ijk}\ga^{l}.
\end{split}
\end{align}
Tracing \eqref{skewnablah} in $il$, relabeling, and substituting in \eqref{qrrelated} yields the first equality of
\begin{align}\label{ddivbt}
4E_{ij} = \nabla_{p}\bt_{ij}\,^{p} - \tfrac{1}{2}\bt H_{ij} = D_{p}\bt_{ij}\,^{p} - 2\ga_{p}\bt_{ij}\,^{p},
\end{align}
while the second follows from the first and tracing \eqref{dga3}. Equation \eqref{ddivbt} plays an important role in deriving consequences of the Einstein equations; see Lemma \ref{eijvanishlemma} below. Recall that the curvature of $D$ is written $\sR_{ijk}\,^{l} = \sR_{h}\delta_{[i}\,^{l}h_{j]k}$. Calculating $\sR_{ijkl} - R_{ijkl}$ using \eqref{dnabladiff} and simplifying yields
\begin{align}
\label{confscal}&\sR_{h} = \uR_{h} + \tfrac{1}{4}|\bt|^{2}_{h}  -2|\det h|^{-1/2}\nabla^{p}\ga_{p} = \uR_{h} + \tfrac{1}{4}|\bt|^{2}_{h} + 2\dad_{h}\ga.
\end{align}
Equation \eqref{confscal} will be used repeatedly throughout the remainder of the paper.

\subsection{}\label{fdhsection}
For a Riemannian signature AH structure $(\en, [h])$ and a representative $h \in [h]$ with associated K\"ahler form $\om_{ij}$ define $\fd_{h} \in \cinf(M)$ by $2F_{ij} = \fd_{h}\om_{ij}$. Equivalently, $\fd_{h} = 2 \star F$. Note that $2|F|_{h}^{2}  = \fd_{h}^{2}$, and $\lie_{J\ga^{\sharp}}\om = d(\imt(J\ga^{\sharp})\om) = -d\ga = F$.  Let $\fd = \fd_{h}|\det h|^{1/2}$, which does not depend on the choice of $h \in [h]$. Decomposing \eqref{ricdecompose} by parts and substituting \eqref{confscal} yields
\begin{align}\label{rparts}
\begin{split}
R_{\al\be} & = -2E_{\al\be},\\
R_{\al\bar{\be}} & = \tfrac{1}{2}(R - \j \fd)H_{\al\bar{\be}} = \tfrac{1}{2}(\uR_{h} - \j \fd_{h})h_{\al\bar{\be}} = \tfrac{1}{2}(\sR_{h} + \tfrac{1}{4}|\bt|_{h}^{2} + 2\dad_{h}\ga - \j \fd_{h})h_{\al\bar{\be}}.
\end{split}
\end{align}
Because of \eqref{rparts} it makes sense to refer to $\cscal \defeq R - \j \fd$ as the \textbf{complex weighted scalar curvature} of the AH structure $(\en, [h])$. There will be written $\cscal_{h} = \uR_{h} - \j \fd_{h}$.

Since on a Riemann surface $M$ a $(0, 2)$ form vanishes, for an operator $\hdelbar$ on a complex vector bundle $E$ to be a holomorphic structure it suffices that it satisfy the Leibniz rule. It follows that if $(M, J)$ is a Riemann surface the most general holomorphic structure on $\cano^{k}$ has the form $\hdelbar = \delbar + \si^{(0, 1)}$ for an arbitrary one-form $\si_{i}$ on $M$ and the holomorphic structure $\delbar$ induced on $\cano^{k}$ by $J$. From \eqref{rparts} it is immediate that $R_{\al\be} = 0$ if and only if $E_{ij} = 0$. These conditions mean that the Ricci curvature of $(\en, [h])$ has type $(1,1)$. Equation \eqref{ddivbt} has the following nice interpretation.
 
\begin{lemma}\label{eijvanishlemma}
For an AH structure $(\en, [h])$ with cubic torsion $\bt_{ij}\,^{k}$ on an oriented surface the following are equivalent:
\begin{list}{(\arabic{enumi}).}{\usecounter{enumi}}
\renewcommand{\theenumi}{Step \arabic{enumi}}
\renewcommand{\theenumi}{\arabic{enumi}}
\renewcommand{\labelenumi}{\textbf{Level \theenumi}.-}
\item\label{er1} The curvature $R_{ij}$ has type $(1,1)$.
\item\label{er3} For any $h \in [h]$ the $(3, 0)$ part of the tensor $B_{ijk} \defeq \bt_{ij}\,^{p}h_{kp}$ is $\hdelbar$-holomorphic for the holomorphic structure $\hdelbar \defeq \delbar - 2\ga^{(0, 1)}$ on $\cano^{3}$, in which $\delbar$ is the holomorphic structure induced on $\cano^{3}$ by the conformal structure $[h]$ and the given orientation, and $4\ga_{i} = |\det h|^{-1}\nabla_{i}|\det h|$.
\end{list}
\end{lemma}
\begin{proof}
For $\tilde{h}_{ij} = fh_{ij}$ let $\tdelbar = \delbar - 2\tilde{\ga}^{(0,1)}$. Since $\tilde{\ga}_{i} = \ga_{i} + \tfrac{1}{2}d\log{f}_{i}$, $\tdelbar = \hdelbar - \delbar \log{f}$. It follows that $\tilde{B}_{ijk} = fB_{ijk}$ satisfies $\tdelbar \tilde{B}^{(3, 0)} = f\tdelbar B^{(3,0)} + \delbar f \tensor B^{(3, 0)} = f\hdelbar B^{(3,0)}$, so that $\tilde{B}^{(3,0)}$ is $\tdelbar$-holomorphic if and only if $B^{(3,0)}$ is $\hdelbar$-holomorphic. This shows that \eqref{er3} has sense. For any $h \in [h]$ it follows from \eqref{delb} that 
\begin{align*}
\begin{split}\hdelbar B^{(3,0)} &= \delbar B^{(3,0)} - 2\ga^{(0,1)}\tensor B^{(3,0)} \\
&= \overline{h^{(1,1)}} \tensor \div_{h}(B)^{(2,0)} - 2\overline{h^{(1,1)}} \tensor (\imt(\ga^{\sharp})B)^{(2,0)} = 4\overline{h^{(1,1)}}\tensor E^{(2,0)}.
\end{split}
\end{align*}
the last equality by \eqref{ddivbt}. Hence $B^{(3,0)}$ is $\hdelbar$-holomorphic if and only if $E_{ij} \equiv 0$.
\end{proof}

For $h \in [h]$ the Chern connection $\hnabla$ on $\cano^{3}$ determined by the Hermitian structure induced on $\cano^{3}$ by $h^{(1,1)}$ and the holomorphic structure $\hdelbar$ is by definition the unique connection on $\cano^{k}$ such that $\hnabla^{0, 1} = \hdelbar$ and for which the induced Hermitian structure is parallel. In terms of the connection induced on $\cano^{3}$ by the Levi-Civita connection $D$ of $h$, $\hnabla$ is expressible by 
\begin{align}
\hnabla = D - 2\j \ga \circ J = D + 2\j \star \ga.
\end{align}
It follows that the difference of the curvatures of $\hnabla$ and $D$ is $2\j d\star \ga = - 2\j \dad_{h}\ga \,\om$.

\subsection{Curvature of the conjugate AH structure}\label{curvatureahpencilsection}
It is convenient to define $\nbt = \bt_{abc}\bt^{abc}$. Evidently $\nbt H_{ij} = |\bt|_{h}^{2}h_{ij}$ for any $h \in [h]$. By Lemma \ref{twodablemma}, $2\bt_{ip}\,^{q}\bt_{jq}\,^{p} = \nbt H_{ij}$. Since $2\bt_{k[i}\,^{p}\bt_{j]lp} - \nbt H_{l[i}H_{j]k}$ has the algebraic symmetries of a Riemannian curvature tensor and is completely trace-free, Lemma \ref{weylcriterion} implies that it is identically zero. That is $2\bt_{k[i}\,^{p}\bt_{j]lp} = \nbt H_{l[i}H_{j]k}$.

Let $(\ben, [h])$ be the AH structure conjugate to $(\en, [h])$. Decorate with a $\,\bar{\,}\,$ the tensors derived from the curvature $\bar{R}_{ijkl}$ of $(\ben, [h])$. By definition 
\begin{align*}
\bar{R}_{ijk}\,^{l} - R_{ijk}\,^{l} = 2\nabla_{[i}\bt_{j]k}\,^{l} - 2\bt_{k[i}\,^{p}\bt_{j]p}\,^{l} = 2\nabla_{[i}\bt_{j]k}\,^{l} - \nbt \delta_{[i}\,^{l}H_{j]k}, 
\end{align*}
so, lowering indices using $2H_{lp}\nabla_{[i}\bt_{j]k}\,^{p} = 2\nabla_{[i}\bt_{j]kl} - 2\bt_{lpl[i}\bt_{j]k}\,^{p} = 2\nabla_{[i}\bt_{j]kl}  + \bt H_{l[i}H_{j]k}$, and using \eqref{skewnablah} there holds
\begin{align}
\label{conjugateblaschkecurvature2}
\begin{split}
\bar{R}_{ijkl} &= R_{ijkl} + 4H_{l[i}E_{j]k} + 4H_{k[i}E_{j]l}.
\end{split}
\end{align} 
Tracing \eqref{conjugateblaschkecurvature2} and substituting \eqref{ricdecompose} and $\bar{F}_{ij} = F_{ij}$ into the result yields $-4\bar{E}_{ij} + \bar{R}H_{ij} = 4E_{ij} + RH_{ij}$, which when traced shows $\bar{R} = R$ and $\bar{E}_{ij} = -E_{ij}$. 

Tensors such as $R$ (resp. $E_{ij}$) unchanged under conjugacy (multiplied by $-1$ under conjugacy) will be called \textbf{self-conjugate} (resp. \textbf{anti-self-conjugate}). Classes of AH structures defined by some condition on the curvatures preserved under conjugacy seem to be of particular interest. By \eqref{rijkl} and \eqref{conjugateblaschkecurvature2}, the self-conjugate and anti-self-conjugate parts of the curvature tensor are
\begin{align}
\begin{split}
\tfrac{1}{2}(R_{ijkl} + \bar{R}_{ijkl}) & = RH_{l[i}H_{j]k} + F_{ij}H_{kl}, \qquad
\tfrac{1}{2}(R_{ijkl} - \bar{R}_{ijkl})  = - 2H_{l[i}E_{j]k} - 2H_{k[i}E_{j]l} .
\end{split}
\end{align}
It follows that $(\en, [h])$ has self-conjugate curvature if and only if $E_{ij} = 0$. Also, for a Weyl structure the curvature is simply $R_{ijkl} = RH_{l[i}H_{j]k} + F_{ij}H_{kl}$. However, it will be evident later that on a compact orientable surface of genus at least $2$ the class of AH structures with self-conjugate curvature tensor is considerably larger than the class of Weyl structures.

\subsection{}
By Lemma \ref{weylcriterion}, the projective Cotton tensor $C_{ijk}$ satisfies $C_{ijk} = 2C_{[i}H_{j]k}$ in which $C_{i} \defeq C_{ip}\,^{p}$. From the definition of $C_{ijk}$ there results
\begin{align}\label{projcottontrace}
\begin{split}
&C_{i}   = - \tfrac{1}{2}\nabla_{i}R - \tfrac{1}{3}\nabla^{p}F_{ip} -2 \nabla^{p}E_{ip} +2\bt_{i}\,^{pq}E_{pq}
\end{split}
\end{align}
From \eqref{projcottontrace} and the easily verified identities
\begin{align*}
&\bnabla_{i}\bar{R} = \nabla_{i}R, && \bnabla^{p}\bar{F}_{ip} = \nabla^{p}F_{ip}, && \bnabla^{p}\bar{E}_{ip} = -\nabla^{p}E_{ip} +\bt_{i}\,^{pq}E_{pq},
\end{align*}
there result
\begin{align}\label{projcotton}
\begin{split}
&\bar{C}_{i} = C_{i} + 4\nabla^{p}E_{ip} -2\bt_{i}\,^{pq}E_{pq},\\
& C_{i} + \bar{C}_{i}  = - \nabla_{i}R - \tfrac{2}{3}\nabla^{p}F_{ip} + 2\bt_{i}\,^{pq}E_{pq}, \\
& C_{i} - \bar{C}_{i}  = -4\nabla^{p}E_{ip} + 2\bt_{i}\,^{pq}E_{pq} =  -4D^{p}E_{ip},
\end{split}
\end{align}
in the last equality of which $D$ is the Levi-Civita connection of any $h \in [h]$. 

\begin{lemma}\label{selfconjugatecottonlemma}
A Riemannian signature AH structure $(\en, [h])$ on an oriented surface $M$ has self-conjugate projective Cotton tensor if and only if $E_{ij}$ is the real part of a holomorphic quadratic differential with respect to the complex structure determined by $[h]$. In particular, if $M$ is a sphere, then $(\en, [h])$ has self-conjugate projective Cotton tensor if and only if $E_{ij} = 0$.
\end{lemma}

\begin{proof} 
Since $E_{ij}$ is trace-free, Lemma \ref{kdifferentialslemma} implies that $E_{ij}$ is the real part of a holomorphic quadratic differential if and only if $\div_{h}(E)_{i} = 0$, and by the last equality of \eqref{projcotton} this holds if and only if $(\en, [h])$ has self-conjugate projective Cotton tensor. 
\end{proof}

\begin{theorem}\label{nonnegselfconjtheorem}
If a Riemannian signature AH structure $(\en, [h])$ on a compact oriented surface has self-conjugate projective Cotton tensor and non-negative weighted scalar curvature, then $E_{ij} = 0$.
\end{theorem}
\begin{proof}
By Lemma \ref{selfconjugatecottonlemma}, $E_{ij}$ is the real part of a holomorphic quadratic differential. By \eqref{keromlap} of Lemma \ref{flatmetriclemma} and \eqref{confscal} for a Gauduchon metric $h \in [h]$ there holds $\lap_{h}|E|_{h}^{2} \geq 2|DE|_{h}^{2} + \tfrac{1}{2}|\bt|_{h}^{2}|E|_{h}^{2} \geq 0$, and the maximum principle then forces $E_{ij} = 0$.
\end{proof}

\subsection{}
Since $R$ and $\nbt $ are $1$-densities their integrals are defined if $M$ is compact. The $L^{2}$-norm $||\bt||_{h}^{2}$ does not depend on the choice of representative $h \in [h]$ and equals $\int_{M}\nbt = \int_{M}|\bt|^{2}_{h}\,d\vol_{h}$.

\begin{theorem}\label{gbtheorem}
If $(\en, [h])$ is a Riemannian AH stucture on a compact, orientable surface $M$, then the Euler characteristic $\chi(M)$ satisfies $4\pi\chi(M) \geq \int_{M}R$, with equality if and only if $(\en, [h])$ is Weyl. In particular, 
\begin{enumerate}
\item If $\int_{M}R \geq 0$, then either $M$ is a sphere, or $M$ is a torus and $(\en, [h])$ is Weyl.
\item If $M$ has genus at least one and $(\en, [h])$ is not Weyl, then $\int_{M} R < 0$.
\end{enumerate}
\end{theorem}

\begin{proof}
By the Gau\ss-Bonnet Theorem, for any $h \in [h]$, integrating \eqref{confscal} yields
\begin{align}\label{gbinter1}
4\pi\chi(M) = \int_{M}\sR_{h}\,d\vol_{h} = \tfrac{1}{4}||\bt||_{h}^{2} + \int_{M} R\geq \int_{M}R .
\end{align}
Equality holds in \eqref{gbinter1} if and only if $\bt_{ij}\,^{k} = 0$. If $\int_{M}R \geq 0$ the Euler characteristic $\chi(M)$ of $M$ must be non-negative, so $M$ must the sphere or the torus. If $M$ is a torus, the Euler characteristic is $0$ and so \eqref{gbinter1} forces $||\bt||_{h}^{2} = 0$, so that $(\en, [h])$ is a Weyl structure. If $M$ has genus at least one and $(\en, [h])$ is not Weyl then $4\pi\chi(M) - \tfrac{1}{4}||\bt||^{2}_{h} < 0$, showing the last claim.
\end{proof}

\section{Einstein equations}\label{einsteinsection}
In this section the Einstein equations for AH structures are defined and their most basic properties described.
\subsection{}
By definition of $E_{ij}$, \eqref{ddivbt}, and Lemma \ref{eijvanishlemma}, the following conditions on an AH structure on a surface having cubic torsion $\bt_{ij}\,^{k}$ are equivalent.
\begin{list}{(\arabic{enumi}).}{\usecounter{enumi}}
\renewcommand{\theenumi}{Step \arabic{enumi}}
\renewcommand{\theenumi}{\arabic{enumi}}
\renewcommand{\labelenumi}{\textbf{Level \theenumi}.-}
\item \label{ne1} The symmetric part of the Ricci tensor is trace free. That is $E_{ij} = 0$. 
\item \label{ne2} The Ricci tensor has type $(1,1)$.
\item \label{ne3} The curvature is self-conjugate.
\item \label{ne3b} For any $h \in [h]$ the $(3, 0)$ part of the tensor $B_{ijk} \defeq \bt_{ij}\,^{p}h_{kp}$ is $\hdelbar$-holomorphic for the holomorphic structure $\hdelbar \defeq \delbar - 2\ga^{(0, 1)}$, in which $\delbar$ is the holomorphic structure induced by the conformal structure $[h]$, and $4\ga_{i} = \nabla_{i}\log|\det h|$.
\end{list}
\begin{definition}
An AH structure on a surface is \textbf{naive Einstein} if it satisfies \eqref{ne1}-\eqref{ne3b}.
\end{definition}

The qualifier \textit{naive} is meant to reflect that, while the most obvious generalization of the usual metric Einstein condition is simply to require the vanishing of the symmetric trace-free Ricci tensor, such an approach turns out to be inadequate.

By Lemma \ref{selfconjugatecottonlemma} a Riemannian signature AH structure $(\en, [h])$ on the two-sphere has self-conjugate projective Cotton tensor if and only if it is naive Einstein. Similarly, by Theorem \ref{nonnegselfconjtheorem}, a Riemannian signature AH structure $(\en, [h])$ on a compact, oriented surface is naive Einstein if it has non-negative weighted scalar curvature.

\begin{definition}
An AH structure $(\en, [h])$ is \textbf{Einstein} if it is naive Einstein and satisfies
\begin{align}\label{conservationcondition}
\nabla_{i}R + 2\nabla^{p}F_{ip} = 0.
\end{align}
\end{definition}
\noindent
The condition \eqref{conservationcondition} will be referred to as the \textbf{conservation condition}. 

Let $h \in [h]$. Using $\nabla_{i}R = D_{i}R + 2\ga_{i}R$ to expand \eqref{conservationcondition}, and using \eqref{gd1}, \eqref{confscal}, and the Ricci identity yields
\begin{align}\label{conserve} 
\begin{split}
|\det h|^{-1/2}&(\nabla_{i}R + 2\nabla^{p}F_{ip})  = D_{i}\uR_{h}  + 2\ga_{i}\uR_{h} + 2h^{pq}D_{p}F_{iq} + 4\ga^{\sharp\,p}F_{ip}\\
 & = D_{i}\uR_{h}  + 2\ga_{i}\uR_{h} - 2\dad_{h}d\ga_{i} - 4\ga^{\sharp\,p}F_{pi} \\
& = D_{i}\left(\sR_{h} - \tfrac{1}{4}|\bt|_{h}^{2}\right)  - \tfrac{1}{2}\ga_{i}|\bt|_{h}^{2} + \sR_{h}\ga_{i} + 2\lap_{h} \ga_{i} -4 \ga_{i}\dad_{h}\ga -4\ga^{\sharp\,p}F_{pi}\\
& = D_{i}\left(\sR_{h} - \tfrac{1}{4}|\bt|_{h}^{2}\right)  - \tfrac{1}{2}\ga_{i}|\bt|_{h}^{2} - 2(\hodge_{h}- \sR_{h}) \ga_{i}  -4 \ga_{i}\dad_{h}\ga  -4\ga^{\sharp\,p}F_{pi}\\
& = D_{i}\left(\sR_{h} - \tfrac{1}{4}|\bt|_{h}^{2} - 4|\ga|_{h}^{2}\right)  - \tfrac{1}{2}\ga_{i}|\bt|_{h}^{2} - 2(\hodge_{h} - \sR_{h}) \ga_{i}  + 4\ga^{\sharp\,p}(\lie_{\ga^{\sharp}}h)_{ip} -4 \ga_{i}\dad_{h}\ga.
\end{split}
\end{align}
While the only explicit direct use of \eqref{conservationcondition}, then in the rewritten form \eqref{conserve}, is made in the proof of Theorem \ref{classtheorem}, its role is fundamental. By \eqref{ddivbt} a Weyl structure on a surface is automatically naive Einstein. In \cite{Calderbank-mobius}, Calderbank showed that taking \eqref{conservationcondition} as the definition of Einstein Weyl yields a nice theory, specializing that in higher dimensions. By definition the Einstein AH equations restrict to Calderbank's Einstein Weyl equations, and Calderbank's definition provided essential motivation for the general case. 

In dimensions $n > 2$ the conservation condition has a more general definition, but it follows from the differential Bianchi identity that a naive Einstein AH structure with self-conjugate curvature is Einstein, satisfying the analogue of \eqref{conservationcondition} with $n$ in place of $2$; see Lemma $4.2$ of \cite{Fox-ahs}.

\subsection{}

Lemma \ref{parallelexactlemma} follows from \eqref{conservationcondition}, Lemma \ref{fcoclosedlemma}, and the discussion at the end of section \ref{faradaysection}.
\begin{lemma}\label{parallelexactlemma} 
A Riemannian signature Einstein AH structure on a surface is closed if and only if it has parallel scalar curvature, in which case either it is proper and exact or it has vanishing weighted scalar curvature. 
\end{lemma}

\begin{lemma}\label{2deinsteinlemma}
For an Einstein AH structure $(\en, [h])$ on a surface, any one of the following statements implies the other two.
\begin{list}{(\arabic{enumi}).}{\usecounter{enumi}}
\renewcommand{\theenumi}{Step \arabic{enumi}}
\renewcommand{\theenumi}{\arabic{enumi}}
\renewcommand{\labelenumi}{\textbf{Level \theenumi}.-}
\item\label{ep1} $(\en, [h])$ is projectively flat.
\item\label{ep2} $(\en, [h])$ is conjugate projectively flat.
\item\label{ep3} The weighted scalar curvature is parallel.
\end{list}
In particular if an Einstein $(\en, [h])$ either is proper or has vanishing scalar curvature then it is projectively flat and conjugate projectively flat. 
\end{lemma}
\begin{proof}
The first claim is immediate from \eqref{projcotton} and the Einstein equations. The last claim follows because a proper Einstein AH structure is exact so has parallel scalar curvature. 
\end{proof}

\subsection{}
The example of affine hypersurfaces gives the primary motivation for the definition of the Einstein equations for AH structures. 
\begin{theorem}\label{sphereeinsteintheorem}
For a non-degenerate positively co-oriented hypersurface immersion into flat three-dimensional affine space the following are equivalent:
\begin{enumerate}
\item The image of the immersion is an affine hypersphere.
\item The AH structure induced on the hypersurface is Einstein.
\end{enumerate}
Moreover, $(1)$ and $(2)$ imply
\begin{enumerate}
\setcounter{enumi}{2}
\item The induced AH structure is projectively flat. 
\end{enumerate}
\end{theorem}
\begin{proof}
 The equivalence of $(1)$ and $(2)$ is proved by the same argument as in the $n > 2$ case, which can be found as Theorem $4.6$ of \cite{Fox-ahs}. The AH structure induced on an affine hypersurface is conjugate projectively flat because the conjugate AH structure is that induced via the conormal Gau\ss\, map from the flat projective structure on oriented projective space, so Lemma \ref{2deinsteinlemma} implies that that the induced AH structure is projectively flat as well. Alternatively, the AH structure induced on a non-degenerate hypersurface in affine space is always closed, and when it is Einstein, this implies the scalar curvature is parallel, so by Lemma \ref{2deinsteinlemma}, that it is projectively flat. 
\end{proof}

\subsection{}
Recall the definitions of $\fd_{h}$ and $\cscal$ from section \ref{fdhsection}. Note that $\cscal$ is a smooth section of the complexification $|\Det \ctm| \tensor_{\rea}\com$ of the line bundle of $1$-densities and that the $(0,1)$ part of the aligned representative $\nabla$ of an AH structure induces a holomorphic structure on $|\Det \ctm| \tensor_{\rea}\com$.

\begin{lemma}\label{einsteinhololemma}
A naive Einstein Riemannian signature AH structure $(\en, [h])$ on an oriented surface $M$ is Einstein if and only if $\cscal$ is a holomorphic section of $|\Det \ctm|\tensor_{\rea}\com$ with respect to the holomorphic structure $\nabla^{0,1}$.
\end{lemma}
\begin{proof}
The claim to be verified is simply $\nabla^{0,1}\cscal = \nabla^{0,1}(R - \j \fd) = 0$. First it is shown that $(\en, [h])$ is Einstein if and only if for any $h \in [h]$ the complex valued function $\uR_{h} - \j \fd_{h}$ is $\hdelbar$-holomorphic for the holomorphic structure $\hdelbar = \delbar + 2\ga^{(0,1)}$ on the trivial line bundle $M \times \com$. 
It is claimed that the conservation condition \eqref{conservationcondition} is equivalent to
\begin{align}\label{holoconservationcondition}
\begin{split}
0 & = \hdelbar(\uR_{h} - \j \fd_{h}) = \delbar(\uR_{h} - \j \fd_{h}) + 2(\uR_{h} - \j \fd_{h})\ga^{(0,1)}\\
& = \delbar(\uR_{h} - 4|\ga|_{h}^{2} - \j \fd_{h}) + 2\uR_{h} \ga^{(0,1)} + 4(\imt(\ga^{\sharp})\lie_{\ga^{\sharp}}h)^{(0,1)}.
\end{split}
\end{align}
Rewriting the first equality of \eqref{conserve} yields
\begin{align}\label{preconserve} 
\begin{split}
|\det h|^{-1/2}&(\nabla_{i}R + 2\nabla^{p}F_{ip}) 
 = d\uR_{h\, i} + 2\ga_{i}\uR_{h} - (\star d\fd_{h})_{i} - 2\fd_{h}(\star \ga)_{i}.
\end{split}
\end{align}
The $(0, 1)$ part of this last expression is
\begin{align}\label{delbarcscal}
\delbar( \uR_{h} - \j \fd_{h}) + 2(\uR_{h} - \j \fd_{h})\ga^{(0, 1)}= \hdelbar( \uR_{h} - \j \fd_{h}),
\end{align}
from which the first claim and the first equality of \eqref{holoconservationcondition} are evident. Taking the $(0,1)$ part of
\begin{align*}
\begin{split}
d_{i} |\ga|_{h}^{2} &= 2\ga^{\sharp\,p}D_{i}\ga_{p} = 2\ga^{\sharp\,p}D_{(i}\ga_{p)} + 2\ga^{\sharp\,p}D_{[i}\ga_{p]}\\
& = \ga^{\sharp\,p}(\lie_{\ga^{\sharp}}h)_{ip} - \tfrac{1}{2}\ga^{\sharp\,p}\fd_{h}\om_{ip} = \ga^{\sharp\,p}(\lie_{\ga^{\sharp}}h)_{ip} + \tfrac{1}{2}\fd_{h}(\star \ga)_{i},
\end{split}
\end{align*}
yields $2\delbar |\ga|_{h}^{2} = \j \fd_{h} \ga^{(0, 1)}+ 2(\imt(\ga^{\sharp})\lie_{\ga^{\sharp}}h)^{(0,1)}$. Substituting this into what comes before yields the second equality of \eqref{holoconservationcondition}. By definition of $\ga_{i}$, for any $h \in [h]$ there holds
\begin{align*}
\begin{split}
\nabla^{0,1}(R - \j \fd) & = \delbar(\uR_{h} - \j \fd_{h})|\det h|^{1/2}  + (\uR_{h} - \j \fd_{h})\nabla^{0, 1}|\det h|^{1/2}\\& = \left(\delbar(\uR_{h} - \j \fd_{h})+ 2(\uR_{h} - \j \fd_{h})\ga^{(0, 1)}\right)|\det h|^{1/2},
\end{split}
\end{align*}
from which the claim follows.
\end{proof}

By Lemma \ref{einsteinhololemma}, a Riemannian AH structure on an oriented surface is Einstein if and only if its Ricci curvature has type $(1,1)$ and its complex weighted scalar curvature is holomorphic.

\begin{lemma}\label{squarehololemma}
Let $(\en, [h])$ be a Riemannian signature Einstein AH structure on an oriented surface $M$. For each Gauduchon metric $h \in [h]$ the square norm $|\cscal_{h}|_{h}^{2} = \uR_{h}^{2} + \fd_{h}^{2}$ of the complex scalar curvature is $h$-harmonic. If $M$ is moreover compact then $  \uR_{h}^{2} + \fd_{h}^{2}$ is constant. If $H^{1}(M; \rea) = \{0\}$ there is a $v \in \cinf(M)$ such that $e^{2\j v}\cscal_{h}$ is holomorphic on $M$.
\end{lemma}

\begin{proof}
Denote by the same notations the lifts to the universal cover $\tilde{M}$ of $M$ of $(\en, [h])$ and all the associated tensors, functions, etc. Let $\ga$ be the Faraday primitive of $h$. By assumption $\star \ga$ is closed, so on $\tilde{M}$ there is a smooth function $v$ such that $\star \ga = -dv$, or, what is equivalent, $\ga = \star dv$. By \eqref{delbarcscal} of Lemma \ref{einsteinhololemma} the function $b = e^{2\j v}(\uR_{h} - \j \fd_{h})$ is holomorphic on $\tilde{M}$, for by construction $\delbar v = -\j \ga^{(0,1)}$. Hence $|b|^{2}_{h} = \uR_{h}^{2} + \fd_{h}^{2}$ is harmonic on $\tilde{M}$, and so the corresponding $\uR_{h}^{2} + \fd_{h}^{2}$ on $M$ must be harmonic as well. If $M$ is moreover compact, then $\uR_{h}^{2} + \fd_{h}^{2}$ is constant by the maximum principle. If $H^{1}(M;\rea) = \{0\}$, then $v$ could be taken initially to be defined on $M$, and so $e^{2\j v}\cscal_{h}$ is holomorphic on $M$.
\end{proof}

\section{Classification of Einstein AH structures by scalar curvature and genus}\label{scalarsection}

\subsection{}
Theorem \ref{classtheorem} is the key technical result for the description of Einstein AH structures on compact orientable surfaces. It generalizes the result for Einstein Weyl structures proved in Theorem $3.7$ of \cite{Calderbank-mobius}. The Killing property of the Gauduchon metric dual of the associated Faraday primitive generalizes to Einstein AH structures a property of the Gauduchon gauge for Einstein Weyl structures first shown in Theorem $2.2$ of \cite{Tod-compact}.

\begin{theorem}\label{classtheorem}
A Riemannian signature AH structure $(\en, [h])$ on a compact orientable surface is Einstein if and only if for a Gauduchon metric $h \in [h]$ with Levi-Civita connection $D$ and Faraday primitive $\ga_{i}$ there are satisfied the equations
\begin{align}\label{confein1} &D_{p}\bt_{ij}\,^{p} = 0, \qquad |\bt|_{h}^{2}\ga_{i} = 0,\qquad (\lie_{\ga^{\sharp}}h)_{ij} = 2D_{(i}\ga_{j)} = 0,\\
\label{confein2} &D_{i}(\sR_{h} - \tfrac{1}{4}|\bt|_{h}^{2} - 4|\ga|_{h}^{2}) = D_{i}(\uR_{h} -4|\ga|_{h}^{2}) = 0.
\end{align}
Moreover, each of $\uR_{h}$, $|\bt|_{h}^{2}$, and $\fd_{h}$ is constant along the flow of $\ga^{\sharp}$.

Conversely, if on a manifold of dimension at least $2$ (not necessarily compact) there are a Riemannian metric $h$ (with Levi-Civita connection $D$), an $h$-Killing field $X^{i}$, and a completely symmetric, completely $h$-trace free tensor $B_{ijk} = B_{(ijk)}$, such that $\ga_{i}\defeq X^{i}h_{pi}$ and $\bt_{ij}\,^{k} \defeq h^{kp}B_{ijp}$ solve \eqref{confein1}-\eqref{confein2}, then $\nabla \defeq D - \tfrac{1}{2}\bt_{ij}\,^{k} - 2\ga_{(i}\delta_{j)}\,^{k} + h_{ij}\ga^{k}$ is the aligned representative of an Einstein AH structure $(\en, [h])$ for which $h$ is a distinguished metric. 
\end{theorem}

\begin{proof}
Let $h \in [h]$ and in this proof raise and lower indices with $h_{ij}$ and $h^{ij}$. Recall that $\lie_{\ga^{\sharp}}h_{ij} = 2D_{(i}\ga_{j)}$. The Ricci identity gives 
\begin{align}\label{bch0}
\begin{split}
2D_{i}D^{p}\ga_{p}&= 2D^{p}D_{i}\ga_{p} - \sR_{h}\ga_{i} = D^{p}(\lie_{\ga^{\sharp}}h)_{ip} - D^{p}F_{ip} - \sR_{h}\ga_{i}.
\end{split}
\end{align}
Let $||\dum||_{h}$ denote the $L^{2}$ norm on tensors with respect to the $h$-volume measure $d\vol_{h}$. Contracting \eqref{bch0} with $\ga^{\sharp}$ and integrating the result gives
\begin{align}\label{bch1}
2||\lie_{\ga^{\sharp}}h||_{h}^{2} - 8||\dad_{h}\ga||^{2}_{h} - ||\fd_{h}||^{2}_{h} + 4\int_{M}\sR_{h}|\ga|_{h}^{2}\,d\vol_{h} = 0. 
\end{align}
Contracting the second line of \eqref{conserve} with $\ga^{\sharp\,i}$, integrating by parts, and substituting in \eqref{confscal} yields
\begin{align}\label{einsteinf}
||\fd_{h}||^{2}_{h} - 2\int_{M}\uR_{h}\dad_{h}\ga \,d\vol_{h} = 4\int_{M}|\ga|_{h}^{2}\uR_{h}\,d\vol_{h} = 4\int_{M}\sR_{h}|\ga|_{h}^{2}\,d\vol_{h} - \int_{M}|\ga|_{h}^{2}|\bt|_{h}^{2}\,d\vol_{h}.
\end{align}
Substituting \eqref{einsteinf} into \eqref{bch1} and taking $h \in [h]$ to be Gauduchon yields
\begin{align}\label{bch2}
2||\lie_{\ga^{\sharp}}h||_{h}^{2} + \int_{M}|\ga|_{h}^{2}|\bt|_{h}^{2}\,d\vol_{h} = 8||\dad_{h}\ga||^{2}_{h} + 2\int_{M}\uR_{h}\dad_{h}\ga \,d\vol_{h} = 0.
\end{align}
Equation \eqref{bch2} implies the first two equalities of \eqref{confein1}. By Lemma \ref{twodablemma} that $|\ga|^{2}_{h}|\bt|^{2}_{h} = 0$ is equivalent to $\ga_{p}\bt_{ij}\,^{p} = 0$. By \eqref{ddivbt} this implies $D_{p}\bt_{ij}\,^{p} = 0$. Wherever $\ga_{i}$ is not zero, the one-form $|\bt|^{2}_{h}\ga_{i}$ is $h$-orthogonal to the linearly indepent one-forms $\ga_{i}$ and $J_{i}\,^{p}\ga_{p}$, so vanishes identically. This completes the proof of \eqref{confein1}. The first equality in \eqref{confein2} is true by \eqref{confscal}. Because $D_{(i}\ga_{j)} = 0$ and $\dad\ga = 0$ there holds $\hodge \ga_{i} = \dad d\ga_{i} = \sR_{h} \ga_{i}$. Substituting the preceeding observations into the last line of \eqref{conserve} gives \eqref{confein2}. By \eqref{confein1} and \eqref{confein2}, $\ga^{\sharp\, i}D_{i}\uR_{h} = 8\ga^{\sharp\,p}\ga^{\sharp\,q}D_{(p}\ga_{q)} = 0$, showing $d\uR_{h}(\ga^{\sharp}) = 0$. Since $\ga^{\sharp}$ is $h$-Killing, there holds $\ga^{\sharp\,p}D_{p}\sR_{h} = 0$, and with $d\uR_{h}(\ga^{\sharp}) = 0$ and \eqref{confscal} this shows $\ga^{\sharp\,i}D_{i}|\bt|^{2}_{h} = 0$. As $\lie_{\ga^{\sharp}}\om = d\imt(\ga^{\sharp})\om = d\star \ga = 0$ there holds $2\lie_{\ga^{\sharp}}F = d\fd_{h}(\ga^{\sharp})\om$. Since $D_{(i}\ga_{j)} = 0$ there holds $d|\ga|^{2}_{h} = 2\ga^{\sharp\,p}D_{i}\ga_{p} = \ga^{\sharp\,p}d\ga_{ip} = \ga^{\sharp\,p}F_{pi}$. Hence $\lie_{\ga^{\sharp}}F = d(\imt(\ga^{\sharp})F) = d(d|\ga|_{h}^{2}) = 0$, showing $2\lie_{\ga^{\sharp}}F = d\fd_{h}(\ga^{\sharp})\om = 0$.

If given $(h, X, B)$ as in the statement of the theorem, then it is straightforward to check that $(\en, [h])$ is an AH structure with cubic torsion $\bt_{ij}\,^{k}$, aligned representative $\nabla$, and Gauduchon metric $h$. The curvatures of $\nabla$ and $D$ are related as in \eqref{confscal}, and there hold \eqref{ddivbt} and \eqref{conserve}. Together \eqref{confein1} and \eqref{ddivbt} show $E_{ij} = 0$, and so show the naive Einstein equations. Finally, substituting $D_{i}|\ga|_{h}^{2} = \ga^{\sharp\,p}d\ga_{ip}$ and $D^{p}d\ga_{pi} = -\sR_{h}\ga_{i}$ into \eqref{conserve} yields $\nabla_{i}R + 2\nabla^{p}F_{ip} = 0$.
\end{proof}

\begin{remark}
In section \ref{naiveexample} it is explained how to construct naive Einstein AH structures which are not Einstein, and which illustrate that not all of the conclusions of Theorem \ref{classtheorem} hold for such structures.
\end{remark}

\begin{corollary}\label{holocorollary}
If $(\en, [h])$ is a Riemannian signature Einstein AH structure on a compact oriented surface $M$ of genus $g$ and $h \in [h]$ is a Gauduchon metric with associated Faraday primitive $\ga_{i}$ and cubic form $B_{ijk} \defeq \bt_{ij}\,^{p}h_{pk}$, then:
\begin{enumerate}
\item With respect to the complex structure determined by $[h]$ and the given orientation, $B^{(3,0)}$ and $\ga^{\sharp\,(1,0)}$ are holomorphic. Moreover, 
\begin{align}
\label{ne4b}&2\lap_{h}|\ga|_{h}^{2} = |d\ga|_{h}^{2} - 2\uR_{h}|\ga|_{h}^{2},& &\text{everywhere},\\
\label{ne4} &\lap_{h}\log |\ga|_{h}^{2} = -\uR_{h},& &\text{wherever $|\ga|_{h}^{2} > 0$},\\
\label{ne4c}&\lap_{h}|B|_{h}^{2} = 4|d|B||_{h}^{2} + 3\sR_{h} |B|_{h}^{2} = 4|d|B||_{h}^{2} + 3\uR_{h} |B|_{h}^{2} + \tfrac{3}{4}|B|_{h}^{4},& &\text{wherever $|B|_{h}^{2} > 0$},\\
\label{ne4d}&\lap_{h}\log |B|_{h}^{2} = 3\sR_{h} = 3\uR_{h} + \tfrac{3}{4}|B|_{h}^{2}& &\text{wherever $|B|_{h}^{2} > 0$}.
\end{align}
\item If $(\en, [h])$ is not Weyl then it is exact, while if $(\en, [h])$ is not exact, then it is Weyl. 
\item If $g \geq 2$ then $(\en, [h])$ is exact. In this case the metric $\sth_{ij} \defeq |B|_{h}^{2/3}h_{ij}$ defined on the open submanifold $M^{\ast}\defeq \{|B|^{2} \neq 0\}$ is flat.
\item If $g = 0$ then $(\en, [h])$ is Weyl.
\end{enumerate}
\end{corollary}
\begin{proof}
By the first and last equalities of \eqref{confein1} and Lemmas \ref{codazzilemma} and \ref{kdifferentialslemma}, $\ga^{\sharp\,(1,0)}$ and $B^{(3,0)}$ are holomorphic. Since by \eqref{confein1} of Theorem \ref{classtheorem}, $d\ga_{ij} = 2D_{i}\ga_{j}$, the second equation of \eqref{keromlap} of Lemma \ref{flatmetriclemma} applied to $\ga^{\sharp}$ (so with $k = -1$) reduces to \eqref{ne4b}, while by the second equation of \eqref{kato} there holds $2\lap_{h}\log |\ga| = -\sR_{h} = - \uR_{h} - \tfrac{1}{2}|\bt|_{h}^{2}$ wherever $\ga$ is non-zero. Multiplying through by $|\ga|^{2}$ and using that, by Theorem \ref{classtheorem}, $|\ga|^{2}|\bt|^{2} = 0$, this yields \eqref{ne4}. Similarly, equations \eqref{ne4c} and \eqref{ne4d} follow from \eqref{keromlap} and \eqref{confscal}.

Since each of $\ga^{\sharp\,(1,0)}$ and $B^{(3, 0)}$ is holomorphic, the zeroes of each are isolated if it is not identically zero, and hence the same is true for $\ga^{\sharp}$ and $B$. By \eqref{confein1}, there holds $|\ga|^{2}_{h}|B|^{2}_{h} = 0$. Because when $\ga^{\sharp}$ or $B$ is not zero its zeroes are isolated, this implies that if $\ga^{\sharp}$ is somewhere not zero then $B$ is identically zero, and if $B$ is not somewhere zero then $X$ is identically zero. Hence if $(\en, [h])$ is not Weyl it is exact and if $(\en, [h])$ is not exact it is Weyl. 

That $g \geq 2$ and $g = 0$ imply that $(\en, [h])$ is respectively exact and Weyl follows from the holomorphicity of $B^{(3,0)}$ and $\ga^{\sharp \,(1,0)}$ and Riemann Roch. By Riemann Roch, when $g > 1$, a holomorphic cubic differential has at most $6(g-1)$ zeroes, so $M^{\ast}$ is the complement of at most $6(g-1)$ points. By \eqref{ne4d} of Corollary \ref{holocorollary} there holds $|B|_{h}^{2/3}\sR_{\sth} = \sR_{h} - \lap_{h}\log |B|_{h}^{2/3} = 0$, so $\sth$ is flat on $M^{\ast}$
\end{proof}

\begin{lemma}\label{sphereconstantlemma}
If $(\en, [h])$ is a Riemannian signature Einstein AH structure on a compact orientable surface $M$, and $h \in [h]$ is any Gauduchon metric, then $\uR_{h}^{2} + \fd_{h}^{2}$ is equal to the constant $(\max_{M}\uR_{h})^{2}$.
\end{lemma}
\begin{proof}
If $(\en, [h])$ is exact then $\uR_{h}$ is constant by \eqref{conservationcondition}, and $\uR_{h}^{2} + \fd_{h}^{2} = \uR_{h}^{2} = (\max_{M}\uR_{h})^{2}$. If $(\en, [h])$ is not exact, then by Corollary \ref{holocorollary}, $M$ is a sphere or a torus. In this case, by Lemma \ref{squarehololemma}, $\be = \uR_{h}^{2} + \fd_{h}^{2}$ is constant. Since $|\ga|_{h}^{2}$ is not identically zero, it assumes somewhere a positive maximum; at such a point there holds $0 = 2d_{i}|\ga|^{2}_{h} = 4\ga^{\sharp\,p}D_{i}\ga_{p} = 2\ga^{\sharp\,p}F_{pi} = \fd_{h}(\star \ga)_{i}$, so at such a point $\fd_{h}$ vanishes, and $\be$ is equal to the value of $\uR_{h}^{2}$ at this point. Since by \eqref{confein2} of Theorem \ref{classtheorem}, $\uR_{h} - 4|\ga|^{2}_{h}$ is constant, the functions $\uR_{h}$ and $|\ga|_{h}^{2}$ assume their maximum values at the same points, so $\be = (\max_{M}\uR_{h})^{2}$. 
\end{proof}

\subsection{}
This section records a geometric interpretation of the integral curves of the Gauduchon metric dual of the Faraday primitive of the Gauduchon class. Recall that the \textbf{magnetic flow} determined on $\ctm$ by a pair $(g, \mu)$ comprising a metric $g$ and a closed two-form $\mu$ on $M$ is the Hamiltonian flow of the function $G(s) = \tfrac{1}{2}g_{\pi(s)}^{ij}s_{i}s_{j}$ on $\ctm$ with respect to the symplectic form $\Omega_{M} - \pi^{\ast}(\mu)$, where $\pi:\ctm \to M$ is the canonical projection and $\Omega_{M}$ is the canonical symplectic form. If $\mu$ is exact the magnetic flow is said to be \textbf{exact}.

The \textbf{magnetic geodesics} are the images in $M$ of the integral curves of this flow, which are the solutions of the equation $D_{\dot{\si}}\dot{\si} = A(\dot{\si})$ where $A$ is the tensor defined by $A_{i}\,^{p}g_{pj} = \mu_{ij}$, $D$ is the Levi-Civita connection of $g$, and $\si(t)$ is a curve in $M$. Along a magnetic geodesic $\si(t)$, the energy $|\dot{\si}|^{2}_{g}$ is constant. If $M$ is an oriented surface then $\mu_{ij} = f\om_{ij}$, where $\om_{ij}$ is the K\"ahler form associated to $g$, and a routine computation shows that the geodesic curvature $\ka_{g}(\si)$ of a magnetic geodesic $\si$ having energy $|\dot{\si}|^{2}_{g} = e^{2c}$ is $e^{-c}f(\si)$.

\begin{theorem}\label{magnetictheorem}
Let $(\en, [h])$ be a non-exact Einstein AH structure on a compact orientable surface $M$. Let $h \in [h]$ be a Gauduchon metric with Faraday primitive $\ga_{i}$. The integral curves of the vector field $\ga^{\sharp}$ are magnetic geodesics for the exact magnetic flow on $\ctm$ determined by $(h, \tfrac{1}{2}d\ga)$. The function $|\ga|_{h}^{2}$ is constant along a non-trivial integral curve of $\ga^{\sharp}$, and the geodesic curvature of such an integral curve is the restriction to the curve of $-\tfrac{1}{4}|\ga|_{h}^{-1}\fd_{h}$, which is constant along the integral curve. The integral curves of $\ga^{\sharp}$ passing through the points where $|\ga|_{h}^{2}$ attains its maximum value are $D$-geodesics; in particular, at least one integral curve of $\ga^{\sharp}$ is a $D$-geodesic. The non-trivial integral curves of $J(\ga^{\sharp}) = (\star \ga)^{\sharp}$ are projective geodesics of $D$.
\end{theorem}
\begin{proof}
By Corollary \ref{holocorollary}, $(\en, [h])$ is Weyl and $M$ is a torus or a sphere. From \eqref{confein1} there follows
\begin{align}\label{maggeo}
-4D_{\ga^{\sharp}}\ga^{\sharp} = \fd_{h}J(\ga^{\sharp}) = \fd_{h}(\star \ga)^{\sharp} = 2(d|\ga|_{h}^{2})^{\sharp}.
\end{align}
The tensor $A_{i}\,^{j}$ such that $A_{i}\,^{p}h_{pj} = -\tfrac{1}{2}F_{ij} = -\tfrac{1}{4}\fd_{h}\om_{ij}$ is $-\tfrac{1}{4}\fd_{h}J_{i}\,^{j}$, so it follows from the first equality of \eqref{maggeo} that the integral curves of $\ga^{\sharp}$ are magnetic geodesics. By the remarks preceeding the statement of the lemma, the geodesic curvature $\ka_{h}(\si)$ of an integral curve $\dot{\si}(t) = \ga^{\sharp}_{\si(t)}$ is $-\tfrac{1}{4}\fd_{h}|\dot{\si}|_{h}^{-1}$, which is constant along $\si$ by Theorem \ref{classtheorem}. Because $(\en, [h])$ is not exact, $|\ga|^{2}_{h}$ assumes a positive maximum at some point of $M$. Since, by \eqref{maggeo}, $2d|\ga|_{h}^{2} = \fd_{h}\star \ga$, at such a point it must be that $\fd_{h} = 0$. Thus $\fd_{h}$ has a zero which is not a zero of $\ga$, and the integral curve of $\ga^{\sharp}$ passing through this zero of $\fd_{h}$ has geodesic curvature $0$, so is a $D$-geodesic. Because $\ga^{\sharp}$ is conformal Killing there holds $\lie_{\ga^{\sharp}}J = 0$, from which follows  $[\ga^{\sharp}, J\ga^{\sharp}] = 0$. Using this and \eqref{maggeo} yields $D_{J\ga^{\sharp}}J\ga^{\sharp} = J(D_{J\ga^{\sharp}}\ga^{\sharp}) = J(D_{\ga^{\sharp}}J\ga^{\sharp}) = -D_{\ga^{\sharp}}\ga^{\sharp} = \tfrac{1}{4}\fd_{h}J\ga^{\sharp}$, which shows that the non-trivial integral curves of $J\ga^{\sharp}$ are projective geodesics of $D$. 
\end{proof}

This suggests viewing these particular exact magnetic flows as \textit{Einstein} and raises the question of whether there is a good notion of Einstein magnetic flows in higher dimensions. 

Theorem \ref{magnetictheorem} was suggested by the explicit models of Einstein Weyl structures described in section \ref{spheretorussection}, from which more precise information can be extracted. From the discussion following equation \eqref{spheresrfh} at the end of section \ref{spheremetricoccurssection} it follows that (with the setting and notations as in the statement of Theorem \ref{magnetictheorem}) in the case of the sphere the vector field $\ga^{\sharp}$ has two zeroes, and its integral curves, which by Theorem \ref{magnetictheorem} are magnetic geodesics, are simple closed circles separating these zeroes. Among these simple closed magnetic geodesics, there is a unique one of maximal energy on which $|\ga|^{2}_{h}$ attains its maximum value and $\fd_{h}$ vanishes, and which is moreover a $D$-geodesic, while for each positive energy below this maximum there are precisely two simple closed magnetic geodesics occurring as integral curves of $\ga^{\sharp}$.

\subsection{}
If a Riemannian signature AH structure $(\en, [h])$ on a compact, oriented surface $M$ is exact and Weyl, then the aligned representative $\nabla$ is the Levi-Civita connection of any Gauduchon metric $h \in [h]$, and \eqref{confscal} shows that $\uR_{h} = \sR_{h}$, so that $(\en, [h])$ is moreover Einstein just when $h$ has constant scalar curvature, that is $h$ is a constant curvature metric. Thus in this case $(\en,[h])$ is naturally identified with a positive homothety class of constant curvature metrics; since $[h]$ contains a unique such class, $(\en, [h])$ can be identified with $[h]$, and such an $(\en, [h])$ will be said simply to be \textit{a conformal structure}. In this case the weighted scalar curvature is parallel, and positive, zero, or negative according to whether the genus $g$ is $0$, $1$, or at least $2$. Alternatively, such an AH structure will be said to be \textit{generated} by its representative constant curvature metric.

\subsection{} 
Since, for a $\binom{p}{q+1}$ tensor $Q$, rescaling $h$ homothetically by $r \in \reap$ causes $|Q|_{h}^{2}|\det h|^{1/2}$ to rescale by $r^{p-q}$, on a compact surface the $L^{2}$-norm $||Q||_{h}^{2} = \int_{M}|Q|_{h}^{2}\,d\vol_{h}$ of a $\binom{p}{p+1}$ tensor $Q$ depends only on the conformal class of $h$, and not on the choice of $h$ itself, so in this case it makes sense to write $||Q||^{2} = ||Q||_{h}^{2}$, dropping the subscript indicating dependence on $h$. Let $(\en, [h])$ be a Riemannian Einstein AH structure $(\en, [h])$ on a compact, oriented surface $M$. Let $h \in [h]$ be a Gauduchon metric with associated Faraday primitive $\ga_{i}$. Then $||\bt||_{h}^{2} = \int_{M} \nbt$ and $||\ga||_{h}^{2}$ are unchanged if $h$ is homothetically rescaled. Although the $L^{2}$ norm of a vector field is changed by rescaling $h$, the $L^{2}$ norm $||\ga^{\sharp}||_{h}^{2} = ||\ga||_{h}^{2}$ is not because $\ga^{\sharp}$ also rescales when $h$ is rescaled, in a way that compensates for the change in norm. Similarly for the cubic form $B_{ijk} \defeq \bt_{ij}\,^{p}h_{pk}$ there holds $||B||_{h}^{2} = ||\bt||^{2}$. Whether there will be written $||\ga||^{2}$, $||\ga^{\sharp}||_{h}^{2}$, $||B||_{h}^{2}$, or $||\bt||^{2}$ will depend on the emphasis desired. From \eqref{confein2} it follows that there is a constant $\ka$ such that $\uR_{h} = 4|\ga|_{h}^{2} + \ka$. This $\ka$ does depend on $h$ in that rescaling $h$ by $r \in \reap$ rescales $\ka$ by $r^{-1}$. Integrating \eqref{confscal} against $d\vol_{h}$ using the Gau\ss-Bonnet theorem yields
\begin{align}\label{vortex}
\begin{split}
8\pi(1-g) & = \int_{M}\sR_{h}\,d\vol_{h} = \tfrac{1}{4}||\bt||^{2} + \int_{M}\uR_{h}\,d\vol_{h} = \tfrac{1}{4}||\bt||^{2} + 4||\ga||^{2} + \ka \vol_{h}(M).
\end{split}
\end{align}
in which $g$ is the genus of $M$. The number $\nv \defeq \ka \vol_{h}(M)$ does not depend on the choice of Gauduchon metric. For reasons explained in section \ref{vortexsection}, it is called the \textbf{vortex parameter} of the Riemannian Einstein AH structure on the compact, oriented surface $M$. Lemma \ref{vortexlemma} follows from \eqref{vortex}.

\begin{lemma}\label{vortexlemma}
For a Riemannian Einstein AH structure $(\en, [h])$ on a compact, orientable surface of genus $g$, the vortex parameter $\nv$ satisfies $\nv \leq 8\pi(1 - g)$, with equality if and only if $(\en, [h])$ is the AH structure generated by a constant curvature metric. 
\end{lemma}

If $(\en, [h])$ is exact Einstein, $\ka = \uR_{h}$ so $\nv = \uR_{h}\, d\vol_{h}(M) = \int_{M}\uR_{h} \,d\vol_{h} = \int_{M}R$ is the total weighted scalar curvature. In particular, this holds if the genus is greater than $1$, or the genus equals $1$ and $(\en, [h])$ is not Weyl. 

\subsection{}
Combining Theorem \ref{classtheorem}, Corollary \ref{holocorollary}, Lemma \ref{vortexlemma}, and making some arguments to relate topological conditions with the sign of the weighted scalar curvature yields Theorem \ref{scalarexacttheorem}, which is the principal structural result about Riemannian Einstein AH structures on compact surfaces.

Recall from the introduction the definition of a (strictly) convex flat real projective structure. 

\begin{theorem}\label{scalarexacttheorem}
For a Riemannian signature Einstein AH structure $(\en, [h])$ on a compact oriented surface $M$ of genus $g$ there holds one of the following mutually exclusive possibilities.
\begin{list}{(\arabic{enumi}).}{\usecounter{enumi}\leftmargin=.2cm}
\renewcommand{\theenumi}{Step (\arabic{enumi})}
\renewcommand{\theenumi}{\{\arabic{enumi}\}}
\renewcommand{\labelenumi}{\textbf{Level \theenumi}.-}
\item $(\en, [h])$ is exact and Weyl, so identified with the unique positive homothety class of constant curvature metrics contained in the underlying conformal structure $[h]$. There holds $\nv = 8\pi(1-g)$. 

\item\label{neg0} $R$ is negative and parallel, $(\en, [h])$ is exact and not Weyl, and $g \geq 2$. $(\en, [h])$ is projectively flat and conjugate projectively flat, and both $\en$ and $\ben$ are strictly convex. A distinguished metric $h \in [h]$ has scalar curvature of the form 
\begin{align}\label{negexactcurv}
\sR_{h} = \tfrac{1}{4}\left(|\bt|_{h}^{2} - \tfrac{||\bt||^{2}}{\vol_{h}(M)}\right) + \tfrac{8\pi(1-g)}{\vol_{h}(M)},
\end{align}
and $\sR_{h}$ is everywhere non-positive and somewhere negative. The cubic differential $B^{(3,0)}$, where $B_{ijk} \defeq \bt_{ij}\,^{p}h_{pk}$, is holomorphic. On the open submanifold  $M^{\ast}\defeq \{|B|^{2} \neq 0\}$, the metric $\sth_{ij} \defeq |B|_{h}^{2/3}h_{ij}$ is flat. There holds $\nv < 8\pi(1-g)$. 

\item\label{neg1}  $R$ is negative and parallel, $(\en, [h])$ is exact and not Weyl, and $M$ is a torus. A distinguished metric $h$ is flat and $B_{ijk} \defeq \bt_{ij}\,^{p}h_{pk}$ is parallel with respect to the Levi-Civita connection of $h$, so, in particular, has constant non-zero norm. $(\en, [h])$ is projectively flat and conjugate projectively flat, and both $\en$ and $\ben$ are convex but not strictly convex. There holds $\nv = -\tfrac{1}{4}||B||_{h}^{2} < 0$. 

\item\label{r0} $R = 0$, $(\en, [h])$ is Weyl and closed but not exact, and $M$ is a torus. A Gauduchon metric $h \in [h]$ is flat, the Faraday primitive $\ga_{i}$ of $h$ is parallel with respect to the Levi-Civita connection of $h$, and $\ga^{(1,0)}$ is holomorphic. The aligned representative of $(\en, [h])$ and is affinely flat, and its $(1,0)$ part is a holomorphic affine connection. There holds $\nv = -4||\ga||_{h}^{2} < 0$.

\item\label{rvarytorus} $R$ is somewhere positive and somewhere negative, $(\en, [h])$ is Weyl and not closed, and $M$ is a torus. For a Gauduchon metric $h \in [h]$ with Faraday primitive $\ga_{i}$ the scalar curvature is $\sR_{h} = 4(|\ga|^{2}_{h} - ||\ga||^{2})$, and there holds $\nv = -4||\ga||^{2} < 0$. 

\item\label{ps6} $R$ is somewhere positive, $(\en, [h])$ is Weyl and not closed, and $M$ is a sphere. A Gauduchon metric $h \in [h]$ with Faraday primitive $\ga_{i}$ has scalar curvature of the form 
\begin{align}\label{poscurv}
\sR_{h} = \uR_{h} = 4\left(|\ga|_{h}^{2} - \tfrac{||\ga||^{2}}{\vol_{h}(M)}\right)+ \tfrac{8\pi}{\vol_{h}(M)}.
\end{align}
On the open submanifold  $M^{\ast}\defeq \{|\ga|^{2} \neq 0\}$, the metric $\sth_{ij} \defeq |\ga|_{h}^{-2}h_{ij}$ is flat. There holds $\nv = 8\pi - 4||\ga||_{h}^{2} < 8\pi$. 
\end{list}
\end{theorem}

\begin{proof}
In this proof $h\in [h]$ is always a Gauduchon metric with Faraday primitive $\ga_{i}$, and $B_{ijk} = \bt_{ij}\,^{p}h_{kp}$.
\newcounter{lcount}
By \eqref{confein2}, $\ka \defeq \uR_{h} -4|\ga|_{h}^{2}$ is constant. The theorem follows by assembling the following claims.

\noindent
\begin{list}{[\alph{enumi}].}{\usecounter{enumi}\leftmargin=.2cm}
\renewcommand{\theenumi}{Step (\alph{enumi})}
\renewcommand{\theenumi}{[\alph{enumi}]}
\renewcommand{\labelenumi}{\textbf{Level \theenumi}.-}

\item If $(\en, [h])$ is not exact then $g < 2$, by Corollary \ref{holocorollary}, and $R$ is not everywhere negative, for, by \eqref{ne4}, at a point at which $|\ga|$ attains its maximum value there holds $\uR_{h} = -2\lap_{h}\log |\ga| \geq 0$.

\item  If $(\en, [h])$ is closed, then, by \eqref{conservationcondition}, $R$ is parallel, and it follows from Lemma \ref{2deinsteinlemma} that $(\en, [h])$ is projectively flat and conjugate projectively flat. If, moreover, $R = 0$, then by \eqref{rijkl} the aligned representatives of $(\en, [h])$ and its conjugate are flat.

\item \label{listb} If $(\en, [h])$ is not exact and $R$ is non-positive, then, by the maximum principle applied to \eqref{ne4b}, $\ga$ is closed, and so by \eqref{confein1}, $\ga$ is parallel. In particular $(\en, [h])$ is closed and $|\ga|^{2}_{h}$ is a non-zero constant. By \eqref{ne4} this forces $R \equiv 0$. Since $R$, $E_{ij}$ and $F_{ij}$ all vanish, the aligned representative $\nabla \in \en$ is affinely flat. By Corollary \ref{holocorollary}, since $(\en, [h])$ is not exact, it is Weyl, and so $\sR_{h} = \uR_{h}$ by \eqref{confscal}, which by Gau\ss-Bonnet forces $g = 1$. Hence $M$ is a torus and a Gauduchon metric is flat with parallel Faraday primitive. By \eqref{db2} there holds $-2\delbar \ga^{(1,0)} = F^{(1,1)} + \dad_{h}\ga\, \bar{h}^{(1,1)} = 0$, so the one-form $\ga^{(1,0)}$ is holomorphic, and hence the complex affine connection $\nabla^{1,0} = D^{1,0} - 2\ga^{(1,0)}$ is holomorphic.

\item \label{listd} If $R$ is non-positive and somewhere negative, then $(\en, [h])$ must be exact by \ref{listb}, and so $R$ is parallel by \eqref{conservationcondition}. Since $R$ is somewhere negative and parallel, it is strictly negative. Since $\uR_{h}$ is constant, $\sR_{h} = \uR_{h} + \tfrac{1}{4}|B|_{h}^{2}$ and $|B|^{2}_{h}$ assume their respective maximum values on $M$ at the same points. At a point at which $|B|^{2}$ takes its maximum value, \eqref{ne4c} yields $0 \geq \lap_{h}|B|^{2}_{h} = 3\sR_{h}|B|^{2}_{h}$, so at such a point $\sR_{h} \leq 0$. Since such a point is also a maximum of $\sR_{h}$, this shows $\sR_{h} \leq 0$. Solving \eqref{vortex} for $\ka$ shows that $\sR_{h}$ has the form \eqref{negexactcurv}. If $\sR_{h}$ is identically zero, integrating \eqref{negexactcurv} shows that $g = 1$; hence if $g > 1$ then $\sR_{h}$ must be somewhere negative. Conversely, if $g = 1$ then since $\sR_{h} \leq 0$, integrating \eqref{negexactcurv} shows that $\sR_{h}$ is identically zero, so that the distinguished metric $h$ is flat. In this case $0 = 4\sR_{h} = 4\uR_{h} + |B|_{h}^{2}$, so $|B|^{2}_{h}$ is constant, equal to $-4\uR_{h}$. 

\item \label{liste} If $g = 1$ and $(\en, [h])$ is exact then $R\leq 0$ by \eqref{vortex}. Since $R$ is parallel, either $R$ is identically $0$, in which case $(\en, [h])$ is Weyl, or $R$ is negative, in which case $\uR_{h}$ is a negative constant. Since $\uR_{h}$ is constant, $\sR_{h} = \uR_{h} + \tfrac{1}{4}|B|_{h}^{2}$ and $|B|_{h}^{2}$ assume their maximum values at the same points. By the same argument as in \ref{listd} there holds $\sR_{h} \leq 0$, and by Gau\ss-Bonnet this means $\sR_{h}$ is identically zero. Since $\uR_{h}$ is constant, this implies $|B|_{h}^{2}$ is constant. Since $B^{(3,0)}$ is holomorphic, it follows from \eqref{kato} of Lemma \ref{flatmetriclemma} that $0 = 2|d|B||^{2} = |DB|_{h}^{2}$, so that $B$ is parallel with respect to $h$. 

\item \label{listh} If $g \geq 2$, $(\en, [h])$ is exact by Corollary \ref{holocorollary} and hence $R$ is parallel. By \eqref{vortex}, $R$ is negative.

\item \label{listg} If $g = 0$, $(\en, [h])$ is Weyl by Corollary \ref{holocorollary}. Solving \eqref{vortex} for $\ka$ and substituting into $\sR_{h} = \uR_{h} = 4|\ga|^{2}_{h} + \ka$ yields \eqref{poscurv}. Since by \eqref{vortex}, $\int_{M}R = 8\pi$, $R$ must be positive somewhere. 

\item \label{listl} 
If $R$ is somewhere negative and somewhere positive then it is not parallel, and so by \eqref{conservationcondition} $(\en, [h])$ is not closed. Hence, by Corollary \ref{holocorollary}, $(\en, [h])$ is Weyl. In this case, since $\uR_{h}$ is somewhere negative, $\nv$ is negative.

\item \label{listf} If $g = 1$ and $(\en, [h])$ is not exact, then it is Weyl. By \eqref{vortex}, $0 = \int_{M}R$, so either $R$ is identically $0$, or $R$ is somewhere positive and somewhere negative. In the latter case, $(\en, [h])$ is not closed by \ref{listl}, and so from \eqref{vortex} it follows that $\nv = - 4||\ga||_{h}^{2} < 0$.

\item \label{listi} If $R$ is non-negative and somewhere positive then by \eqref{confscal} the same is true of $\sR_{h}$, and so $g = 0$ by Gau\ss-Bonnet and thus $(\en, [h])$ is Weyl by Corollary \ref{holocorollary}. If $(\en, [h])$ is moreover closed then the Einstein condition implies $R$ is parallel, and so everywhere positive, and hence $(\en, [[h])$ is exact by Lemma \ref{parallelexactlemma}. Thus in this case $(\en, [h])$ is exact Weyl and a distinguished metric has constant curvature.

\item \label{listj} That the flat projective structures in \eqref{neg0} and \eqref{neg1} are properly convex follows from Theorem \ref{convextheorem}, which is deferred to section \ref{convexsection} because its proof uses a point of view more conveniently introduced later. That the projective structures of \eqref{neg0} are strictly convex while those of \eqref{neg1} are not follows from Theorem $1.1$ of Y. Benoist (\cite{Benoist-convexesdivisiblesI}), that a discrete group which divides some properly convex domain in the projective sphere is Gromov hyperbolic if and only if the domain is strictly convex, coupled with the observation that the fundamental group of a surface of genus $g > 1$ is Gromov hyperbolic, while that of the torus is not. \qedhere
\end{list}
\end{proof}

In section \ref{constructionsection} it will be shown that all the possibilities identified in Theorem \ref{scalarexacttheorem} actually occur.

\subsection{}
Item \ref{listd} of the proof of Theorem \ref{scalarexacttheorem} shows that for an Einstein AH structure on a surface of genus $g > 1$ the norm squared of the cubic torsion with respect to a Gauduchon metric $h \in [h]$ satisfies the equivalent pointwise bounds
\begin{align}\label{calabibound}
&|B|_{h}^{2} = -4\uR_{h} + 4\sR_{h}  \leq -4\uR_{h}, & &|B|_{h}^{2} - ||B||_{h}^{2}\vol_{h}(M)^{-1} \leq 2\pi(g-1)\vol_{h}(M)^{-1}.
\end{align}
The pull back of such an Einstein AH structure to the universal cover of $M$ can be identified with that induced on an affine hypersphere asymptotic to the cone over the developing map image of the universal cover (see section \ref{convexsection}), and then the non-positivity of the scalar curvature of a Gauduchon metric is the conclusion of the main theorem (Theorem $5.1$) of Calabi's \cite{Calabi-completeaffine}. Here this non-positivity has been given a direct, autonomous proof. Based on these considerations it seems reasonable to expect the following.

\begin{theorem}
On an oriented surface $M$, let $(\en, [h])$ be an exact Einstein AH structure with negative weighted scalar curvature for which a distinguished metric is complete. Then the curvature of a distinguished metric is non-positive. 
\end{theorem}
\begin{proof}
Since the curvature of a distinguished metric $h \in [h]$ is $\uR_{h} + \tfrac{1}{4}|B|_{h}^{2}$ (the notations are as usual), it is equivalent to show that the function $u \defeq |B|_{h}^{2}$ is not greater than the constant $-4\uR_{h}$. There holds \eqref{ne4c}, and so $u$ satisfies a differential inequality of the form $\lap_{h}u \geq Bu^{2} - Au$ where $A = -3\uR_{h}$ and $B = 3/4$. A theorem of Cheng and Yau, from \cite{Cheng-Yau-maximalspacelike} and \cite{Cheng-Yau-affinehyperspheresI}, a statement and proof of which can be found as Theorem $6.6$ of \cite{Fox-ahs}, shows that for a complete metric $g$ with Ricci curvature bounded from below, a differential inequality of the form $\lap_{g}v \geq bv^{1+\si} - av$, with positive $b$ and holding where the non-positive smooth function $v$ is positive, implies an upper bound $v \leq |a/b|^{\si}$. Applying this with $\si = 1$, $v = u$, $a = A$ and $b = B$ yields the desired bound.
\end{proof}

There is no similar lower bound for $\sR_{h}$ in the case \ref{ps6} of Theorem \ref{scalarexacttheorem}. Suppose the Einstein AH structure $(\en, [h])$ is not exact, so is Weyl, and $M$ must be a sphere or a torus. If the constant $\ka = \uR_{h} - 4|\ga|_{h}^{2}$ is negative then the constant function $\log(4^{-1} |\ka|)$ solves
\begin{align}\label{wga}
\lap_{h}\phi + \ka + 4e^{\phi} = 0.
\end{align}
By \eqref{ne4} the function $\psi = \log |\ga|_{h}^{2}$ solves \eqref{wga} on the complement $M^{\ast}$ of the zero set of $\ga$, which is discrete because $h^{ip}\ga_{p}$ is the real part of a holomorphic vector field. Since $\psi$ tends to $-\infty$ at the zeroes of $\ga$ it is tempting to conclude that $\psi$ is bounded from above by $\log(4^{-1} |\ka|)$, but the standard comparison argument fails because in the operator $\lap_{h}\phi + \ka + 4e^{\phi}$ the zeroth order term is increasing in $\phi$, while the routine application of the maximum principle needs it to be non-increasing. In section \ref{spheremetricoccurssection} it is shown explicitly how to construct an Einstein Weyl structure on the two sphere $\sphere$ for which $\ka$ has an arbitrarily negative value.

\section{Relation with the Abelian vortex equations}\label{vortexsection}
In this section it is shown that an exact Einstein AH structure determines a solution of the Abelian vortex equations. On the other hand, an Einstein Weyl structure gives rise to a solution of equations superficially similar to the vortex equations, but differing from them by a change of sign in one term.

Let $M$ be a compact manifold and let $(g, J, \om)$ be a K\"ahler structure on $M$. Let $\lne$ be a smooth complex line bundle over $M$. The \textbf{Abelian vortex equations} with parameter $\tau$ are the following equations for a triple $(\nabla, k, s)$, in which $k$ is a Hermitian metric on $\lne$, $\nabla$ is a Hermitian connection on $(\lne, k)$, and $s$ is a smooth section of $\lne$:
\begin{align}\label{vortexeqns}
&\Om^{(0,2)} = 0,& &\delbar_{\nabla}s = 0,& &\j \dlef(\Om) + \tfrac{1}{2}|s|_{k}^{2} = \tfrac{1}{2}\tau.
\end{align}
Here $\Om_{ij}$ is the curvature of $\nabla$, viewed as a real-valued two-form on $M$; $\delbar_{\nabla}$ is the $(0,1)$ part of $\nabla$; and $\dlef$ is the dual Lefschetz operator given on $(1,1)$ forms by $\dlef(A) = -\om^{\al\bar{\be}}A_{\al\bar{\be}} = -\j A_{\si}\,^{\si}$. The first two equations say that $\delbar_{\nabla}$ is a holomorphic structure on $\lne$ with respect to which $s$ is a holomorphic section, while the third equation is something like an Einstein equation. A solution of \eqref{vortexeqns} is \textbf{non-trivial} if $s$ is not identically zero. The trivial solutions correspond to holomorphic structures on $\lne$; a precise statement is Theorem $4.7$ of \cite{Bradlow}.

The basic theorem about the Abelian vortex equations on a surface is the following.
\begin{theorem}[\cite{Noguchi}, \cite{Bradlow}, \cite{Garcia-Prada}]\label{vortextheorem}
Let $M$ be a compact surface equipped with a K\"ahler metric $(g, J)$. Let $\lne$ be a smooth complex line bundle with a fixed Hermitian metric $k$. Let $D$ be an effective divisor of degree equal to $\deg(\lne)$. There exists a non-trivial solution $(s, \nabla)$ of the vortex equations \eqref{vortexeqns}, unique up to gauge equivalence, if and only if $4\pi\deg(\lne) < \tau \vol_{g}(M)$. Moreover the holomorphic line bundle and section canonically associated to $D$ are $(\lne, \delbar_{\nabla})$ and $s$.
\end{theorem}

The space of effective divisors on $M$ of a given degree $r$ is the symmetric product $S^{r}(M)$ of $M$ and Theorem \ref{vortextheorem} shows that $S^{\deg(\lne)}(M)$ is in bijection with the moduli space of gauge equivalence classes of vortex solutions on $\lne$. It is shown in \cite{Garcia-Prada} by symplectic reduction, that this moduli space carries a K\"ahler structure.

Note that \textit{a priori} the K\"ahler metric $g_{ij}$ and the Hermitian metric $k$ on $\lne$ are unrelated. However, in what follows the interest will be on solutions for which $\lne$ \textit{qua} holomorphic bundle is identified with a power $\cano^{p}$ of the canonical bundle, and $k$ and $\nabla$ are the Hermitian metric and Hermitian connection induced by the underlying K\"ahler metric and connection on $M$. Equivalently, the corresponding effective divisor is a $p$-canonical divisor, that is it is in the $p$-fold product of the canonical divisor class. This motivates the following definition. A solution $(s, \nabla)$ of the Abelian vortex equations on $(M, g, J)$ is \textbf{$p$-canonical} if the divisor of $s$ is $p$-canonical; equivalently $s$ is a section of $\cano^{p}$ holomorphic with respect to the complex structure induced by $J$, and $\nabla$ is the Hermitian connection induced on $\cano^{p}$ by the Levi-Civita connection of $g$. 

Let $M$ be a compact orientable surface of genus $g$ and let $(\en, [h])$ be an exact Riemannian signature Einstein AH structure. Let $h \in [h]$ be a Gauduchon metric and $B_{ijk} = \bt_{ij}\,^{p}h_{pk}$. Let the K\"ahler structure on $M$ be that determined by $h$ and the given orientation. Let $k = h$ be the Hermitian metric induced on $\cano^{3}$ by $h$. The Levi-Civita connection $D$ of $h$ induces a Hermitian connection, also denoted by $D$, on $\cano^{3}$, for which the induced holomorphic structure is the canonical one. It is claimed that for $\tau = -3\uR_{h}$, the section $s = (3/2)^{1/2}B^{(3,0)}$ of $\cano^{3}$ solves the vortex equations \eqref{vortexeqns}. Note that $|s|_{h}^{2} = (3/2)|B^{(3,0)}|^{2} = (3/4)|\bt|^{2}_{h}$. The second equation of \eqref{vortexeqns} is valid by construction. The curvature of $D$ on $\cano^{3}$ is $\Om_{ij} = (3\j \sR_{h}/2)\om_{ij}$, so that $\Om^{(2,0)} = 0$ and $\dlef(\Om) = (3\j/2)\sR_{h}$. Because $(\en, [h])$ is exact, $\uR_{h}$ is constant, and it follows from \eqref{confscal} that $3\sR_{h} = 3(\uR_{h} + \tfrac{1}{4}|\bt|^{2}_{h}) = -\tau + |s|^{2}_{h}$. Hence
\begin{align*}
\j \dlef(\Om) + (1/2)|s|_{k}^{2} - \tfrac{1}{2}\tau = -(3/2)\left(\sR_{h}- \tfrac{1}{4}|B|_{h}^{2} - \uR_{h}\right) = 0,
\end{align*}
which shows the claim. Note that the resulting solution to the vortex equations is trivial exactly when $(\en, [h])$ is simply a conformal structure. In this case the vortex parameter $\nv$ is $-\tau/3$, which explains the terminology for $\nv$.

By Corollary \ref{holocorollary} if a Riemannian signature Einstein AH structure on a compact surface is not exact then it is Weyl. In this case, let $h \in [h]$ be a Gauduchon metric and $X^{i} = h^{ip}\ga_{p}$, where $\ga_{i}$ is the Faraday primitive of $h$. Let the K\"ahler structure on $M$ be that determined by $h$ and the given orientation. Let $k = h$ be the Hermitian metric induced on $\cano^{-1}$ by $h$. The Levi-Civita connection $D$ of $h$ induces a Hermitian connection, also denoted by $D$, on $\cano^{-1}$, for which the induced holomorphic structure is the canonical one. By Theorem \ref{classtheorem} there is a constant $\ka$ such that $\sR_{h} = \uR_{h} = \ka + 4|X|^{2}_{h}$. Let $\tau = \ka$ and consider the section $s = 2^{3/2}X^{(1,0)}$ of $\cano^{-1}$. The curvature of $D$ on $\cano^{-1}$ is $\Om_{ij} = -(\j \sR_{h}/2)\om_{ij}$, so that $\Om^{(2,0)} = 0$ and $\dlef(\Om) = -(\j/2)\sR_{h}$. Hence
\begin{align}\label{ewvortexlike}
\j \dlef(\Om) - (1/2)|s|_{k}^{2} - \tfrac{1}{2}\tau = (1/2)\left(\sR_{h}  - 4  |X|_{k}^{2} - \ka\right) = 0.
\end{align}
The equation \eqref{ewvortexlike} differs from \eqref{vortexeqns} by the change of sign on the $|s|_{k}^{2}$ term. This affects the interpretation of the equations. The vortex equations are the adaptation to compact surfaces of the Landau-Ginzburg equations modeling the macroscopic behavior of superconducting materials. In the context of superconductors the wrong sign on the $|s|_{k}^{2}$ term corresponds to a physically unreasonable negative sign on the quartic term in the free energy. While the equations \eqref{ewvortexlike} by themselves make sense, and, as follows from section \ref{spheretorussection}, have solutions, they are not the usual vortex equations, and a more satisfactory derivation of them is needed before any significance can be attached to their similarity with the vortex equations.

Nonetheless, the preceeding can be given the following uniform, albeit unmotivated description. Let the setup be as in the two preceeding paragraphs. Let $q$ be either $-1$ or $3$ and let $\si$ be $X$ or $B$. Let $\ka$ be the constant such that $\uR_{h} = \ka + 4|\ga|^{2}_{h}$. Let $\tau = -q\ka$ and $s = 2^{(2-q)/2}|q|^{1/2} \si^{(k,0)}$. Then $s$ is a section of $\cano^{q}$ solving the following modified vortex equations with respect to the holomorphic structure and Hermitian metric and connection on $\cano^{q}$ induced by $h$:
\begin{align}\label{modifiedvortexeqns}
&\Om^{(0,2)} = 0,& &\delbar_{\nabla}s = 0,& &2\j \dlef(\Om) + \sign(q)|s|_{k}^{2} = \tau.
\end{align}

Note that distinct exact Einstein AH structures need not determine gauge inequivalent solutions of \eqref{vortexeqns}. If $B^{(3, 0)}$ is replaced by $e^{\j\theta}B^{(3, 0)}$ for a constant $\theta$, the resulting vortex solutions are gauge equivalent. In section \ref{gaugesection} it is shown that the real part of $e^{\j\theta}B^{(3, 0)}$ is the cubic torsion of an Einstein AH structure with the same underlying conformal structure and Gauduchon metric. The solutions to the Abelian vortex equations with $\deg(D) = 3\chi(M)$ which arise in this way from exact Einstein AH structures are exactly the $3$-canonical solutions. This essentially means that on a compact, orientable surface of genus $g > 1$ the quotient of the moduli space of Einstein AH structures by the action of $\comt$ is the space of gauge equivalence classes of $3$-canonical Abelian vortices; see section \ref{constructionsection} for related discussion.

The existence of $p$-canonical Abelian vortices is not immediate from Theorem \ref{vortextheorem}; it is demonstrated in Corollary \ref{vortexcorollary}. In particular, the existence of Einstein AH structures as in \ref{neg0} of Theorem \ref{scalarexacttheorem} is not immediate from the existence theorem for Abelian vortices; the complication is the additional requirement of compatibility of the underlying K\"ahler structure (the putative Gauduchon metric) and the holomorphic and metric structure on the line bundle.

However, note that the constraint on the vortex parameter given by Lemma \ref{vortexlemma} also follows from Theorem \ref{vortextheorem} and the reinterpretation of the Einstein AH equations in terms of the vortex equations.

\section{Einstein AH structures on compact orientable surfaces of genus at least two}\label{constructionsection}
In this section and the next the classification of Riemannian Einstein AH structures on compact, orientable surfaces is completed by showing that all the possibilities identified in Theorem \ref{scalarexacttheorem} are realized. Throughout $M$ is  a compact, orientable surface of genus $g$. 

Theorem \ref{scalarexacttheorem} shows that a Riemannian Einstein AH structure on $M$ must determine data of one of the following forms
\begin{enumerate}
\item A constant curvature metric. (Any $g$). 
\item A conformal structure and a non-trivial holomorphic vector field. ($g \in \{0, 1\}$).
\item A conformal structure and a non-trivial holomorphic cubic differential. ($g \geq 1$).
\end{enumerate} 
In order to complete the classification it is necessary to show how to construct from data of type $(2)$ or $(3)$ a Riemannian Einstein AH structure, and to analyze when two Riemannian Einstein AH structures are equivalent modulo $\diff^{+}(M)$. The second step is basically straightforward because of the uniformization theorem, and so the content of this section and the next is the analysis of the first step. For whatever $g$ the given data of a holomorphic section and a conformal structure determine an elliptic PDE for the conformal factor expressing a putative Gauduchon metric in terms of a constant scalar curvature background metric. 

For $g \geq 2$ the existence of a unique solution for the resulting PDE is a straightforward application of standard elliptic PDE techniques. This is explained in the present section. For $g \in \{0,1\}$ the uniqueness fails because the conformal automorphism group is large, but the holomorphic vector field induces an $S^{1}$ symmetry which can be used to reduce the PDE to an ODE which is easily solved. This is described in section \ref{spheretorussection}.

\subsection{}
On a smooth surface $M$, associate to a triple $(h, F, B)$ comprising a Riemannian metric $h$, a smooth function $F \in \cinf(M)$, and the real part $B$ of a smooth $k$-differential, the differential operator defined for $\phi \in \cinf(M)$ by
\begin{align}\label{wdefined}
\ahop(h, F, B)(\phi) \defeq  \lap_{h}\phi - \sR_{h} + Fe^{\phi} + 2^{1-k}e^{(1-k)\phi}|B|_{h}^{2}.
\end{align}
It is convenient to include the constant factor $2^{1-k}$. More generally, there can be several $B$'s, for different $k$'s, e.g. for a holomorphic one-form $X^{(1,0)}$, and a holomorphic cubic differential $B^{(3,0)}$, it is convenient to write $\ahop(h, F, X, B)(\phi) = \lap_{h}\phi - \sR_{h} + Fe^{\phi} + 4e^{2\phi}|X|_{h}^{2} + \tfrac{1}{4}e^{-2\phi}|B|_{h}^{2}$. The metric $h$ will be called the \textbf{background} metric of the equation. Note that in the analysis of $\ahop$ the holomorphicity or not of $B^{(k,0)}$ plays \textit{no} role; it is important only for the properties of the objects constructed from the solutions. For $\mu, \la \in \cinf(M)$ there holds the scaling rule 
\begin{align}\label{wrescale}
\begin{split}
e^{\mu} \ahop(e^{\mu} h, e^{\la} F, e^{(1-k)\la/2}B)(\phi - \mu - \la) = \ahop(h, F, X, B)(\phi) - \lap_{h}\la.
\end{split}
\end{align}

\subsection{}
Let $(\en, [h])$ be an AH structure on a compact orientable surface $M$ of genus $g$. Let $\nabla \in \en$ be its aligned representative and let $\tilde{h}_{ij} = e^{-\phi}h_{ij} \in [h]$ and $h \in [h]$ have Levi-Civita connections related by $\tD = D + 2\si_{(i}\delta_{j)}\,^{k} - h_{ij}\si^{k}$, in which $2\si_{i} \defeq -d\phi_{i}$ and $\si^{i} \defeq h^{ip}\si_{p}$. Recall the notational conventions established in section \ref{riemannsurfacehodgestarsection}. Let $\ga_{i}$ be the Faraday primitive of $h$ and let $B_{ijk} \defeq \bt_{ij}\,^{p}h_{pk}$. Using \eqref{conformalscalardiff}, \eqref{confscal}, $\dad_{h}\ga = e^{-\phi}\dad_{\tilde{h}}\ga$, and $\lap_{h}\phi = e^{-\phi}\lap_{\tilde{h}}\phi$ there results
\begin{align}\label{weqn}
\begin{split}
\ahop(\tilde{h}, \uR_{h}, B)(\phi) + \dad_{\tilde{h}}\ga& =  \lap_{\tilde{h}}\phi - \sR_{\tilde{h}} + e^{\phi}\uR_{h} + \tfrac{1}{4}e^{-2\phi}|B|_{h}^{2} + \dad_{\tilde{h}}\ga\\
& =  e^{\phi}\left(\lap_{h}\phi + \uR_{h} + \tfrac{1}{4}|\bt|^{2}_{h} + \dad_{h}\ga \right) - \sR_{\tilde{h}}= 0.
\end{split}
\end{align}
Equation \eqref{weqn} can be used to construct AH structures. The metric $\tilde{h}$ will be treated as the background metric and $h$ (equivalently $\phi$) as the unknown. If $\tilde{h}_{ij}$, $\ga_{i}$, $\uR_{h}$, and $B_{ijk}$ are given, then solving \eqref{weqn} for $\phi$ yields an AH structure with cubic torsion $\bt_{ij}\,^{k} = h^{kp}B_{ijp}$ and scalar curvature $R = |\det h|^{1/2}\uR_{h}$, for which $h = e^{\phi}\tilde{h}$ is a representative metric with Faraday primitive $\ga_{i}$. It is convenient to seek Gauduchon $h$, in which case it must be that $\dad_{\tilde{h}}\ga = e^{\phi}\dad_{h}\ga = 0$, and the equation to be solved reduces to $\ahop(\tilde{h}, \uR_{h}, B)(\phi) = 0$. If the AH structure obtained by solving \eqref{weqn} is to be Einstein then it must be that $B^{(3,0)}$ is holomorphic, and there must be a constant $\ka$ and a holomorphic vector field $X^{(1,0)}$ such that $X^{p}h_{ip} = \ga_{i}$ is the Faraday primitive of $h$ and $\uR_{h} = 4|X|_{h}^{2} + \ka = 4e^{\phi}|X|^{2}_{\tilde{h}} + \ka$. In \eqref{weqn} this yields
\begin{align}\label{geneq}
\ahop(\tilde{h}, \ka, X, B )(\phi) = \lap_{\tilde{h}}\phi - \sR_{\tilde{h}} +\ka e^{\phi} + 4e^{2\phi}|X|_{\tilde{h}}^{2} + \tfrac {1}{4}e^{-2\phi}|B|^{2}_{\tilde{h}} = 0,
\end{align}
as the equation to be solved. Of course, by Corollary \ref{holocorollary}, if there is to be a solution only one of $X$ and $B$ can be non-zero. 

It follows from \eqref{wrescale} that $\phi$ solves \eqref{geneq} if and only if for any $\mu \in \cinf(M)$ and $r \in \reap$ the function $\psi = \phi - \log r - \mu$ solves $\ahop(e^{\mu}h, r\ka, rX, r^{-1}B)(\psi) = 0$. The resulting metrics $\bar{h}_{ij} = e^{\psi}e^{\mu}\tilde{h}_{ij} = r^{-1}e^{\phi}\tilde{h}_{ij}$ and $h_{ij} = e^{\phi}\tilde{h}_{ij}$ are positively homothetic, while the resulting tensors $\ga_{i} = X^{p}h_{pi} = rX^{p}\tilde{h}_{pi}$ and $\bt_{ij}\,^{k} = h^{kp}B_{ijp} = \bar{h}^{kp}r^{-1}B_{ijp}$ are the same, so that the AH structures resulting from these solutions are the same. Thus such rescaling is trivial from the point of view of constructing AH structures, and it is natural to restrict the allowable $\phi$ by imposing some normalization which eliminates this freedom. 

The scaling in $\tilde{h}$ given by $\mu$ is eliminated by fixing a convenient background metric $\tilde{h}$. The freedom in $r$ is most naturally eliminated by imposing some condition on the resulting metric $h = e^{\phi}\tilde{h}$, e.g. fixing the minimum of $\sR_{h}$. There are several possibilities. An obvious one is to require that $\uR_{h}$ take a specific value; this amounts to fixing $\ka$. Another is to demand that $\vol_{h}(M) = \int_{M}e^{\phi}\,d\vol_{\tilde{h}}$ have some prespecified value; as will be made precise below, this puts some conditions on $\ka$ necessary for the existence of solutions. The curvature normalization is probably more natural from the geometric point of view, while the volume normalization is probably more natural from the point of view of partial differential equations. A function $\phi \in \cinf(M)$ is \textbf{volume normalized} with respect to $\tilde{h}$ if $\int_{M}e^{\phi}d\vol_{\tilde{h}} = \vol_{\tilde{h}}(M)$, in which case $\vol_{h}(M) = \vol_{\tilde{h}}(M)$. From Lemma \ref{vortexlemma} it follows that for there to exist a volume normalized $\phi$ solving \eqref{geneq} it is necessary that $\ka \leq 4\pi\chi(M)/\vol_{\tilde{h}}(M)$. Moreover, for equality to hold it is necessary that both $X$ and $B$ be identically zero, in which case \eqref{geneq} becomes simply $\sR_{e^{\phi}\tilde{h}} = e^{-\phi}\left(\sR_{\tilde{h}} - \lap_{\tilde{h}}\phi\right) = 4\pi\chi(M)/\vol_{\tilde{h}}(M)$. This is simply the equation that $e^{\phi}\tilde{h}$ have constant curvature; in this case there is by the uniformization theorem a unique normalized $\phi$ solving the equation. Such a solution to \eqref{geneq} will be referred to as \textbf{uniformizing}. 

In summary, if $\phi$ is a $\cinf$ solution of \eqref{geneq} then $h_{ij} \defeq e^{\phi}\tilde{h}_{ij}$ is the Gauduchon metric with associated Faraday primitive $\ga_{i} = X^{p}h_{ip}$ of an Einstein AH structure having cubic torsion $\bt_{ij}\,^{k} = h^{kp}B_{ijp}$, scalar curvature $\uR_{h} = 4|X|_{h}^{2} + \ka$, and vortex parameter $\nv = \ka \vol_{h}(M)$. The natural choices for the background metric are the metric $\tilde{h} \in [h]$ having constant scalar curvature in $\{0, \pm 2\}$, and the metrics $\sth = |X|_{\tilde{h}}^{-2}\tilde{h}$ and $\sth = |B|_{\tilde{h}}^{2/3}\tilde{h}$ (which do not depend on the choice of $\tilde{h} \in [h]$), which are flat by Lemma \ref{flatmetriclemma}. The metrics $\sth$ are defined only on an open subset of $M$, so if they are to used there have to be imposed boundary conditions on $\phi$. 

\subsection{} 
The existence of a unique solution to  $\ahop(h, k, B)(\phi) = 0$ when $\chi(M) \leq 0$, $B$ is somewhere non-zero, and $k \leq 4\pi\chi(M)/\vol_{h}(M)$ follows from Lemmas \ref{wuniqlemma} and Lemma \ref{wexistencelemma}. These lemmas are, however, stated more generally, so as to yield also Corollary \ref{vortexcorollary}, which shows how to associate to a Riemann surface equipped with a non-trivial holomorphic $k$-differential a canonical smooth metric in the given conformal structure, by solving the Abelian vortex equations. For a cubic holomorphic differential, the existence and uniqueness statements of Lemmas \ref{wuniqlemma} and \ref{wexistencelemma} for a metric of constant scalar curvature $-2$ on a compact oriented surface $M$ of genus greater than $1$ and $k = 2$, are Theorem $4.0.2$ of Loftin's \cite{Loftin-affinespheres}, and the proofs given here are very similar to Loftin's, which are themselves based on standard arguments as in \cite{Kazdan-Warner} or, particularly, \cite{Wolf-teichmuller}. The generality of allowing nonconstant $k$ is convenient for the example of section \ref{naiveexample}.

\begin{lemma}\label{wuniqlemma}
If $M$ is a compact, orientable surface with $\chi(M) \leq 0$, $h$ is a Riemannian metric, and, for an integer $p \geq 1$, $B$ is the real part of a non-trivial smooth $p$-differential, then the equation $\ahop(h, k, B)(\phi) = 0$ has at most one solution in $\cinf(M)$ for each $k \in \cinf(M)$ satisfying $k \leq 0$. 
\end{lemma}
\begin{proof}
Suppose $\phi, \psi \in \cinf(M)$ solve $\ahop(h, k, B)(\phi) = 0$ for some $k \in \cinf(M)$. Then
\begin{align}\label{lapdiff}
\lap_{h}(\phi - \psi) = - k (e^{\phi} - e^{\psi}) - 2^{(1-p)/2}(e^{(1-p)\phi} - e^{(1-p)\psi})|B|^{2}_{h}.
\end{align}
If $k$ is strictly negative the claim follows from the maximum principle applied to \eqref{lapdiff}. If $k \leq 0$, it follows from \eqref{lapdiff} that
\begin{align}\label{uniqineq}
\begin{split}
\lap_{h}(\phi& - \psi)^{2}  = 2(\phi - \psi)\lap_{h}(\phi - \psi) + 2|d(\phi - \psi)|^{2}_{h} \\
& = -2k(\phi - \psi)(e^{\phi} - e^{\psi}) - 2^{2-p}|B|^{2}_{h}(\phi - \psi)(e^{(1-p)\phi} - e^{(1-p)\psi})\geq 2|d(\phi - \psi)|^{2}_{h} \geq 0,
\end{split}
\end{align}
so that $(\phi - \psi)^{2}$ is subharmonic and therefore constant, since $M$ is compact. Write $\psi = \phi + c$ for a constant $c$, and substitute this into \eqref{lapdiff} to obtain $e^{-p\phi}|B|_{h}^{2}(1-e^{(1-p)c}) = 2^{p-1}k(e^{c} - 1)$. If not both $p = 1$ and $k \equiv 0$, then, since the signs of $e^{c} - 1$ and $1- e^{-2c}$ are the same if $c \neq 0$, and, by hypothesis, $k \leq 0$ and $B$ is not identically zero, this can be only if $c = 0$. If both $p =1$ and $k \equiv 0$, there is no solution, for if $\ahop(h, k, B)(\phi) = 0$, then integration yields $0 \geq 4\pi\chi(M) = \int_{M}\sR_{h}\,d\vol_{h} = ||B||^{2}_{h}$, contradicting the assumption that $B$ is not identically zero.
\end{proof}

\begin{lemma}\label{genuszerolemma}
Let $M$ be a smooth torus. Let $h$ be a flat Riemannian metric, and let $B$ be the real part of a non-trivial holomorphic $p$-differential for $p \geq 1$. Then $B$ is parallel and $|B|_{h}^{2}$ is constant. For each constant $\ka < 0$ the unique solution $\phi$ to the equation $\ahop(h, \ka, B)(\phi) = 0$ is the constant function $\tfrac{1}{p}\log\left(2^{1-p}|\ka|^{-1}|B|_{h}^{2}\right)$.
\end{lemma}

\begin{proof}
  By Corollary \ref{rrcorollary}, $B$ is $h$-parallel, and hence $|B|_{h}^{2}$ is a non-zero constant, so the given $\phi$ solves $\ahop(h, \ka, B)(\phi) = 0$. This is the unique solution by Lemma \ref{wuniqlemma}. 
\end{proof}

Lemma \ref{rootestimatelemma} is needed in the proof of Lemma \ref{wexistencelemma}.
\begin{lemma}\label{rootestimatelemma}
Let $0 < p \in \integer$ and $0 < a, b \in \rea$. The unique positive root $r_{1}$ of $f(r) = r^{p} - ar^{p-1} - b$ satisfies $a \leq r_{1} \leq a + b^{1/p}$.
\end{lemma}
\begin{proof}
For $p = 1$ the positive root of $g$ is $a + b$, while for $p = 2$ it is $(a + \sqrt{a^{2} + 4b})/2 \leq a + b^{1/2}$, so suppose $p > 2$. Since $f$ is negative at $ r= 0$, has a negative minimum at $a(1 - 1/p)$, is monotone decreasing on $(0, a(1 - 1/p))$ and monotone increasing on $(a(1-1/p), \infty)$, and satisfies $\lim_{r \to \infty}f(r) = \infty$, it has a unique positive real root $r_{1}$ which is greater than $a(1 - 1/p)$, and, since $a + b^{1/p} > a > a(1 - 1/p)$, it suffices to observe $f(a) = -b < 0$ and
\begin{align*}
\begin{split}
f(a + b^{1/p}) & = (a + b^{1/p})^{p-1}b^{1/p} - b = \sum_{s = 0}^{p-2}\binom{p-1}{s}b^{(s+1)/p}a^{p-1-s} > 0. 
\end{split}\qedhere
\end{align*}
\end{proof}

\begin{lemma}\label{wexistencelemma}
Let $M$ be a compact, orientable surface with $\chi(M) < 0$. Let $h$ be a Riemannian metric, let $B$ be the real part of a smooth $p$-differential not everywhere zero, and let $k \in \cinf(M)$ be everywhere negative. Let $q = \min_{M}\sR_{h}$, $Q= \max_{M}\sR_{h}$, $P = \max_{M}|B|^{2}_{h}$, $\ka = \min_{M}k$ and $K = \max_{M}k$. If the curvature of $h$ is negative, then there is a unique solution $\phi$ to the equation $\ahop(h, \ka, B)(\phi) = 0$, and it satisfies
\begin{align}\label{bbound}
\ka^{-1}(\max_{M}\sR_{h}) \leq e^{\phi} 
\leq K^{-1}(\min_{M}\sR_{h}) + 2^{(1-p)/p}|K|^{-1/p}(\max_{M}|B|^{2/p}_{h}).
\end{align}
\end{lemma}

\begin{proof}
The uniqueness follows from Lemma \ref{wuniqlemma}. Following the proof of Theorem $4.0.2$ in \cite{Loftin-affinespheres} the existence is demonstrated by applying Theorem V.$1.1$ of \cite{Schoen-Yau}, which shows that if a semi-linear elliptic equation $\lap_{h}u + F(x, u)$ with $F \in \cinf(M\times \rea)$ on a compact manifold $M$ admits a $C^{2}$ supersolution $u^{+}$ and a $C^{2}$ subsolution $u^{-}$ such that $u^{-} \leq u^{+}$ then it admits a $\cinf$ solution $u$ such that $u^{-} \leq u \leq u^{+}$ (this is proved by a modification of a standard iteration argument closely related versions of which can be found in section $2$ of the appendix to the fourth chapter of \cite{Courant-Hilbert-II} and section $9$ of \cite{Kazdan-Warner}). Suppose $Q < 0$. Since both $Q$ and $\ka$ are negative, $Q/\ka$ is a positive constant, and $\ahop(h, k, B)(\log (Q/\ka)) \geq -Q + kQ/\ka \geq  0$, so $\log(Q/\ka)$ is a subsolution. On the other hand, the polynomial $f(r) = K r^{3} - qr^{2} + \tfrac{1}{4}P$ has the unique real root $r_{1}$, which is easily seen to be positive and no smaller than $q/K$, which in turn is no smaller than $Q/\ka$. As $\ahop(h, \ka, B)(\log r_{1}) \leq r_{1}^{-1}f(r_{1}) \leq 0$, $\log r_{1}$ is a supersolution. There follows $\log (Q/\ka) \leq \phi \leq \log r_{1}$. The bound \eqref{bbound} follows from Lemma \ref{rootestimatelemma} with $a = q/K$ and $b = 2^{1-p}|K|^{-1}P$.
\end{proof}

There follows the existence of $p$-canonical Abelian vortices, as claimed in section \ref{vortexsection}.

\begin{corollary}\label{vortexcorollary}
Let $(M, [h], J)$ be a compact Riemann surface of genus $g > 1$, let $0 < p \in \integer$, let $B^{(p, 0)}$ be a non-trivial holomorphic $p$-differential, and let $\ka$ be a negative constant. There is a unique representative $h \in [h]$ such that $s = 2^{1 -p/2}p^{1/2}B^{(p,0)}$ solves the Abelian vortex equations \eqref{vortexeqns} for $\tau = -p\ka$ and the Hermitian structure induced on $\cano^{p}$ by $h$.
\end{corollary}

\begin{proof}
Let $\tilde{h} \in [h]$ and write $h = e^{\phi}\tilde{h}$. Arguing as in section \ref{vortexsection} the putative $h$ must satisfy
\begin{align*}
\begin{split}
0 & = \tfrac{2}{p}\left(\j \dlef(\Omega) + \tfrac{1}{2}|s|^{2}_{h} - \tfrac{1}{2}\tau\right) = -\sR_{h} + 2^{2-p}|B^{(p,0)}|^{2}_{h} + \ka\\
& = e^{-\phi}\left(\lap_{\tilde{h}}\phi - \sR_{\tilde{h}} + \ka e^{\phi} + 2^{1-p}e^{(1-p)\phi}|B|^{2}_{\tilde{h}}\right) = e^{-\phi}\ahop(\tilde{h}, \ka, B)(\phi).
\end{split}
\end{align*}
For $\tilde{h}$ such that $\sR_{\tilde{h}} = -2$, the existence of a unique solution follows from Lemmas \ref{wuniqlemma} and \ref{wexistencelemma}.
\end{proof}

For a compact orientable surface, Corollary \ref{vortexcorollary} associates to each $(J, B) \in \prekhodge(M)$ a distinguished metric representing the conformal structure. If $(J, B)$ is constructed from a flat metric $\sth$ with conical singularities as in section \ref{singularmetricsection}, then this associates to such a metric a smooth conformal metric in a canonical manner. This observation might be of use in the study of such metrics.

\begin{lemma}\label{factorboundlemma}
Let $1 \leq p \in \integer$. Let $M$ be a compact, orientable surface with $\chi(M) < 0$, let $\tilde{h}$ be a Riemannian metric with negative scalar curvature satisfying $\min_{M}\sR_{\tilde{h}} = -2$, let $B^{(p,0)}$ be holomorphic, and let $\ka$ be a negative constant. If $\phi$ solves $\ahop(\tilde{h}, \ka, B)(\phi) = 0$, and $h_{ij} = e^{\phi}\tilde{h}_{ij}$ then 
\begin{align}\label{calabibound2}
& -2\ka^{-1} + 2^{(1-p)/p}|\ka|^{-1/p}(\max_{M}|B|^{2/p}_{\tilde{h}})  \geq  e^{\phi} \geq 2^{(1-p)/p}|\ka|^{-1/p}|B|_{\tilde{h}}^{2/p},\\
\label{volumebound}
\begin{split}
(-2\ka^{-1}& + 2^{(1-p)/p}|\ka|^{-1/p}(\max_{M}|B|^{2/p}_{\tilde{h}})  )\vol_{\tilde{h}}(M) \geq\\
&\vol_{h}(M) \geq 2^{(1-p)/p}|\ka|^{-1/p}\int_{M}|B|^{2/p}_{\tilde{h}}\,d\vol_{\tilde{h}} = 2^{(1-p)/p}|\ka|^{-1/p}\vol_{\sth}(M).
\end{split}
\end{align}
Here $\sth \in [h]$ is the singular flat metric $\sth = |B|^{2/p}_{\tilde{h}}\tilde{h}$.
\end{lemma}
\begin{proof}
In the $p = 3$ case, by \eqref{calabibound} applied to the resulting exact Einstein AH structure with Gauduchon metric $h = e^{\phi}\tilde{h}$, there holds $-4\ka = -4\uR_{h} \geq |B|_{h}^{2} = e^{-3\phi}|B|_{\tilde{h}}^{2}$, which is the second inequality of \eqref{calabibound2}. Alternatively, the following proof, applicable for all $p \geq 1$, can be viewed as giving another proof of \eqref{calabibound}. By \eqref{keromlap} of Lemma \ref{flatmetriclemma}, $\psi = \log \left( 2^{(1-p)/p}|\ka|^{-1/p} |B|_{\tilde{h}}^{2/p}\right)$ solves $\ahop(\tilde{h}, \ka, B)(\psi) = 0$ on the complement of the zero set of $B$. (Essentially this observation is key in both \cite{Noguchi} and the proof of Proposition $1$ of \cite{Loftin-cubic}). Since $\psi$ goes to $-\infty$ on the zero set of $B$, there is $\ep > 0$ so that on the boundary $\pr M^{\ep}$ of the complement $M^{\ep}$ of an $\ep$ neighborhood of the zero set of $B$ there holds $\psi \leq \min_{M}\phi$. Since the zeroth order part of $\ahop(\tilde{h}, \ka, B)(u)$ is non-increasing in $u$, it follows that 
\begin{align*}
\lap_{\tilde{h}}(\psi - \phi) = -\ka(e^{\psi} - e^{\phi}) - 2^{1-p}(e^{(1-p)\psi} - e^{(1-p)\phi})|B|_{\tilde{h}}^{2}
\end{align*}
 is non-negative on the domain $U = \{x \in M^{\ep}: \psi(x) > \phi(x)\}$. Since by the choice of $\ep$ the closure of $U$ is contained properly in $M^{\ep}$, the maximum principle implies $U$ is empty, showing that $\psi \leq \phi$ on $M^{\ep}$. Letting $\ep \to 0$ yields the second inequality of \eqref{calabibound2}. The first inequality of \eqref{calabibound2} follows from \eqref{bbound}. Integrating \eqref{calabibound2} yields the volume bound \eqref{volumebound}.
\end{proof}

\subsection{Example: naive Einstein AH structures which are not Einstein}\label{naiveexample}
On a compact surface $M$ of genus $g > 1$, let $\tilde{h}$ be a Riemannian metric with scalar curvature $-2$, and let $B^{(3,0)}$ be a cubic holomorphic differential. Let $k$ be a smooth function on $M$ which satisfies $k < 0$. By Lemmas \ref{wuniqlemma} and \ref{wexistencelemma} the equation $\ahop(\tilde{h}, k, B)(\phi) = 0$ admits a unique smooth solution $\phi$. Let $h_{ij} = e^{\phi}\tilde{h}_{ij}$, $\bt_{ij}\,^{k} = h^{kp}B_{ijp}$, and $\nabla = D  - \tfrac{1}{2}\bt_{ij}\,^{k}$, in which $D$ is the Levi-Civita connection of $h_{ij}$. Then $\nabla$ is the aligned representative of the AH structure which it generates with $[h]$, which is exact, and $h$ is a distinguished representative of $[h]$. There holds $D_{p}\bt_{ij}\,^{p} = h^{pq}D_{p}B_{ijq} = 0$ because $B$ is the real part of a holomorphic differential, and so by \eqref{ddivbt}, $(\en, [h])$ is naive Einstein. On the other hand, by \eqref{confscal}, there holds $\uR_{h} = \sR_{h} - \tfrac{1}{4}|\bt|_{h}^{2} = \sR_{h} - \tfrac{1}{4}|B|_{h}^{2}$, while by \eqref{conformalscalardiff} and the construction of $\phi$ there holds $\sR_{h} = e^{-\phi}(-2 - \lap_{\tilde{h}}\phi) = k + \tfrac{1}{4}|B|_{h}^{2}$, so that $\uR_{h} = k$. Thus if $k$ is not constant, then $\uR_{h}$ is not constant, so by \eqref{confein1}, $(\en, [h])$ is not Einstein. 

\subsection{}
Let $M$ be a compact, oriented surface. Recall the spaces defined in section \ref{khodgesection}. Let $\epremod(M)$ be the space of Riemannian signature Einstein AH structures on $M$. The underlying conformal structure of each element of $\epremod(M)$ determines a complex structure inducing the given orientation, and so there is an evidently sujective map $\epremod(M) \to \jpremod(M)$. The group $\diff^{+}(M)$ acts on $\epremod(M)$ by pullback. Define the \textbf{deformation space} $\tmod(M) = \epremod(M)/\diff_{0}(M)$ of Einstein AH structures on $M$ and the \textbf{moduli space} $\emod(M) = \epremod(M)/\diff^{+}(M)$ of Einstein AH structures on M. The oriented mapping class group $\map^{+}(M)$ acts on $\tmod(M)$ with quotient $\emod(M)$. The canonical map $\epremod(M) \to \jpremod(M)$ descends to a map $\tmod(M) \to \teich(M)$, and similarly at the level of moduli spaces. 

\subsection{}\label{modulisection}
Suppose the compact, oriented surface $M$ has genus $g > 1$. For each $a \in \reap$ define a map $\prebij^{a}:\epremod(M) \to \precubic(M)$ by $\prebij^{a}(\en, [h]) = (J, B)$ in which $J$ is the complex structure determined by $[h]$ and the given orientation, and $B$ is defined by $B_{ijk} = \bt_{ij}\,^{p}h_{pk}$ for the unique distinguished metric $h \in [h]$ such that $\uR_{h} = -2a^{-1}$. Alternatively, $||B||_{h}^{2}$ does not depend on the choice of $h \in [h]$ and the choice of $h$ is determined by requiring that $2\vol_{h}(M) = a(4^{-1}||B||_{h}^{2} - 4\pi\chi(M))$.

Evidently $\prebij^{a}$ is $\diff^{+}(M)$ equivariant, so descends to a  $\map^{+}(M)$-equivariant map $\bij^{a}:\tmod(M) \to \cubic(M)$, which covers the identity on $\teich(M)$. It is convenient to write $\bij^{a}(\en, [h]) = [\prebij^{a}(\en, [h])]$ for the image of an equivalence class $[\en, [h]] \in \tmod(M)$ (the notation indicating the equivalence class of $(\en, [h])$ was omitted on the lefthand side).

Running through the definitions shows that, for $r \in \reap$, $\prebij^{ra} = (J, rB)$ if $\prebij^{a} = (J, B)$, which in part accounts for using $a^{-1}$ in the definition of $\prebij^{a}$.  Since $(\prebij^{r})^{-1}(J, B) = (\prebij^{1})^{-1}((J, r^{-1}B))$ it in general suffices to work with the fixed map $\prebij \defeq \prebij^{1}$ and its inverse. The duality on $\epremod(M)$ given by conjugacy of AH structures corresponds under $\prebij$ to replacing the holomorphic cubic differential $B$ by $-B$, in the sense that $\prebij(\ben, [h]) = (J, -B)$ if $\prebij(\en, [h]) = (J, B)$.

\begin{theorem}\label{2dmodulitheorem}
For a compact, oriented surface $M$ of genus $g > 1$, the $\diff^{+}(M)$-equivariant map $\prebij:\epremod(M) \to \precubic(M)$ is a bijection, and so the $\map^{+}(M)$-equivariant map $\bij:\tmod(M) \to \cubic(M)$ is a bijection as well. 
\end{theorem} 

\begin{proof}
For $(J, B) \in \precubic(M)$ let $[h]$ be the conformal structure determined by $J$. For the unique representative $\tilde{h} \in [h]$ having constant scalar curvature $-2$, there is by Lemma \ref{wexistencelemma} a unique solution $\phi$ to $\ahop(\tilde{h}, -2, B)(\phi) = 0$. The metric $h = e^{\phi}\tilde{h} \in [h]$ is the distinguished representative of an exact Einstein AH structure $(\en, [h])$ having aligned representative $\nabla = D - \tfrac{1}{2}h^{kp}B_{ijp}$ and $\uR_{h} = -2$, so that $\prebij(\en, [h]) = (J, B)$. This shows that $\prebij$ is onto.

Now suppose that $\prebij(\en, [h]) = (J, B) = \prebij(\ten, [g])$. Since $[g]$ and $[h]$ induce the same complex structure $J$, they are equal. Let $\nabla \in \en$ and $\tnabla \in \ten$ be the aligned representatives, and let $h \in [h]$ and $\bar{h} \in [h]$ be the distinguished represenatives of $(\en, [h])$ and $(\ten, [h])$ such that the corresponding curvatures $\uR_{h}$ and $\tilde{\uR}_{\bar{h}}$ equal $-2$. Let $\tilde{h} \in [h]$ be the unique representative with constant scalar curvature $\sR_{\tilde{h}} = -2$ and write $h_{ij} = e^{\phi}\tilde{h}_{ij}$ and $\bar{h}_{ij} = e^{\bar{\phi}}\tilde{h}_{ij}$. Then $\phi$ and $\bar{\phi}$ solve $\ahop(\tilde{h}, -2, B)(\psi) = 0$, and so, by Lemma \ref{wuniqlemma}, $\bar{\phi} = \phi$. This shows that $\bar{h}_{ij} = h_{ij}$. It follows that $\tnabla = \nabla$, and so $\ten = \en$, and the injectivity of $\prebij$ has been shown.
\end{proof}

The content of the surjectivity statement in Theorem \ref{2dmodulitheorem} is not essentially different from that of Theorem $3.4$ of C.~P. Wang's \cite{Wang}, in which it is shown how to construct an affine hypersphere from a conformal metric and a cubic holomorphic differential. As is explained in section \ref{convexsection}, $\tmod(M)$ can be identified with the deformation space $\dproj(M)$ of convex flat real projective structures on $M$. The key point is that an exact Einstein AH structure is determined by its underlying flat projective structure. The resulting identification of $\cubic(M)$ with $\dproj(M)$ is due independently to F. Labourie and J. Loftin.

\subsection{}\label{gaugesection}
The results of this section can be viewed as preliminary steps in the direction of understanding the action of $GL^{+}(2, \rea)$ on $\tmod(M)$ described in section \ref{singularmetricsection}. In the remainder of the section $M$ is a compact, oriented surface of genus at least two and $B\in \Ga(S^{3}_{0}(\ctm))$ is the real part of a cubic holomorphic differential supposed not identically zero. Sense can be made of all the results of this section for $1 \leq p \in \integer$ in place of $3$ if Gauduchon metrics are interpreted as the representatives given by Corollary \ref{vortexcorollary} and the numbers $3$ and $2$ are replaced by $p$ and $(p-1)$, where appropriate.

Specializing the action of $GL^{+}(2, \rea)$ on $\precubic(M)$ to $\comt = C^{+}O(2)$, write $z \cdot (J, B^{(3,0)}) = (z \cdot J, z \cdot B^{(3, 0)})$. Then, for $z = re^{\j \theta} = e^{t}e^{\j\theta} \in\comt$, $z \cdot J = J$ and $2\re(z \cdot B^{(3,0)}) = 2\re(z^{3}B^{(3,0)}) = e^{3t}(\cos (3\theta) B + \sin (3\theta) \bj(B))$. Let $(\en,[h]) = \prebij^{-1}(J, B)$ and $(\ten, [h]) = \prebij^{-1}(J, B(t, \theta))$ where $B(t, \theta) = e^{3t}(\cos(\theta)B + \sin(\theta)\bj(B))$ for some fixed $t > 0$ and $\theta$. Retaining the factor of $3$ in the scale factor, while dropping it in the rotational part simplifies formulas appearing later. Both $(\en, [h])$ and $(\ten, [h])$ are exact Einstein AH structures. Let $h \in [h]$ and $\tilde{h} \in [h]$ be the respective distinguished metrics such that $\uR_{h} = -2 = \uR_{\tilde{h}}$, and write $\tilde{h}_{ij} = e^{\phi}h_{ij}$. By construction the cubic torsions are $\bt_{ij}\,^{k} = h^{kp}B_{ijp}$ and 
\begin{align*}
\begin{split}
\tilde{\bt}_{ij}\,^{k} &=  \tilde{h}^{kp}(e^{t}e^{\j\theta} \cdot B)_{ijp} = e^{3t-\phi}h^{kp}(\cos (\theta) B_{ijp} + \sin (\theta) \bj(B)_{ij}\,^{p})\\
& = e^{3t-\phi}(\cos (\theta) \bt_{ij}\,^{k} + \sin (\theta) J_{i}\,^{p}\bt_{pj}\,^{k}),
\end{split}
\end{align*}
(recall \eqref{btj} of Lemma \ref{2dcomplexlemma}). By construction $\phi$, and so also the Gauduchon representative of $\prebij^{-1}(\ten, [h])$, does not depend on $\theta$, so in analyzing the dependence on $t$ it may as well be assumed that $\theta = 0$. The dependence on $t$ was partly analyzed in Proposition $1$ of \cite{Loftin-cubic}, and that result, as well as similar ones for quadratic differentials in the context of Teich\"uller space proved in section $5$ of \cite{Wolf-teichmuller}, provided motivation for Lemma \ref{phigrowthlemma} and Theorem \ref{liptheorem}.

Since, by Theorem \ref{classtheorem}, $\sR_{h} = -2 + 4^{-1}|B|_{h}^{2}$, by construction $\phi$ solves
\begin{align}\label{operphi}
\oper(\phi) \defeq \ahop(h, -2, e^{3t}B)(\phi) = \lap_{h}\phi + 2 - 2e^{\phi} + 4^{-1}(e^{2(3t-\phi)} - 1)|B|_{h}^{2} = 0.
\end{align}
\begin{lemma}\label{phigrowthlemma}
Let $(J, B) \in \precubic(M)$ for a compact smooth orientable surface $M$ of genus $g > 1$. Let $h$ and $\tilde{h} = e^{\phi}h$ be the Gauduchon metrics of $\prebij^{-1}(J, B)$ and $\prebij^{-1}(J, B(t, \theta))$ such that $\uR_{h} = -2 = \uR_{\tilde{h}}$. Then $\phi$ is the unique solution of \eqref{operphi}, and satisfies
\begin{align}\label{phibound}
&\max\{0, 2t + (2/3)\log|B|_{h} - \log 2\} \leq \phi \leq 2t,& & \text{if}\,\, t \geq 0.
\end{align}
(In particular $\phi$ is identically $0$ if $t = 0$).
\end{lemma}
\begin{proof}
By \eqref{calabibound}, $|B|_{h}^{2} \leq 8$, and for $t > 0$ there results
\begin{align}\label{ologr}
\begin{split}
\oper(2t) &= 2 - 2e^{2t} + 4^{-1}(e^{2t} - 1)|B|_{h}^{2} = 4^{-1}(e^{2t} - 1)(|B|_{h}^{2} - 8) \leq 0,   
\end{split}
\end{align}
showing that the constant function $2t$ is a supersolution of $\oper$. Since $\oper(0) \geq 0$ also, there exists a solution $\psi$ of $\oper(\psi) = 0$ such that $0 \leq \psi \leq 2t$. If $\varphi$ is a second solution then 
\begin{align}\label{omax}
\lap_{h}(\varphi - \psi)^{2} = 2|d(\varphi - \psi)|^{2} +  4(\varphi - \psi)(e^{\varphi} - e^{\psi}) - 2^{-1}e^{6t}(\varphi - \psi)(e^{-2\varphi} - e^{-2\psi}) \geq 0.
\end{align}
By the maximum principle there is a constant $c$ such that $\varphi - \psi = c$, and in \eqref{omax} this yields $8c(e^{c} - 1)e^{3\varphi} = ce^{6t}(e^{-2c} - 1)|B|^{2}_{h}$. Since $B$ is not identically zero this can be only if $c = 0$. Thus $\phi$ is the unique solution of $\oper(\phi) = 0$, and $0 \leq \phi \leq 2t$. By \eqref{ne4d}, the function $\psi = \log(2^{-1}e^{2t}|B|_{h}^{2/3})$ satisfies $\oper(\psi) = 0$ on the complement $M^{\ast}$ of the zero set of $B$. Since $\psi$ tends to $-\infty$ at the zeroes of $B$ and $\phi$ is bounded on $M$, the closure of the set $U = \{x \in M^{\ast}:\psi(x) > \phi(x)\}$ is contained properly in $M^{\ast}$. On $U$ there holds
\begin{align*}
\lap_{h}(\psi - \phi) = 2(e^{\psi} - e^{\phi}) - 4^{-1}e^{6t}(e^{-2\psi} - e^{-\phi})|B|_{h}^{2} \geq 0,
\end{align*}
which by the maximum principle contradicts that $U$ be non-empty. This proves \eqref{phibound}. 
\end{proof}

\begin{theorem}\label{liptheorem}
Let $M$ be a smooth compact orientable surface of genus $g > 1$. Fix $(J, B) \in \precubic(M)$ such that $B$ is not identically zero. For $t \geq 0$, let $\pwr{t}{h}$ be the Gauduchon metric of $\prebij^{-1}(J, e^{3t}B)$ such that $\uR_{\pwr{t}{h}} = -2$. Write $h = \pwr{1}{h}$ and $\pwr{t}{h} = e^{\phi_{t}}h$. Then $\phi_{t}$ is pointwise non-decreasing and Lipschitz as a function of $t \in(0, \infty)$, with Lipschitz constant $2$.
\end{theorem}
\begin{proof}
Applying Lemma \ref{phigrowthlemma} with $t = t_{2} - t_{1}$ and $\phi = \phi_{t_{2}} - \phi_{t_{1}}$, so that $\pwr{t_{2}}{h} = e^{\phi}\,\pwr{t_{1}}{h}$, yields
\begin{align}\label{lipt1}
\max\{0, 2(t_{2} - t_{1}) + (2/3)\log|e^{3t_{1}}B|_{e^{\phi_{t_{1}}}h} - \log 2\} \leq \phi_{t_{2}} - \phi_{t_{1}} \leq 2(t_{2} - t_{1}).
\end{align}
which, after simplifying and rearranging terms is
\begin{align}\label{lipest} 
\max\{\phi_{t_{1}}, 2t_{2} + \tfrac{2}{3}\log |B|_{h} - \log 2\} \leq \phi_{t_{2}} \leq \phi_{t_{1}} + 2(t_{2} - t_{1}).
\end{align}
The first inequality of \eqref{lipest} shows $\phi_{t_{2}} \geq \phi_{t_{1}}$, so $\phi_{t}$ is non-decreasing for $t \in (0, \infty)$. From the second inequality of \eqref{lipest} it follows that $0 \leq \phi_{t_{2}} - \phi_{t_{1}} \leq 2|t_{2} - t_{1}|$ for $t_{1}, t_{2} \in (0, \infty)$. 
\end{proof}

\begin{corollary}
Let $M$ be a smooth compact orientable surface of genus $g > 1$. Fix $(J, B) \in \precubic(M)$ such that $B$ is not identically zero. For $t \geq 0$ let $\pwr{t}{h}$ be the Gauduchon metric of $\prebij^{-1}(J, e^{3t}B)$ such that $\uR_{\pwr{t}{h}} = -2$ and let $\pwr{t}{g} = (e^{-2t})\pwr{t}{h}$ (which is also a Gauduchon metric for $\prebij^{-1}(J, e^{3t}B)$). Then the limits $\lim_{t\to \infty}\vol_{\pwr{t}{g}}(M)$ and $\lim_{t\to \infty}||B||_{\pwr{t}{g}}^{2}$ exist and satisfy
\begin{align}\label{rescaledlimits}
0 < \tfrac{1}{2}\vol_{\sth}(M) = \tfrac{1}{2}\int_{M}|B|^{2/3}_{h} \,d\vol_{h} \leq \lim_{t\to \infty}\vol_{\pwr{t}{g}}(M) = \tfrac{1}{8}\lim_{t\to \infty}||B||_{\pwr{t}{g}}^{2} \leq \vol_{h}(M).
\end{align}
Here $\sth = |B|^{2/3}_{h}h$. The curvature $\sR_{\pwr{r}{g}}$ satisfies $0 \geq \sR_{\pwr{t}{g}} \geq 4^{-1}e^{2t}(|B|^{2}_{h} - 8)$.
\end{corollary}
\begin{proof}
Write $h = \pwr{1}{h}$ and $\pwr{t}{h} = e^{\phi_{t}}h$. From the second inequality of \eqref{lipt1} of the proof of Theorem \ref{liptheorem} it is immediate that $\vol_{\pwr{t}{g}}(M)$ is non-increasing for $t>0$, and by \eqref{phibound} of Lemma \ref{phigrowthlemma}, $\vol_{h}(M) \geq e^{-2t}\vol_{\pwr{t}{h}}(M) = \vol_{\pwr{t}{g}}(M) \geq \tfrac{1}{2}\int_{M}|B|^{2/3}_{h} \,d\vol_{h}$. This shows the existence of $\lim_{t\to \infty}\vol_{\pwr{t}{g}}(M)$ and the bounds of \eqref{rescaledlimits}. From \eqref{vortex} it follows that $\vol_{\pwr{t}{g}}(M) - 8^{-1}||B||_{\pwr{r}{g}}^{2} = 4\pi(g-1)e^{-2t}$, from which the the equality of the limits in \eqref{rescaledlimits} is apparent. The final estimate on the curvature follows from $e^{-2t}\sR_{\pwr{t}{g}} =\sR_{\pwr{t}{h}} = -2 + 4^{-1}e^{6t}|B|^{2}_{\pwr{r}{h}} = -2 + 4^{-1}e^{6t}e^{-3\phi_{t}}|B|_{h}^{2}$ coupled with \eqref{phibound}.
\end{proof}

Lemma \ref{continuitylemma} does not appear to be particularly useful, but it illustrates how the bounds of Lemma \ref{factorboundlemma} can be used, and confirms a naive expectation.
\begin{lemma}\label{continuitylemma}
Let $M$ be a compact orientable smooth surface. Fix a complex structure $J$ and equip the space of holomorphic cubic differentials with the sup norm and $\cinf(M)$ with the sup norm. Then $\prebij^{-1}((J, \dum))$ is continuous with respect to these sup norms.
\end{lemma}
\begin{proof}
Let $\tilde{h}$ be the metric of curvature $-2$ representing the conformal structure determined by $J$. Let $B^{(3,0)}$ and $C^{(3, 0)}$ be holomorphic cubic differentials and suppose $|B - C|_{\tilde{h}}\leq \ep$ on $M$. Let $\phi = \prebij^{-1}((J, B))$ and $\psi = \prebij^{-1}((J, C))$. By construction $4\lap_{\tilde{h}}(\psi - \phi) = 8(e^{\psi} - e^{\phi}) + e^{-2\phi}|B|_{\tilde{h}}^{2} - e^{-2\psi}|C|_{\tilde{h}}^{2}$. Suppose $(\psi - \phi)^{2}$ assumes a positive maximum at $p \in M$. Then $|\psi - \phi|$ also assumes a positive maximum at $p$. Without loss of generality it may be supposed that $\psi(p) > \phi(p)$, so that $\max_{M}|\psi - \phi| = \psi(p) - \phi(p)$. At $p$ there holds
\begin{align}\label{cont1}
\begin{split}
0  \geq 2(\psi - \phi)^{-1}\lap_{\tilde{h}}(\psi - \phi)^{2} &\geq 8(e^{\psi} - e^{\phi}) + (e^{-2\phi}|B|^{2}_{\tilde{h}} - e^{-2\psi}|C|_{\tilde{h}}^{2})\\
& = 8(e^{\psi} - e^{\phi}) + (e^{-2\phi} - e^{-2\psi})|B|^{2}_{\tilde{h}} + e^{-2\psi}(|B|^{2}_{\tilde{h}} - |C|_{\tilde{h}}^{2}).
\end{split}
\end{align}
The inequality \eqref{cont1} forces that the value at $p$ of $|C|^{2}_{\tilde{h}}$ is greater than the value at $p$ of $|B|^{2}_{\tilde{h}}$. Using $|B|_{\tilde{h}}^{2} \leq 8$, it follows that at $p$ there holds 
\begin{align}\label{cont1b}
\begin{split}
e^{-2\psi}\ep(\ep + 16) &\geq e^{-2\psi}(\ep + 16)|C - B|_{\tilde{h}}\geq e^{-2\psi}(|C|_{\tilde{h}} - |B|_{\tilde{h}})(2|B|_{\tilde{h}} + \ep) \\
&\geq e^{-2\psi}(|C|^{2}_{\tilde{h}} - |B|_{\tilde{h}}^{2})  \geq 8(e^{\psi} - e^{\phi}) + (e^{-2\phi} - e^{-2\psi})|B|_{\tilde{h}}^{2}. 
\end{split}
\end{align}
By the second inequality of \eqref{calabibound2} of Lemma \ref{factorboundlemma}, $e^{-2\psi}|C|_{\tilde{h}}^{4/3} \leq 4$. In \eqref{cont1b} this yields that at $p$,
\begin{align}\label{cont2}
\begin{split}
4\ep(\ep + 16) & \geq 8|C|_{\tilde{h}}^{4/3}(e^{\psi} - e^{\phi}) 
\geq 8|C|_{\tilde{h}}^{4/3}(\psi(p) - \phi(p)) =8|C|_{\tilde{h}}^{4/3}\max_{M}|\psi - \phi| .
\end{split}
\end{align} 
Since, by \eqref{cont1}, $|C|_{\tilde{h}}^{2}$ is positive at $p$, this shows the claimed continuity.
\end{proof}

\section{Einstein Weyl structures on the sphere and torus}\label{spheretorussection}
In this section the deformation spaces of Einstein Weyl structures on the two sphere $\sphere$ and torus $\torus$ and some geometric properties of their members are described. In \cite{Calderbank-mobius} and \cite{Calderbank-twod} Calderbank found explicit descriptions of Einstein Weyl structures on $\sphere$ and $\torus$. While he did not explicitly address the description of the deformation spaces, all the necessary information is at least implicit in what he writes. He finds a local normal form for solutions of the two-dimensional Einstein Weyl equations, and shows that on $\sphere$ or $\torus$ the solutions are defined globally. Essentially he writes the underlying conformal structure in isothermal coordinates and uses the Killing field provided by the Gauduchon gauge to reduce the Einstein-Weyl equations to an ODE; the reduction and solution of the resulting ODE given in \cite{Calderbank-twod} is different than that given in \cite{Calderbank-mobius}. The description given here is similar to that given in \cite{Calderbank-mobius}, though the solutions that are found are given a bit more explicitly, in terms of elementary trigonometric functions (rather than elliptic functions). The scalar curvatures and Faraday curvatures are computed explicitly, and their values are related to the parameterization of the deformation space.

\subsection{}\label{twopointssection}
An $[h]$-conformal Killing field $X$ is \textbf{inessential} if there is some $h \in [h]$ for which $X$ is Killing, i.e. $\lie_{X}h = 0$, and is \textbf{essential} otherwise. By Theorem \ref{classtheorem} the Gauduchon metric dual of the Faraday primitive of the Gauduchon class of an Einstein AH structure on a compact orientable surface is inessential. 

By Lemma \ref{flatmetriclemma}, every conformal Killing vector field on $\torus$ is parallel (and so Killing) for a flat representative of the conformal structure, and so, if non-trivial, is inessential and nowhere vanishing. A flat metric on $\torus$ may be represented as that induced by the Euclidean metric on the quotient of $\rea^{2}$ by a rank two lattice $\Ga$. A non-trivial conformal Killing field is parallel in this flat metric and so can be written as a linear combination of the constant vector fields generating $\Ga$. Such a vector field is \textbf{rational} if it is a real multiple of a linear combination of the generators of the lattice with integer coefficients, and \textbf{irrational} otherwise. A conformal Killing field on $\torus$ is rational if and only if its orbits are simple closed curves, while when irrational it has a single, dense orbit.
 
\begin{lemma}\label{rationallemma}
A vector field $X$ Killing for a non-flat Riemannian metric $h$ on $\torus$ is rational.
\end{lemma}
\begin{proof}
Since $\torus$ is compact, the group $G$ of isometries of $h$ is a compact Lie group containing the one-parameter subgroup generated by $X$. Since the infinitesimal generator of any one-parameter subgroup of $G$ is $h$-Killing, it is parallel in a flat metric conformal to $h$, so has no fixed points. If $\dim G \geq 2$ then in a neighborhood of each point of $\torus$ its action generates linearly independent $h$-Killing fields, and it follows that $h$ is flat. Since $h$ is not flat, $G$ is a union of circles. Since the one-parameter subgroup generated by $X$ is connected, it must be a circle. Its non-trivial orbits must be simple closed curves, so $X$ is rational.
\end{proof}

Since the interest here is in the deformation space of Einstein Weyl structures and any two K\"ahler structures on $\sphere$ are equivalent, in considering $\sphere$ it will suffice to regard it as the Riemann sphere $\proj^{1}(\com)$ with its standard K\"ahler structure. A holomorphic vector field on $\sphere$ has either one multiplicity two zero, or two multiplicity one zeroes.

\begin{lemma}\label{spherelemma}
An inessential conformal Killing field on $\sphere$ has two zeroes.
\end{lemma}
\begin{proof}
By Lemma $0.1$ of \cite{Chen-Lu-Tian}, a non-trivial inessential conformal Killing vector field $X$ on $\sphere$ generates an isometric $S^{1}$ action on $\sphere$ fixing some zero $p$ of $X$. Since the orbit of a point $q$ distinct from $p$ but close to $p$ is a loop comprising points equidistant from $p$, the index of $X$ at $p$ is $1$, and so by the Hopf index theorem $X$ must have a second zero. 
\end{proof}

\subsection{}
On $\sphere$ not every conformal Killing field arises as the Gauduchon dual of the Faraday primitive of an Einstein Weyl structure because not every conformal Killing field is inessential. While on a torus every conformal Killing field is inessential, and every such vector field arises from a closed but non-exact Einstein Weyl structure, if such a vector field arises from a non-closed Einstein Weyl structure then by Lemma \ref{rationallemma} it must be rational. One consequence of what follows is to show that these necessary conditions on $X$ are sufficient for it to arise in this way from an Einstein Weyl structure. The remainder of the section is dedicated to analyzing when the solutions obtained in this way are equivalent modulo $\diff_{0}(M)$, and to describing explicitly the geometry of the resulting solutions.

\subsection{}\label{reducedsection}
Suppose $(M, J, [h])$ is a K\"ahler structure on a compact surface of genus $g$ equal to zero or one, $\tilde{h} \in [h]$ has constant scalar curvature $2(1-g)$, and $X \in \Ga(TM)$ is an inessential conformal Killing field. Given some background metric $g \in [h]$ it is desired to find $\phi \in \cinf(M)$ such that $h = e^{\phi}g$ is a Gauduchon metric  of an Einstein Weyl structure $(\en, [h])$ with aligned representative $\nabla = D - 2\ga_{(i}\delta_{j)}\,^{k} + h_{ij}h^{kp}\ga_{p}$, in which $D$ is the Levi-Civita connection of $h$ and $\ga_{i} = X^{p}h_{ip}$. Since $X^{i}$ is to be $h$-Killing, by Lemmas \ref{rationallemma} and \ref{spherelemma} it must generate a circle action on $M$. The initial idea is to use this circle action to reduce $\ahop(\tilde{h}, \ka, X)(\phi) = 0$ to an ODE. In the case of the torus this works, but in the case of the sphere, it is not obvious how to make the reduction, and the substitute is to work with the flat metric $\sth = |X|_{h}^{-2}h$ on the complement $M^{\ast}$ of the zero set of $X$ (when $g = 1$, the metric $\sth$ depends on the choice of $\tilde{h}$). So there will be sought $\phi \in \cinf(M^{\ast})$ such that the metric $h = e^{\phi}\sth$ extends $\cinf$ smoothly to $M$ and is otherwise as above. The difference between the torus and sphere cases is in the boundary conditions necessary for $\phi$ to extend smoothly to all of $M$. For the torus $M^{\ast} = M$, while for the sphere, $M^{\ast}$ is the complement of two points. Note that $|X|_{\sth}^{2} = 1$ and $e^{\phi} = |X|_{h}^{2}$. By Theorem \ref{classtheorem}, there must be a constant $\ka$ such that $\uR_{h} - 4|X|_{h}^{2} = \uR_{h} - 4|\ga|_{h}^{2} = \ka$. Rewriting this last equation in terms of $\sth$ using \eqref{confscal} shows that the desired Einstein Weyl structure can be found if $\phi$ solves
\begin{align}\label{torussphereeq}
\lap_{\sth} \phi + \ka e^{\phi} + 4e^{2\phi} = 0,
\end{align}
where in the spherical case \eqref{torussphereeq} is supplemented by the condition that $\phi$ extends smoothly to $M$. 

\subsection{}
The Killing property of $X$ implies $e^{\phi}d\phi(X) = \lie_{X}(|X|_{h}^{2}) = (\lie_{X}h)(X, X) = 0$, so that $d\phi(X) = 0$. Since $\lie_{X}J = 0$ there holds $[X, JX] = 0$. Let $r$ and $s$ be parameters for the flows of $X$ and $JX$, respectively; then $X = \pr_{r}$, $JX = \pr_{s}$. Since $X$ and $JX$ are complete on $M$, their flows on $M^{\ast}$ exist for all time, so $r$ and $s$ are global coordinates on the universal cover of $M^{\ast}$, which is the Euclidean plane. Note that each of $r$ and $s$ is determined only up to a translation. Since $d\phi(X) = 0$, $\phi$ is a function of $s$ alone, and since in these coordinates $\sth = dr^{2} + ds^{2}$, \eqref{torussphereeq} becomes the ODE
\begin{align}\label{torusode}
\ddot{\phi} + \ka e^{\phi} + 4e^{2\phi} = 0,
\end{align}
in which a dot indicates differentiation with respect to $s$. If $\phi$ solves \eqref{torusode} then there is a some constant $c$ such that that 
\begin{align}\label{torusode2}
c = (\dot{\phi})^{2} + 2\ka e^{\phi} + 4 e^{2\phi} = (\dot{\phi})^{2} + (2e^{\phi} + \ka/2)^{2} - \ka^{2}/4.
\end{align}
Then the function $u = e^{-\phi}$ must solve
\begin{align}\label{tode}
(\dot{u})^{2} = cu^{2} -2 \ka u - 4 = \sign(c)(\sqrt{|c|} u - \tfrac{\sign(c)\ka}{\sqrt{|c|}})^{2} - (4 + \tfrac{\ka^{2}}{c}).
\end{align}
(For the second equality of \eqref{tode}, assume $c \neq 0$). Conversely, if $u$ is a $\cinf$ positive solution of \eqref{tode} then $\phi = -\log{u}$ solves $\dot{\phi}(\ddot{\phi} + \ka e^{\phi} + 4e^{2\phi}) = 0$. There is a unique $\cinf$ smooth solution $\phi$ of \eqref{torusode} with prescribed $1$-germ at any given point. From \eqref{torusode} it follows that $4c + \ka^{2} > 0$ with equality if and only if $\phi$ is constant, equal to $\log(-\ka/4)$, in which case it must be $\ka < 0$. If $\ka \geq 0$ then $\ddot{\phi} < 0$, so any zero of $\dot{\phi}$ is isolated and, moreover, $\phi$ has no local minimum, so $M^{\ast}$ is not compact. If $\ka < 0$ and $\phi$ is not constant, then $\ka^{2} > -4c$. If $\dot{\phi}(s_{0}) = 0 = \ddot{\phi}(s_{0})$ then by \eqref{torusode} and \eqref{torusode} there holds $\ka e^{\phi(s_{0})} = c $; in particular $c < 0$. Substituting this into \eqref{torusode2} gives $c(4c + \ka^{2}) = 0$, which is a contradiction since neither $c$ nor $\ka^{2} + 4c$ is zero. Hence if $\phi$ is not constant, then, whatever is $\ka$, the zeroes of $\dot{\phi}$ are isolated, so that, by continuity, a solution of \eqref{tode} determines a solution of \eqref{torusode}. This shows that to solve \eqref{torusode} it suffices to find $\cinf$ positive solutions of \eqref{tode}. Differentiating \eqref{tode} shows that where $\dot{u}$ is not zero, $u$ solves
\begin{align}\label{tode2}
\ddot{u} = cu - \ka.
\end{align}
Since a positive $\cinf$ solution of \eqref{tode} can be written as $u = e^{-\phi}$ for some $\phi$ solving \eqref{torusode}, it follows from the isolation of the zeroes of $\dot{\phi}$ that the zeroes of $\dot{u}$ are isolated, and so $u$ is a positive $\cinf$ solution of \eqref{tode} then it solves \eqref{tode2} subject to \eqref{tode}, viewed as a constraint. 

The constant $c$ in \eqref{tode} can be interpreted as follows. If $\phi$ solves \eqref{torusode} then $\sR_{h} = \ka + 4e^{\phi}$. Since $e^{\phi}$ must extend smoothly to $M$ it assumes a maximum value on $M$, and $\sR_{h}$ assumes its maximum value on $M$ at this same point. Since $e^{\phi}$ must vanish off of $M^{\ast}$, it must assume its maximum in the interior of $M^{\ast}$, and so $\phi$ must assume a maximum in $M^{\ast}$ as well. From \eqref{torusode2} it follows that at such a point there holds 
\begin{align}\label{cnorm}
4c + \ka^{2} = (\max_{M}\sR_{h})^{2}.
\end{align} 

\subsection{}\label{cnormalizationsection}
The Einstein Weyl structure $(\en, [h])$ resulting from a solution $\phi$ of \eqref{torussphereeq} given $X$ and $\ka$ will be said to be \textit{determined by $(X, \ka, \phi)$}. Note that solutions of \eqref{torussphereeq} need not be unique, and it can occur that the same Einstein Weyl structure is determined by various triples $(X, \ka,\phi$). This possibility will now be illustrated. Let $(\en, [h])$ be determined by $(X, \ka,\phi)$. If instead of $X$ there is considered $\bar{X} = e^{-\la}X$, then $\bsth = e^{2\la}\sth$, and $\bar{\phi} = \phi - \la$ is a solution of \eqref{torussphereeq} with $\bar{\ka} = e^{-\la}\ka$ in place of $\ka$ and $\bsth$ in place of $\sth$. The resulting Gauduchon metric $\bar{h} = e^{\bar{\phi}}\bsth = e^{\la}h$ is positively homothetic to $h$; since the resulting one-form $\bar{X}^{p}\bar{h}_{ip}$ is equal to $X^{p}h_{ip}$, the solution determined by $(\bar{X}, \bar{\ka}, \bar{\phi})$ is the same as that determined by $(X, \ka, \phi)$.

The parameters $\bar{r}$ and $\bar{s}$ corresponding to $\bar{X}$ and $J\bar{X}$ are related to $r$ and $s$ by $\bar{r} = e^{\la}r$ and $\bar{s} = e^{\la}s$. The function $\bar{u}(\bar{s})$ defined by $\bar{u}(\bar{s}) \defeq e^{-\bar{\phi}}(\bar{s})$ is $\bar{u}(\bar{s}) = e^{\la}u(\bar{s}) = e^{\la}u(e^{\la}s)$ and solves $(\tfrac{\pr \bar{u}}{\pr \bar{s}})^{2} = \bar{c}\bar{u}^{2} - 2e^{-\la}\ka\bar{u} - 4$ with $\bar{c} = e^{-2\la}c$. If $u$ solves \eqref{tode} then $(X, \ka, -\log{u})$ determines an Einstein Weyl structure. The preceeding shows that the same Einstein Weyl structure is determined by the triple $(\bar{X}, \bar{\ka}, -\log{\bar{u}})$ resulting from the solution $\bar{u}$ of \eqref{tode} with $\bar{\ka}$ and $\bar{c}$ in place of $\ka$ and $c$. Thus solving \eqref{tode} for distinct values of $c$ related by a positive constant will not result in inequivalent Einstein Weyl structures. It follows that when $c \neq 0$ the value of $c$ can without loss of generality be normalized to be any given non-zero number having the same sign as has $c$. By \eqref{cnorm} the geometric meaning of such a normalization is to fix the maximum value of the scalar curvature of $h$.

\subsection{}\label{sslidesection}
The Einstein Weyl structures determined by $(X, \ka, \phi)$ and $(X, \ka, \bar{\phi})$, in which $\bar{\phi}(s) = \phi(s - a)$ need not be the same, but they are always equivalent by an orientation-preserving diffeomorphism isotopic to the identity. Let $\tau_{t}:M \to M$ be the flow of $-JX$, the fixed points of which are the zeroes of $X$. The restriction to $M^{\ast}$ of $\tau_{t}$ is given in $(r, s)$ coordinates by $\tau_{t}(r, s) = (r, s - t)$. Evidently $\tau_{t}$ is an isometry of $\sth$ preserving $X$ and satisfying $\bar{\phi} = \phi \circ \tau_{a}$. It follows that the resulting Gauduchon metrics $\bar{h} = e^{\bar{\phi}}\sth$ and $h = e^{\phi}\sth$ are related by $\bar{h} = \tau_{a}^{\ast}(h)$, and so the resulting Einstein Weyl structures determine the same point in the deformation space $\tmod(M)$. The consequence of this observation relevant in the sequel is that the parameter $s$ can always be modified by a translation, as is convenient.

\subsection{}\label{torusconstantsolutionssection}
Consider now the case of the torus. There are discussed first the constant solutions of \eqref{torusode}, and then the nonconstant ones. By Theorem \ref{scalarexacttheorem}, $\nv$ must be negative, so $\ka < 0$ in \eqref{torussphereeq}. The only relevant solutions of \eqref{tode} are those for which $u$ is non-constant, positive, and bounded from above. For each $(X, \ka)$, there is a constant solution to \eqref{torussphereeq}, namely $\phi = \log(-\ka/4)$ (equivalently, the constant function $-4/\ka$ solves \eqref{tode} with $c = -\ka^{2}/4$). Hence $h = -(\ka/4)\sth$ is simply a flat representative of the given conformal structure $[h]$. The resulting $\ga_{i}$ is $-(\ka/4)X^{p}\sth_{ip}$, and the parameter $\nv$ of the resulting Einstein Weyl structure is $\nv = \ka \vol_{h}(M) = -\ka^{2}\vol_{\sth}(M)/4$. Since, by the discussion in section \ref{cnormalizationsection}, $(e^{-\la}X, e^{-\la}\ka, \log(-\ka/4) - \la)$ determines the same Einstein Weyl structure as does $(X, \ka, \log(-\ka/4))$, there can be imposed a normalization fixing $X$, and a convenient one is to require that $\vol_{\sth}(M) = 4\pi^{2}$, so that $\nv = -\pi^{2}\ka^{2}$. The resulting Gauduchon metric having volume $4\pi^{2}$ is simply $\sth$, and the dual to $\ga_{i}$ in this metric is $Y^{i} = \sth^{ip}\ga_{p} = -(\ka/4)X^{i}$.

Replacing $X$ by $X^{\theta} = \cos \theta X + \sin \theta JX$ for any $\theta \in [0, 2\pi)$ gives rise to an Einstein Weyl structure with the same $\nv$. Since the group of conformal automorphisms of a torus is conjugate to the elliptic modular group, which acts discretely and properly discontinuously on the upper half space, which is the Teichm\"uller space of the smooth torus, the Einstein Weyl structures determined by $(X, \ka, \log(-\ka/4))$ and $(X^{\theta}, \ka, \log(-\ka/4))$ are equivalent modulo an element of $\diff_{0}(M)$ if and only if $\theta = 0$. It follows that distinct elements of the deformation space $\tmod(M)$ are obtained as $(\ka, \theta)$ varies over $(-\infty, 0)\times[0, 2\pi)$, that is as $-\ka e^{\j \theta}$ varies over $\comt$. All the possibilities can be determined in terms of a given fixed $X$ as follows. Let $\sth \in [h]$ be the flat representative of volume $4\pi^{2}$. Then to each $-\ka e^{\j \theta} \in \comt$ there corresponds an Einstein Weyl structure such that the $\sth$ dual of the Faraday primitive is $Y^{i} = -(\ka/4)X^{\theta}$, and these Einstein Weyl structures determine distinct elements of the deformation space $\tmod(M)$ (they might not be distinct in the moduli space $\emod(M)$; for instance for the square and hexagonal tori the modular group fixes some points). These Einstein Weyl structures are characterized by having weighted scalar curvature equal to zero. The preceeding shows that the trivial complex line bundle over the upper half space parameterizes the Einstein Weyl structures on a torus having zero weighted scalar curvature considered up to deformation.

An intrinsic formulation of the preceeding goes as follows. As explained in the proof of Theorem \ref{scalarexacttheorem} an Einstein Weyl structure on the torus which is closed but not exact determines a holomorphic affine connection. Namely, in this case a Gauduchon metric $h \in [h]$ with Levi-Civita connection $D$ and Faraday primitive $\ga_{i}$ is flat, and the $(1,0)$ part $\nabla^{1,0} = D^{1,0} - 2\ga^{(1,0)}$ is holomorphic because $\ga^{(1,0)}$ is holomorphic. Moreover, every holomorphic affine connection on $M$ has the form $D^{1,0} - 2\ga^{(1,0)}$ for some holomorphic one-form $\ga^{(1,0)}$ and the Levi-Civita connection $D$ of a flat metric $h \in [h]$, so that the space of holomorphic affine connections on torus with a fixed complex structure is parameterized by $H^{0}(M, \cano^{1}) \simeq \com$, the origin corresponding to $D^{1,0}$, where $D$ is the Levi-Civita connection of any flat metric representing the conformal structure determined by the given complex structure. This description of the holomorphic affine connections on a two torus seems to have been first explicitly observed by A. Vitter in \cite{Vitter}. 

The preceeding can be summarized as saying that on the smooth torus the space of Einstein Weyl structures which are closed but not exact considered up to equivalence modulo $\diff_{0}(M)$ is in bijection with the complement of the zero section in the bundle over the upper half space (the Teichm\"uller space of $M$) the fiber of which over a given conformal structure on $M$ comprises the one complex dimensional vector space of holomorphic one forms. 

\subsection{}\label{torusmetricsection}
Now there will be sought non-constant, positive, bounded solutions to \eqref{tode} on the torus. By \eqref{cnorm}, $-\ka^{2} \leq 4c$. The equation \eqref{tode2} can have positive bounded solutions only if $c < 0$, in which case the general solution is $u(s) = c^{-1}(\ka - \sqrt{\ka^{2} + 4c}\cos(\sqrt{|c|} s - \al))$, in which $\al$ is arbitrary. By the discussion in section \ref{sslidesection}, it can with no loss of generality be supposed that $\al = 0$. Thus $h = -c(\sqrt{\ka^{2} +4c}\cos(\sqrt{|c|}s) - \ka)^{-1}(dr^{2} + ds^{2})$, which has period $2\pi/\sqrt{|c|}$ in $s$.

$X$ and $JX$ are linearly independent commuting vector fields on $\torus$ preserving $\sth$, the composition of their flows defines an isometric action of $\rea^{2}$ on $\torus$ for which the stabilizers of any two points are conjugate, and the stabilizer of a given point is a discrete subgroup of $\rea^{2}$, so a lattice. The lifts to the universal cover of the flows of $X$ and $JX$ are by translations parallel to the generators of the lattice. A fundamental domain for the action of $\pi_{1}(\torus)$ is a half-closed parallelogram. If $h$ is to define a metric on $\torus$, it must be that the value of $h$ is the same at every pre-image of any $p \in \torus$. This forces the periodicity of $h$ in $s$ to be commensurate with that of the fundamental domain, so that one closed side of the fundamental domain lies on the $s$ axis, and the length of this side is a multiple of $2\pi/\sqrt{|c|}$ by a positive integer $m$. Since $\vol_{\sth}(\torus)$ is given by integrating $dr \wedge ds$ over the fundamental domain, the length of the side lying on the $r$-axis must be $\ell = \sqrt{|c|}\vol_{\sth}(\torus)/2\pi m$. Using that for $b > 1$ a primitive of $(\cos{x} + b)^{-1}$ is $\tfrac{2}{\sqrt{b^{2} - 1}}\arctan\left(\sqrt{\tfrac{b-1}{b+1}}\tan{\tfrac{x}{2}} \right)$ yields
\begin{align*}
\begin{split}
\vol_{h}(\torus) &= \int_{\torus}u^{-1}\,dvol_{\tilde{h}} 
 = \tfrac{2m\ell\sqrt{|c|}}{\sqrt{\ka^{2} + 4c}}\int_{0}^{\pi} \left(\cos t - \tfrac{\ka}{\sqrt{\ka^{2} + 4c}}\right)^{-1}\,dt 
= \pi m \ell= \sqrt{|c|}\vol_{\sth}(\torus)/2.
\end{split}
\end{align*}
By the discussion in section \ref{cnormalizationsection} fixing either of $c$ or $\vol_{\sth}(\torus)$ determines the other, and such a choice can be made as is convenient without changing the equivalence class of the Einstein Weyl structure which results; make the normalization $c = -4$, so that the metric $h$ is
\begin{align}\label{torusmetric}
h = 4\left(\sqrt{\ka^{2} - 16}\cos(2s) - \ka\right)^{-1}(dr^{2} + ds^{2}),
\end{align}
in which $\ka < -4$. By construction $\nv = \ka \vol_{h}(\torus) = \ka \vol_{\sth}(\torus)$. Computing the scalar curvature of $h$ using \eqref{conformalscalardiff} gives
\begin{align}
-\sqrt{\ka^{2} -16} = \min_{\torus}\sR_{h}  \leq \sR_{h} = \sqrt{\ka^{2} -16}\left(\frac{\ka \cos(2s) - \sqrt{\ka^{2} -16}}{\sqrt{\ka^{2} -16}\cos(2s) - \ka}\right) \leq \max_{\torus}\sR_{h} = \sqrt{\ka^{2} - 16}.
\end{align}
For $\ka = -4$ the metric $h$ corresponds to the solutions described in section \ref{torusconstantsolutionssection} having $R \equiv 0$, while for $\ka < -4$ these yield solutions have weighted scalar curvature which is somewhere positive and somewhere negative, as in \ref{rvarytorus} of Theorem \ref{scalarexacttheorem}. The metrics $\gu$ and $\hu$ homothetic to $h$ and defined by the requirements that $\vol_{\gu}(\torus) = 4\pi^{2}$ and $\max_{\torus}\sR_{\hu} = 2$, are $\gu = 4\pi^{2}\vol_{h}(\torus)^{-1}h$ and 
\begin{align*}
\begin{split}
&\hu = 2^{-1}\sqrt{\ka^{2} - 16}\,h = 2\left(\cos(2s) - \frac{\ka}{\sqrt{\ka^{2} - 16}}\right)^{-1}\left(dr^{2} + ds^{2}\right). 
\end{split}
\end{align*}
As $\nv \to -\infty$ (equivalently $\ka \to -\infty$) the metrics $\hu$ tend pointwise to the degenerate metric $\hi = 2(1 + \cos 2s)^{-1}(dr^{2} + ds^{2})$, which is the hyperbolic metric of constant scalar curvature $-2$ on its domain of non-degeneracy, which is the strip $\rea \times (-\pi/2, \pi/2)$. Since $\sth = dr^{2} + ds^{2}$ and $|X|_{\sth}^{2} = 1$, 
\begin{align}
\begin{split}
\ga &= 4(\sqrt{\ka^{2} - 16}\cos (2s) - \ka)^{-1} dr, \quad |\ga|_{h}^{2} = 4(\sqrt{\ka^{2} - 16}\cos (2s) - \ka)^{-1},\\
F & = \frac{8}{\sqrt{\ka^{2} - 16}}\frac{\sin(2s)}{(\cos(2s) - \tfrac{\ka}{\sqrt{\ka^{2} - 16}})^{2}} dr \wedge ds, \quad \fd_{h} = \frac{4\sin(2s)}{\cos(2s) - \tfrac{\ka}{\sqrt{\ka^{2} - 16}}}.
\end{split}
\end{align}
Note that $\sR_{h} - 4|\ga|_{h}^{2} = \ka$, as must be the case by construction, and that $\sR_{h}^{2} + \fd_{h}^{2} = \ka^{2} - 16$, as required by Lemma \ref{squarehololemma}. For $b > 1$ a primitive for $(\cos t + b)^{-2}$ is 
\begin{align*}
&\frac{1}{b^{2} - 1}\left(  \frac{2b}{\sqrt{b^{2} - 1}}\arctan\left(\sqrt{\frac{b-1}{b+1}}\tan\left(\frac{t}{2}\right) \right) - \frac{\sin t}{\cos t + b } \right),&
\end{align*}
so that
\begin{align*}
&\int_{0}^{\pi}\left(\cos{t} - \frac{\ka}{\sqrt{\ka^{2} - 16}}\right)^{-2}\, dt =  -\pi \ka (\ka^{2} - 16)/64,
\end{align*}
from which follows $4||\ga||_{h}^{2} = -\ka \vol_{\sth}(M) = -\nv$, as required by \ref{rvarytorus} of Theorem \ref{scalarexacttheorem}. 

This completes the analysis of the torus case. The Einstein Weyl structures on the torus which are not closed are determined up to the action of orientation-preserving diffeomorphisms by the choice of a conformal structure (a point in the hyperbolic plane) and a rational number (determining an inessential conformal Killing field). However, it is not clear how to describe the deformation/moduli space of solutions in a conceptual, geometric way (as was done in section \ref{torusconstantsolutionssection} for the closed but not exact Einstein Weyl structures).

\subsection{}
Now consider the case $g = 0$. In considering \eqref{tode} on the sphere, the only difference with the analysis in the case of the torus is that the boundary conditions are different. The requirement that $|X|_{h}^{2} = e^{\phi}$ forces that $\phi$ tends to $-\infty$ at the zeroes of $X$, or, what is the same, that $u$ tends to $+\infty$ at the zeroes of $X$. Thus the only solutions of \eqref{tode} that are relevant are those for which $u$ is positive and tends to $+\infty$ at the boundary of $\sphereast$. Where $\dot{u} \neq 0$ there holds $\ddot{u} = cu - \ka$; if $c < 0$ then the solutions are bounded, so it must be that $c \geq 0$. 

\subsection{}
For $c = 0$ the general solution of \eqref{tode} is $u(s) = -\tfrac{\ka}{2} s^{2} + a s + b$ for some constants $a$ and $b$. In this case either $u$ is equal to a positive constant, or $\ka$ is negative, for otherwise $u$ would be somewhere negative. Evaluating \eqref{tode} at $s = 0$ shows $a^{2} = -2b\ka - 4$. Since the coordinate $s$ is determined only up to translation it can be supposed that $a = 0$ by making a translation in $s$. In this case $b = - 2/\ka$, so $u = -\tfrac{\ka}{2}s^{2} - \tfrac{2}{\ka}$. Since $\sphereast$ is the complement of two points, it is topologically a cylinder, and so its universal cover is the plane with the global coordinates $r$ and $s$; in this case the metric $h = u^{-1}(dr^{2} + ds^{2})$ can descend to a metric on the cylinder obtained by quotienting by translations by $2\pi$ in the $r$ direction. This metric gives rise to an Einstein Weyl structure on the cylinder which is not exact and which has a distinguished complete metric representative $h \in [h]$ of infinite volume. Because this metric gives the cylinder infinite volume, there is no way it extends smoothly to $\sphere$. Thus these solutions of \eqref{tode} do not yield Einstein Weyl structures on the sphere.

\subsection{}\label{spheremetricoccurssection}
Now suppose that $c > 0$. The general solution of \eqref{tode} which tends to $+\infty$ as $s \to \pm\infty$ is $u(s) = c^{-1}(\sqrt{\ka^{2}+ 4c} \cosh(\sqrt{c} s + \al) + \ka)$, in which $\al$ is arbitrary. By the discussion in section \ref{sslidesection}, it can with no loss of generality be supposed that $\al = 0$. Following the discussion in section \ref{cnormalizationsection}, make the normalization $c = 4$. It will be convenient to write also $\mu = 2\ka/\sqrt{\ka^{2} + 16}$, which ranges over $(-2, 2)$ as $\ka = 4\mu/\sqrt{4 - \mu^{2}}$ ranges over $\rea$, and to introduce the coordinates $x = e^{s}\cos r$ and $y = e^{s}\sin r$, and to write $\rho = \sqrt{x^{2} + y^{2}} = e^{s}$. Then 
\begin{align}\label{spheremetric}
\begin{split}
h &= \frac{4(dr^{2} + ds^{2})}{\sqrt{\ka^{2} + 16}\cosh(2 s ) + \ka} = \frac{8}{\sqrt{\ka^{2} + 16}}\frac{dx^{2} + dy^{2}}{1 + \mu \rho^{2} + \rho^{4}}= \frac{2}{\sqrt{\ka^{2} + 16}}\frac{1 + 2\rho^{2} + \rho^{4}}{1 + \mu \rho^{2} + \rho^{4}}\tilde{h}.
\end{split}
\end{align}
(Recall $\tilde{h} = 4(1 + \rho^{2})^{-2}(dx^{2} + dy^{2})$ is the metric of scalar curvature $2$ and volume $4\pi$). In this form it is evident that $h$ extends smoothly through the origin of the $x$ and $y$ coordinates, which corresponds to the fixed point as $s \to -\infty$. Replacing $s$ by $-s$ in the preceeding gives a coordinate system in which it is evident that $h$ is smooth when $s \to \infty$. Thus the metric $h$ extends to a smooth metric on all of the two-sphere which by construction generates an Einstein Weyl structure of the sort in \ref{ps6} of Theorem \ref{scalarexacttheorem}. 

It will now be shown that for distinct values of $\ka$ the Einstein-Weyl structures so obtained are inequivalent. It is convenient to introduce the function
\begin{align}
\tau(z) = 2z\arctan(\sqrt{z^{2} + 1} - z)  = z(\pi/2 - \arctan(z)), 
\end{align}
which is a $\cinf$ orientation preserving diffeomorphism of $\rea$ onto $(-\infty, 1)$. The parameter which does not depend on the scaling of $h$ is $\nv = \ka \vol_{h}(\sphere)$. The K\"ahler form associated to $h$ is 
\begin{align}
\om =  8(\ka^{2} + 16)^{-1/2}\rho(1 + \mu\rho^{2} + \rho^{4})^{-1}dr \wedge d\rho  = dr \wedge d\left(\arctan\left(\ka/4 + \sqrt{\ka^{2}/16 + 1}\rho^{2}\right)\right),
\end{align}
and so 
\begin{align}\label{nvkavol}
\begin{split}
\nv & = \ka\vol_{h}(\sphere) = \ka\int_{\sphere}\om = 2\pi\ka\left(\pi/2 - \arctan(\ka/4)\right) = 8\pi\tau(\ka/4).
\end{split}
\end{align}
Letting $\ka$ run over $\rea$, this shows that all values of $\nv < 8\pi$ are realized by some Einstein Weyl structure on the sphere.

Henceforth, $\ka$ will be viewed as a function of $\nv$ via $\ka = 4\tau^{-1}(\nv/8\pi)$. However, because of the implicit nature of this definition, it will be convenient to continue writing $\ka$ and $\mu$ in formulas. Using \eqref{conformalscalardiff} the scalar curvature $\sR_{h}$ of $h$ can be computed by finding the Euclidean Laplacian of the conformal factor in \eqref{spheremetric}. The result is
\begin{align}\label{kasr}
\sR_{h} = \frac{\sqrt{\ka^{2} + 16}}{2}\left(\frac{\mu + 4\rho^{2} + \mu \rho^{4}}{1 + \mu\rho^{2} + \rho^{4}}\right) = \ka  + \frac{32}{\sqrt{\ka^{2} + 16}}\left(\frac{\rho^{2}}{1 + \mu\rho^{2} + \rho^{4}}\right),
\end{align}
For fixed $\mu$, $\sR_{h}$ takes its maximum value on the equatorial circle $\rho = 1$, where its value is $\sqrt{\ka^{2} + 16}$, while $\sR_{h}$ attains no minimum in the plane, tending to $\ka$ as $\rho$ tends to either $0$ or $\infty$. In particular,
\begin{align}\label{kaest}
\begin{split}
&\ka = \min_{\sphere} \sR_{h} \leq \sR_{h} \leq \max_{\sphere}\sR_{h} = \sqrt{\ka^{2} + 16}.
\end{split}
\end{align}
Observe that if $\nv$ is positive or non-negative, then the scalar curvature $\sR_{h}$ has the same property. Let $\hu$ and $\gu$ be the Gauduchon metrics of the Einstein Weyl structure $(\en, [h])$ corresponding to $\nv$ distinguished respectively by the requirements that $\max_{\sphere}\sR_{\hu} = 2$ and $\vol_{\gu}(\sphere) = 4\pi$. Then $\gu = 4\pi\vol_{h}(\sphere)^{-1}h $ while  
\begin{align}
\hu = (\sqrt{\ka^{2} +16}/2)h =  4\left(1 + \mu\rho^{2} + \rho^{4}\right)^{-1}\left(dx^{2} + dy^{2}\right) = \frac{1 + 2\rho^{2} + \rho^{4}}{1 + \mu\rho^{2} + \rho^{4}}\tilde{h},
\end{align}
Hence
\begin{align}\label{huest}
\begin{split}
&-2 \leq \mu = \min_{\sphere} \sR_{\hu} \leq \sR_{\hu} = \mu + (4 - \mu^{2})\left(\frac{\rho^{2}}{1 + \mu \rho^{2} + \rho^{4}}\right)  \leq \max_{\sphere}\sR_{\hu} = 2,\\
&4\pi  \leq \vol_{\hu}(\sphere) = \tfrac{4\pi\sqrt{\ka^{2}+16}}{\ka}\tau(\ka/4).
\end{split}
\end{align}
As $\nv \to 8\pi$ (so $\ka \to \infty$ and $\mu \to 2$) the metric $\hu$ tends pointwise to the spherical metric $\tilde{h}$. As $\nv \to -\infty$ (so $\ka \to -\infty$ and $\mu \to -2$), the volume of $\hu$ goes to $+\infty$. On either of the disks complementary in $\sphere$ to the equatorial circle $\rho = 1$, the metric $\hu$ converges pointwise to the hyperbolic metric $4(1 - \rho^{2})^{-2}(dx^{2} + dy^{2})$ of constant scalar curvature $-2$. The family $\hu$ interpolates between the spherical metric and the hyperbolic metric. The positive curvature concentrates on the equatorial circle as $\mu$ nears $-2$, while the negative curvature concentrates on the complementary disks. Precisely, for $\mu < 0$ the curvature is positive for $(- 2 + \sqrt{4 - \mu^{2}})/\mu < \rho^{2} < -(2 +\sqrt{4 - \mu^{2}})/\mu$.

By construction the vector field $X^{i}$ is $x\pr_{y} - y\pr_{x} = \pr_{r}$. Explicit expressions for $\ga_{i} = X^{p}h_{pi}$ and its norm are
\begin{align}\label{normga}
&\ga = \frac{8}{\sqrt{\ka^{2} + 16}(1 + \mu \rho^{2} + \rho^{4})}\left(xdy - y dx\right), && |\ga|_{h}^{2} = \frac{8\rho^{2}}{\sqrt{\ka^{2} + 16}(1 + \mu \rho^{2} + \rho^{4})},
\end{align}
the second of which in any case follows from $\sR_{h} - \ka = 4|\ga|_{h}^{2}$ in conjunction with \eqref{kasr}. Consequently,
\begin{align}\label{fdrho}
&F = -d\ga = \frac{16}{\sqrt{\ka^{2} + 16}}\frac{\rho^{4} - 1}{(1 + \mu \rho^{2} + \rho^{4})^{2}}dx \wedge dy, & &\fd_{h}^{2} = \frac{16(1- \rho^{4})^{2}}{(1 + \mu \rho^{2} + \rho^{4})^{2}}.
\end{align}
Equation \eqref{fdrho} shows that the Faraday curvature is non-zero except along the equatorial circle, and is largest at the poles. There can be verified the equivalent
\begin{align}\label{spheresrfh}
&\sR_{h}^{2} + \fd_{h}^{2} = \ka^{2} + 16, & &\sR_{\gu}^{2} + \fd_{\gu}^{2} = 16\mu^{-2}\tau(\ka/4)^{2},
\end{align}
confirming in this case the claim of Lemma \ref{sphereconstantlemma}. By Theorem \ref{magnetictheorem} the equatorial circle, along which vanishes $F$, must be a geodesic; its $h$-length and $\hu$-length are respectively $\pi(\sqrt{\ka^{2} + 16} - \ka)^{1/2}$ and $2\sqrt{2}\pi(\ka^{2} + 16)^{1/4}(\sqrt{\ka^{2} + 16} + \ka)^{-1/2}$.

A map associating to an Einstein Weyl structure on the sphere the real part of a holomorphic vector field is determined by a manner of choosing a normalized Gauduchon metric. For example, the metric $h$ is normalized by specifying the minimum value assumed by its scalar curvature, which is the parameter $\ka$. Similarly, $\hu$ is determined by the requirement that $\max_{\sphere}\sR_{\hu} = 2$. The associated vector field $Y^{i} = \hu^{ij}\ga_{j}$ is $\tfrac{16}{\sqrt{\ka^{2} + 16}}(x\pr_{y} - y\pr_{x})$; observe that this vector field is associated to two distinct Einstein Weyl structures, those corresponding to $\pm \ka$. The normalization like that used for surfaces of higher genus selects the metric $\gu$ having volume $4\pi$. The associated vector field $Z^{i} = \gu^{ij}\ga_{j}$ is $2\ka^{-1}\tau(\ka/4)(x\pr_{y} - y\pr_{x})$. Since $\ka \to 2\ka^{-1}\tau(\ka/4)$ is an orientation-reversing diffeomorphism of $\rea$ onto $(0, \pi/2)$, in this case there is a unique vector field associated to each Einstein Weyl structure. 

It is noted in passing that the claims about magnetic geodesics made immediately following the statement of Theorem \ref{magnetictheorem} follow straightforwardly from the explicit form of $Y^{i}$ and \eqref{normga} and \eqref{fdrho}.

In the coordinate system $(x, y)$ let $z = x + \j y$ and view $z$ as the inhomogeneous coordinate on the standard chart in $\proj^{1}(\com)$. To an element $\Phi = \begin{pmatrix} a & b \\ c & d \end{pmatrix} \in \sll(2, \com)$ associate the holomorphic vector field $X^{\Phi}$ on $\proj^{1}(\com)$ defined by $X^{\Phi}_{p} = \tfrac{d}{dt}_{| t = 0}\exp(t\Phi)p$ which in the inhomogeneous coordinate $z$ has the form $(b + (a-d)z - cz^{2})\pr_{z}$. An element $\Psi \in PSL(2, \com)$ is elliptic if and only if $(\tr \Psi)^{2}/\det \Psi \in [0, 4)$, in which case $\Psi$ is conjugate in $PSL(2, \com)$ to a unique element of the form $\exp \Psi_{\theta}$ where $\Psi_{\theta} \defeq \begin{pmatrix} \j \theta & 0 \\ 0 & -\j \theta\end{pmatrix}$ with $\theta \in (0, \pi/2)$. The real part of $X^{\Psi_{\theta}} = 2\j \theta z \pr_{z}$ is $\theta(x\pr_{y}- y\pr_{x})$. It follows that the vector field $Z$ described in the preceeding paragraph is the real part of $X^{\Psi_{2\tau(\ka/4)/\ka}}$. Since $\ka \to 2\tau(\ka/4)/\ka$ is a diffeomorphism of $\rea$ onto $(0, \pi/2)$, this shows that each Einstein Weyl structure on the sphere determines a unique conjugacy class of elliptic elements in $PSL(2, \com)$. Precisely, the vector field $Z$ is the real part of the holomorphic vector field generated by the elliptic transformation given in the inhomogeneous coordinate by $z \to e^{4\j \tau(\ka/4)/\ka}z = \tfrac{\ka + 4\j}{\sqrt{\ka^{2} + 16}}z$. The parameter $\theta$ is expressed in terms of the scale invariant parameter $\nv$ by $\theta = \nv/16\pi\tau^{-1}(\nv/8\pi)$. Using \eqref{nvkavol} it is straightforward to see that $\ka$ and $\nv$ are expressed in terms of $\theta$ by
\begin{align}\label{kanvtheta}
&\ka = 4\cot 2\theta,& &\nv = 16\pi\theta \cot 2\theta.
\end{align}
Note that no Einstein Weyl structure on the two sphere gives rise to the conjugacy class of elliptic elements corresponding to $\theta = \pi/2$. What distinguished these elliptic elements is that the associated homography leaves invariant a circle. There is an Einstein Weyl structure corresponding to $\theta = \pi/2$, namely that generated by the hyperbolic metric and corresponding to the degeneration as $\ka \to -\infty$. This picture was already described in the last paragraph of section $5$ of \cite{Calderbank-twod}. Here the equivalence problem is resolved explicitly. In particular it is shown that two Einstein Weyl structures on the sphere are equivalent if and only if their vortex parameters $\nv$ are the same, and that the corresponding extended elliptic homography can be explicitly described by a rotation by an angle explicitly expressible in terms of $\nv$, namely via \eqref{kanvtheta}.

An Einstein Weyl structure $(\en, [h])$ on the sphere having vortex parameter $\nv$ is equivalent to one represented by the Gauduchon metric $\gu$ of volume $4\pi$ where the explicit expressions for $\gu$, $\sR_{\gu}$, and $\ga$ in terms of $\theta \in [0, \pi/2)$ defined from $\nv$ by \eqref{kanvtheta} are
\begin{align}\label{sphericalreps}
\begin{split}
\gu  &= \frac{2\sin 2\theta}{\theta( 1+ 2\rho^{2}\cos 2\theta  + \rho^{4})}\left(dx^{2} + dy^{2}\right),\qquad \ga  = \frac{2\sin 2\theta}{1 + 2\rho^{2}\cos 2\theta + \rho^{4}}\left(xdy - ydx\right).\\
\sR_{\gu} & 
= \frac{\theta}{\sin 2\theta}\left(\frac{4\cos 2\theta + 8\rho^{2} + 4\rho^{4}\cos 2\theta}{1 + 2\rho^{2}\cos 2\theta  + \rho^{4}}\right).
\end{split}
\end{align}

By the uniformization theorem, $\teich(\sphere)$ is a point. Since $\map^{+}(\sphere)$ is trivial the moduli space $\emod(\sphere)$ equals the deformation space $\tmod(\sphere)$. Since under the usual definition the identity element of $PSL(2, \com)$ is not considered to be elliptic, it will be convenient to call \textbf{extended elliptic} a homography $\Psi$ that is either elliptic or the identity. For such $\Psi$ there holds $(\tr \Psi)^{2}/\det \Psi \in [0, 4]$. Note, however, that $(\tr \Psi)^{2}/\det \Psi \in [0, 4]$ need not imply that $\Psi$ is extended elliptic because while the function $(\tr \Psi)^{2}/\det \Psi$ distinguishes the conjugacy classes of non-identity elements, it does not distinguish a parabolic transformation from the identity. Let $\elp(PSL(2, \com))$ denote the space of conjugacy classes of extended elliptic elements in $PSL(2, \com)$. Suppose given an Einstein-Weyl structure $(\en, [h])$ on the two sphere which is not that generated by the uniformizing conformal structure. After a diffeomorphism it may be supposed that $[h]$ is the standard conformal structure. Associate to $(\en, [h])$ the element $\Psi \in \sll(2, \com)$ such that the $(1,0)$ part of the vector field $X^{i} = h^{ip}\ga_{p}$, where $h \in [h]$ is the Gauduchon metric with volume $4\pi$ and $\ga$ is the corresponding Faraday primitive, is equal to $X^{\Psi}$. The element $\exp\Psi$ is elliptic and conjugate to $\exp\Psi_{\theta}$ for $\theta = \nv/16\tau^{-1}(\nv/8\pi)$. Applying the same construction to the pullback of $(\en, [h])$ by an element of $PSL(2, \com)$ yields an element $\exp\Psi^{\prime}$, conjugate to $\exp\Psi_{\theta^{\prime}}$ for some $\theta^{\prime}$. However, since $(\en, [h])$ and its pullback determine the same parameter $\nv$ there holds $\theta^{\prime} = \nv/16\tau^{-1}(\nv/8\pi) = \theta$, and so $\Psi^{\prime}$ is conjugate to $\Psi$. Hence the map associating to $(\en, [h])$ the infinitesimal generator $\Psi$ descends to an injection $\emod(\sphere) \to \elp(PSL(2, \com))$; the image omits the conjugacy class of the simple inversions corresponding to $\theta = \pi/2$. The map sending $[\Psi] \in \elp(PSL(2, \com))$ to the unique $\theta \in [0, \pi/2]$ such that $4\cos^{2}\theta = (\tr \Psi)^{2}/\det \Psi$ identifies the space $\elp(PSL(2, \com))$ with the interval $[0, \pi/2]$, and the subspace $[0, \pi/2)$ corresponds to $\emod(\sphere)$. While the topologies on the spaces $\emod(\sphere)$ and $\elp(PSL(2, \com))$ have not been discussed, it makes sense to regard the interval $[0, \pi/2]$ as the compactification of $\emod(\sphere)$. Theorem \ref{spheremodulitheorem} summarizes the preceeding discussion.

\begin{theorem}\label{spheremodulitheorem}
The following spaces are in pairwise bijection
\begin{list}{(\arabic{enumi}).}{\usecounter{enumi}}
\renewcommand{\theenumi}{Step \arabic{enumi}}
\renewcommand{\theenumi}{(\arabic{enumi})}
\renewcommand{\labelenumi}{\textbf{Level \theenumi}.-}
\item The moduli space $\emod(\sphere)$ of Einstein Weyl structures on $\sphere$.
\item The space of conjugacy classes of extended elliptic elements of $PSL(2, \com)$ which leave invariant no circle. 
\item The half-open interval $[0, \pi/2)$.
\end{list}
Precisely, to $\theta \in [0, \pi/2)$ correspond the conjugacy class of the homography $z \to e^{2\j \theta}z$ and the equivalence class $\{\en, [h]\}_{\theta}$ of Einstein Weyl structures having vortex parameter $\nv = 16\pi \theta \cot 2\theta \in (-\infty, 8\pi]$. To $\theta = 0$ corresponds the Einstein Weyl structure generated by the standard round metric. The equivalence class $\{\en, [h]\}_{\theta}$ is represented by an Einstein Weyl structure $(\en, [h])_{\theta}$ for which $[h]$ is the standard conformal structure on $\proj^{1}(\com)$, and the Faraday primitive $\ga$ of the Gauduchon class and the Gauduchon metric $\gu \in [g]$ of volume $4\pi$ are as in \eqref{sphericalreps}. For every $\theta \in (0, \pi/2)$ the zero set of the Faraday curvature of the associated Einstein Weyl structure is the equatorial circle. For $(\en, [h])_{\theta}$ the Gauduchon metric $\hu \in [h]$ such that $\max_{\sphere}\sR_{\hu} = 2$ is $\hu = \tfrac{\theta}{2\sin 2\theta}\gu$. As $\theta \to 0$, $\hu$ tends pointwise to the Fubini-Study metric $\tilde{h}$, while as $\theta \to \pi/2$, the restriction to either connected component of the complement of the zero set of the Faraday curvature of the metric $\hu$ tends pointwise to the hyperbolic metric of constant scalar curvature $-2$.
\end{theorem}

\subsection{}
Differentiating the family \eqref{spheremetric} of metrics with respect to the parameter $t \in (-\pi/2, 0)$ defined by $\ka = -2\cot{2t}$ and comparing the result with \eqref{kasr} shows that the one-parameter family of metrics $h(t)$ constitutes a solution to the Ricci flow, $\tfrac{d}{dt}h = -\sR_{h}h$. The one-parameter family $s(t) = -\j h(\j t)$ of metrics on the sphere obtained from $h(t)$ after conjugation by a rotation in the parameter space is an ancient solution of the Ricci flow having positive curvature discovered in \cite{Fateev-Onofri-Zamolodchikov}, where these metrices were called \textit{sausage metrics}; they are known to mathematicians as the King-Rosenau metrics. Similarly, the family \eqref{torusmetric} of metrics on the torus is a solution of the Ricci flow with respect to the parameter $t \in (0, \infty)$ defined by $\ka = -4\coth{4t}$. It would be interesting to explain conceptually why solutions of the Ricci flow arise naturally from Einstein-Weyl structures.

\section{Convexity and Hessian metrics}
\label{hessianmetricsection}
In this section it is shown that the cone over an exact Einstein AH structure with negative scalar curvature carries particularly nice Riemannian and Lorentzian Hessian metrics, and there is proved Theorem \ref{hessiantheorem}. These constructions are used in section \ref{convexsection}, where it is explained that such an Einstein AH structure is determined by its underlying flat projective structure, and there is proved Theorem \ref{summarytheorem}.

\subsection{}\label{kahleraffinesection}
Following \cite{Cheng-Yau-realmongeampere}, a pair $(\hnabla, g)$ comprising a flat affine connection $\hnabla$ and a pseudo-Riemannian metric $g$ such that around each point there are local affine coordinates $x^{i}$ (meaning that the $dx^{i}$ are $\hnabla$-parallel) and a \textbf{potential function} $F$ such that $g_{ij} = \hnabla_{i}dF_{j}$, is called a \textbf{K\"ahler affine metric}. A K\"ahler affine metric will be said to be a \textbf{Hessian metric} if the potential function is globally defined, as will be the case in the examples constructed in what follows. (This terminological distinction between \textit{Hessian metric} and \textit{K\"ahler affine metric} is not standard).

Let $(\hnabla, g)$ be a K\"ahler affine metric. Let $\mu = dx^{1}\wedge \dots \wedge dx^{n+1}$ be the volume form determined by the choice of local affine coordinates. Then $\det_{\mu}g \defeq (\det g)/\mu^{2}$ is a function. It is not well-defined, but its logarithmic differential $d\log \det_{\mu}g$ is, because changing the choice of affine coordinates only changes $\mu$ by a constant factor. The \textbf{Ricci tensor} of a K\"ahler affine metric is defined to be $-\hnabla_{i}d\log(\det_{\mu}g)_{j}$. Note that this Ricci tensor is not in general the Ricci tensor of either $\hnabla$ or $g$. Rather it is defined in analogy with the Ricci form of a K\"ahler metric. A K\"ahler affine metric is said to be Einstein if its Ricci tensor is a constant multiple of $g$. Ricci flat K\"ahler affine metrics should be seen as real analogues of Calabi-Yau manifolds. Here they will be called \textbf{Monge-Amp\`ere metrics}, as in \cite{Kontsevich-Soibelman}.
 
\subsection{}
Let $M$ be a surface. Let $\rho:\hat{M} \to M$ be the principal $\reat$ bundle over $M$ such that the third power of $\hat{M}$ is equal to the complement $\Det \ctm \setminus \{0\}$ in $\Det \ctm$ of the zero section, viewed as a principal bundle. Let $\emf$ be the line bundle associated to $\hat{M}$ the sections of which are identified with homogeneity $1$ functions on the total space of $\hat{M}$. A section of $\emf^{3}$ is naturally viewed as a section of $\Det TM$. A section of $\emf^{k}$ will be said to have \textbf{weight $k$}. Let $R_{r}$ denote dilation in the fibers of $\hat{M}$ by $r \in \reat$, and let $\rad$ be the vector field generating the flow $R_{e^{t}}$. If $u \in \Ga(\emf^{k})$ then the associated equivariant function $\tilde{u} \in \cinf(\hat{M})$ has homogeneity $k$; in particular $d\tilde{u}(\rad) = k\tilde{u}$. On the total space of $\hat{M}$ there is a tautological $2$-form $\mu$ defined for $X, Y \in T_{\theta}\emf$ by $\mu_{\theta}(X, Y) = \theta^{3}(T\rho(\theta)(X), T\rho(\theta)(Y))$, in which $\theta^{3}$ is viewed as a $2$-form on $T_{\rho(\theta)}M$. It is straightforward to check that $\Psi \defeq d\mu$ is a volume form. 

The existence part of Theorem \ref{thomastheorem} is due to T.Y. Thomas; see \cite{Thomas}. That $M$ have dimension $2$ in Theorem \ref{thomastheorem} is unimportant, but greater generality is not needed here. It is convenient to use uppercase Latin letters as abstract indices on $\hat{M}$.

\begin{theorem}\label{thomastheorem}
Let $M$ be a smooth surface equipped with a projective structure $\en$. There is a unique torsion-free affine connection $\hnabla$ on $\hat{M}$ satisfying
\begin{enumerate}
\item $\hnabla_{I}\rad^{J} = \delta_{I}\,^{J}$
\item $\hnabla \Psi = 0$.
\item $\hnabla$ is Ricci flat.
\item $R_{r}^{\ast}(\hnabla) = \hnabla$ for all $r \in \reat$.
\item The inverse image in $\hat{M}$ of a projective geodesic of $M$ is a totally geodesic surface in $\hat{M}$ tangent to $\rad$.
\item The curvature $\hat{R}_{IJK}\,^{L}$ of $\hnabla$ satisfies $\hnabla_{Q}\hat{R}_{IJK}\,^{Q} = 0$.
\end{enumerate}
The connection $\hnabla$ is the \textbf{Thomas connection} associated to $\en$. The assignment $\en \to \hnabla$ is equivariant with respect to the action of $\diff(M)$ on $M$ in the sense that for $\phi \in \diff(M)$ there holds $\lin(\phi)^{\ast}(\hnabla) = \widehat{\phi^{\ast}(\en)}$ in which $\lin(\phi)$ is the unique principal bundle automorphism of $\hat{M}$ covering $\phi$ and preserving the tautological two-form $\mu$. Moreover, $\lin(\phi) \in \Aut(\hnabla,\Psi)$ if and only if $\phi \in \Aut(\en)$.
\end{theorem}

\begin{proof}[Indication of proof]
A principal connection on $\hat{M}$ determines a principal connection on $\Det \ctm$ and vice-versa. A torsion-free affine connection induces a principal connection on $\Det \ctm$ and hence on $\hat{M}$. It is easily seen that each principal connection $\be$ on $\hat{M}$ is induced by a unique torsion-free $\nabla$ representing $\en$. Fix a principal connection $\be$ on $\hat{M}$ and let $X \to \hat{X}$ be the horizontal lift of $X \in \Ga(TM)$ to $\hat{X} \in \Ga(\hat{M})$ determined by $\be$. Note that $[\hat{X}, \hat{Y}] = \widehat{[X, Y]} - \rho^{\ast}(\om(X, Y))\rad$, in which $\rho^{\ast}(\om) = d\be$. Let $\nabla \in \en$ be the representative determined by $\be$. Define $\hnabla$ by requiring that it be torsion-free, that it satisfy $\hnabla_{I}\rad^{J} = \delta_{I}\,^{J}$, and that for any $X, Y \in \Ga(TM)$ there hold $\hnabla_{\hat{X}}\hat{Y} = \widehat{\nabla_{X}Y} + P(X, Y)\rad$, in which $P_{ij} = -R_{(ij)} - \tfrac{1}{3}R_{[ij]} = -R_{(ij)} - \tfrac{1}{2}\om_{ij}$. That $\hnabla$ verifies all the stated conditions is straightforward. That it does not depend on the choice of $\be$ can be verified directly, but it is probably easier, and conceptually better, to deduce this as a consequence of the claimed uniqueness. Verifying the uniqueness is a bit more involved; in this regard note that it is straightforward to construct examples on $\rea^{3} \setminus\{0\} \to \sphere$ showing the necessity, for the uniqueness, of the condition $R_{r}^{\ast}(\hnabla) = \hnabla$.
\end{proof}

The curvature $\hat{R}_{IJK}\,^{L}$ of $\hnabla$ is given by
\begin{align}\label{thomascurvature}
\hat{R}_{IJK}\,^{L} = \rho^{\ast}(C)_{IJK}\rad^{L},
\end{align}
in which $C_{ijk} = 2\nabla_{[i}P_{j]k}$ is the projective Cotton tensor of $\en$. In particular $\hnabla$ is flat if and only if $\en$ is projectively flat. For later use note that if $\be$ is the principal connection on $\hat{M}$ induced by $\nabla \in \en$ then there holds
\begin{align}\label{hnablabe}
\hnabla_{I}\be_{J} = -\be_{I}\be_{J} - \rho^{\ast}(P)_{IJ},
\end{align}
in which $P_{ij}$ is the modified Ricci tensor of $\nabla$.

\subsection{}\label{mongesection}
For a section $u \in \Ga(\emf)$ there holds $\tnabla_{i}u = \nabla_{i}u + \ga_{i}u$ if $\tnabla = \nabla + 2\ga_{(i}\delta_{j)}\,^{k}$. Using this it can be verified that the operator $\B(u)_{ij} \defeq \nabla_{i}\nabla_{j}u - P_{ij}u$ is projectively invariant in the sense that it does not depend on the choice of $\nabla \in \en$. If $\tilde{u}$ is the homogeneity $1$ function on $\hat{M}$ corresponding to $u$ then the symmetric tensor $\hnabla_{I}d\tilde{u}_{J}$ has the properties that $\rad^{I}\hnabla_{I}d\tilde{u}_{J} = 0$ and that for all $X, Y \in \Ga(TM)$, $\hat{X}^{I}\hat{Y}^{J}\hnabla_{I}d\tilde{u}_{J}$ is the homogeneity $1$ function on $\hat{M}$ corresponding to $\B(u)_{ij}X^{i}Y^{j}$. Similarly, if $v = u^{2}$ (equivalently $\tilde{v} = \tilde{u}^{2}$), then
\begin{align}\label{mongelift}
\tfrac{1}{2}\hnabla_{I} d\tilde{v}_{J} = \begin{pmatrix} u\B(u)_{ij} + \nabla_{i}u\nabla_{j}u & u\nabla_{i}u\\ u \nabla_{j}u & u^{2} \end{pmatrix},
\end{align}
in which the righthand side is a notationally abusive shorthand utilizing the splitting of $T\hat{M}$ determined by $\be$ and signifying, for example, that $\hat{X}^{I}\rad^{J}\hnabla_{I}d\tilde{v}_{J}$ is the equivariant function corresponding to $2u\nabla_{X}u$. Since $\det \B(u)$ is naturally viewed as a section of $(\Det \ctm)^{2} \tensor \emf^{2} \simeq \emf^{-4}$, the operator $\monge(u) \defeq u^{4}\det \B(u)$ is function-valued. Calculating $\det \hnabla d\tilde{v}$ by applying elementary row and column operations to \eqref{mongelift} yields
\begin{align}\label{dethnablav}
\det \hnabla d\tilde{v} = 8\rho^{\ast}(\monge(u))\Psi^{2}.
\end{align}
It follows immediately that $\monge(u)$ is constant if and only if $\det \hnabla d\tilde{v}$ is $\hnabla$-parallel.

Applying \eqref{thomascurvature} and the Ricci identity yields
\begin{align}
2\hnabla_{[I}\hnabla_{J]}d\tilde{v}_{K} =  -\hat{R}_{IJK}\,^{Q}d\tilde{v}_{Q} = -2\tilde{v}\rho^{\ast}(C)_{IJK},
\end{align}
so that, if $\hnabla_{I}d\tilde{v}_{J}$ is non-degenerate, then it forms with $\hnabla$ a Hessian metric if $\en$ is projectively flat.

\subsection{}\label{hessiansection}
Let $M$ be a surface and let $(\en, [h])$ be an exact Riemannian signature Einstein AH structure with negative weighted scalar curvature $R$.  Let $u = (R/2)^{-1/3} \in \Ga(\emf)$ and $v = u^{2} = 2^{2/3}R^{-2/3}$. Because $\nabla_{i}u = 0$ for the aligned representative $\nabla \in \en$, the principal connection $\be_{I}$ induced by $\nabla$ is $\tilde{u}^{-1}d\tilde{u}_{I} = \tfrac{1}{2}\tilde{v}^{-1}d\tilde{v}_{I}$. Let $F = -\tfrac{3}{2}\log \tilde{v} = -3\log |\tilde{u}|$, which is logarithmically homogeneous in the sense that $R_{e^{t}}^{\ast}(F) = F - 3t$. Define covariant symmetric two-tensors $g_{IJ}$ and $f_{IJ}$ on $\hat{M}$ by $g_{IJ} \defeq \tfrac{1}{2}\hnabla_{I}d\tilde{v}_{J}$ and $f_{IJ} = \hnabla_{I}dF_{J}$. It will be shown that with $\nabla$ the $(1, -2)$ signature metric $g_{IJ}$ and the Riemannian metric $f_{IJ}$ form with $\hnabla$ Hessian metrics which are respectively Monge-Amp\`ere and Einstein K\"ahler affine. 

Let $\tilde{R}$ be the homogeneity $-3$ function on $\hat{M}$ corresponding to $R$, and note that $d\tilde{R}(\hat{X})$ is the homogeneous function corresponding to $\nabla_{X}R$, so equals $0$. By definition and \eqref{hnablabe} there hold
\begin{align}\label{gmetric}
\begin{split}
 g_{IJ} &= \tfrac{1}{2}\hnabla_{I}d\tilde{v}_{J} = \tilde{v}\left(\hnabla_{I}\be_{J} + 2\be_{I}\be_{J}\right) = \tilde{v}\left(\be_{I}\be_{J} - \rho^{\ast}(P)_{IJ}\right)\\
& = 2^{2/3}\tilde{R}^{-2/3}\left(\be_{I}\be_{J} + \tfrac{1}{2}\rho^{\ast}(RH)_{IJ}\right),\\
f_{IJ} & = -3\hnabla_{I}\be_{J}= -3\tilde{v}^{-1}g_{IJ} + 6\be_{I}\be_{J} = 3\rho^{\ast}(P)_{IJ} + 3\be_{I}\be_{J}.
\end{split}
\end{align}
As $RH_{ij}$ is an unweighted tensor, its pullback $\rho^{\ast}(RH)_{IJ}$ has sense. From \eqref{gmetric} it is apparent that $g_{IJ}$ has signature $(1, -2)$ and $f_{IJ}$ is Riemannian. By \eqref{projcotton} $\en$ is projectively flat, and so each of $g_{IJ}$ and $f_{IJ}$ forms with $\hnabla$ a Hessian metric. Because $\nabla_{i}u = 0$ there holds $\B(u) = - P_{ij}u = (R/2)H_{ij}u = u^{-2}H_{ij}$, and so $\monge(u) = 1$. Hence by \eqref{dethnablav} and the definition of $g_{IJ}$ there holds 
\begin{align}\label{detgpsi}
\det g = \Psi^{2}. 
\end{align}
In particular $\det g$ is $\hnabla$-parallel, so that $(\hnabla, g_{IJ})$ is a Monge-Amp\`ere Hessian metric of signature $(1, -2)$. Since $\rad^{P}g_{IP} = \tilde{v}\be_{I}$ there holds $g^{IP}\be_{P} = \tilde{v}^{-1}\rad^{I}$, and so $f_{IP}g^{JP} = -3\tilde{v}^{-1}\delta_{I}\,^{J} + 6\tilde{v}^{-1}\be_{I}\rad^{J}$, from which it follows that 
\begin{align}\label{detf}
\det f = 27 \tilde{v}^{-3}\det g = 27e^{2F}\Psi^{2},
\end{align}
so that $(\hnabla, f_{IJ})$ is a Riemannian signature Einstein Hessian metric.

\subsection{}
Alternatively, $\hnabla$ and $g_{IJ}$ generate an exact flat AH structure $([\hnabla], [g])$ with $\hnabla$ as the aligned representative and $g_{IJ}$ as a distinguished metric. The cubic torsion $\hbt_{IJ}\,^{K}$ is $\hbt_{IJ}\,^{K} = g^{KQ}\nabla_{I}g_{JQ} = g^{KQ}\nabla_{Q}g_{IJ}$. 
It is convenient to write $\hbt_{IJK} = \hbt_{IJ}\,^{Q}g_{KQ} = \nabla_{I}g_{JK} = \tfrac{1}{2}\nabla_{I}\nabla_{J}d\tilde{v}_{K}$. The Levi-Civita connection $\hD$ of $g_{IJ}$ is $\hD = \hnabla + \tfrac{1}{2}\hbt_{IJ}\,^{K}$, and the Levi-Civita connection $\tilde{D}$ of $f_{IJ}$ is
\begin{align}
\tilde{D} = \hD - 2\be_{(I}\delta_{J)}\,^{K} - \tfrac{1}{3}f_{IJ}\rad^{k} + 2\be_{I}\be_{J}\rad^{K}.
\end{align}
Since $\rad^{I}\hnabla_{I}g_{JK} = (\lie_{\rad}g)_{JK} -2g_{JK} = 0$, there holds $\hD_{I}\rad^{J} = \delta_{I}\,^{J}$, and, since $\rad^{P}f_{IP} = 3\be_{I}$, there hold $\tilde{D}_{I}\rad^{J} = 0$ and $\tilde{D}_{I}\be_{J} = 0$. 

Because $\hnabla_{I}d\tilde{v}_{J}$ is non-degenerate, and because of the form of $g_{IJ}$ it is evident that the submanifolds $M_{c^{2}} = \{p \in \hat{M}: \tilde{v}(p) = c^{2}\}$ are immersed and spacelike for $c > 0$. In particular, the induced metric $\tilde{v}^{\ast}(g)_{ij}$ is $c^{-1}h_{ij}$. Because $\rad^{I}\rad^{J}g_{IJ} = -\tilde{v}$, the vector field $N^{I} = \tilde{u}^{-1}\rad^{I}$ is a unit normal along $M_{c^{2}}$. Since $\hD_{I}N^{J} = \tilde{u}^{-1}\delta_{I}\,^{J} - \be_{I}N^{J}$, the hypersurface $M_{c^{2}}$ is totally umbilic with constant mean curvature with respect to $g_{IJ}$. Similarly, because $\rad^{I}\rad^{J}f_{IJ} = 3$, the $\tilde{D}$-parallel vector field $\tfrac{1}{\sqrt{3}}\rad^{I}$ is an $f$-unit normal to the submanifolds $M_{c^{2}}$, which are therefore totally geodesic with respect to $f_{IJ}$. The preceeding is summarized in Theorem \ref{hessiantheorem}.
\begin{theorem}\label{hessiantheorem}
On a smooth surface $M$, let $(\en, [h])$ be an exact Riemannian signature Einstein AH structure with negative weighted scalar curvature $R$. The metric $g_{IJ}$ defined in \eqref{gmetric} forms with the Thomas connection $\hnabla$ on $\hat{M}$ a Lorentzian signature Monge-Amp\`ere Hessian metric structure, such that the level surfaces of $\tilde{R}$ are smoothly immersed and, with respect to $g_{IJ}$ are spacelike, umbilic, and have constant mean curvature. The metric $f_{IJ}$ defined in \eqref{gmetric} forms with $\hnabla$ a Riemannian signature Einstein Hessian structure, such that the level surfaces of $\tilde{R}$ are totally geodesic surfaces with respect to $f_{IJ}$. 
\end{theorem}

\subsection{}\label{dustsection}
Let $\hsR_{IJK}\,^{L}$ be the curvature of $\hD$. From the flatness of $\hnabla$ and 
\begin{align}\label{skewhnablahbt}
\hnabla_{[I}\hbt_{J]K}\,^{L} = - \tfrac{1}{2}\hat{R}_{IJK}\,^{L} - \tfrac{1}{2}\hat{R}_{IJ}\,^{L}\,_{K} - \hbt_{Q[I}\,^{L}\hbt_{J]K}\,^{Q}= - \hbt_{Q[I}\,^{L}\hbt_{J]K}\,^{Q},
\end{align}
there results $\hsR_{IJK}\,^{L} = -\tfrac{1}{2}\hbt_{Q[I}\,^{L}\hbt_{J]K}\,^{Q}$, and so $\hsR_{IJ} = \tfrac{1}{4}\hbt_{IA}\,^{B}\hbt_{JB}\,^{A}$ and $\hsR_{g} = \tfrac{1}{4}|\hbt|_{g}^{2}$. Since $\hat{M}$ is $3$-dimensional the metric $g$ is conformally flat, and so its curvature is completely determined by its Ricci curvature $\hsR_{IJ}$.

Since $R\bt_{ijk}$ is an unweighted tensor, it makes sense to write $\rho^{\ast}(R\bt)_{IJK}$ for its pullback to $\hat{M}$. Differentiating \eqref{gmetric} yields
\begin{align}\label{hbt}
\hbt_{IJK} = \hnabla_{I}g_{JK} = \tilde{v}\left(\be_{I}\rho^{\ast}(RH)_{JK} + \be_{(J}\rho^{\ast}(RH)_{K)I} + \tfrac{1}{2}\hnabla_{I}\rho^{\ast}(RH)_{JK}\right).
\end{align}
Since $R$ is parallel, there holds $\hat{X}^{I}\hat{Y}^{J}\hat{Z}^{K}\hnabla_{I}\rho^{\ast}(RH)_{JK} = \rho^{\ast}(RH(\bt(X, Y), Z))$, and there hold 
\begin{align*}
\begin{split}
&\rad^{I}\hnabla_{I}\rho^{\ast}(RH)_{JK} = (\lie_{\rad}\rho^{\ast}(RH))_{JK} - 2\rho^{\ast}(RH)_{JK} = -2\rho^{\ast}(RH)_{JK},\\
&\rad^{J}\hnabla_{I}\rho^{\ast}(RH)_{JK} = - \rho^{\ast}(RH)_{IK}, 
\end{split}
\end{align*}
so that the righthand side of \eqref{hbt} is equal to $\tfrac{1}{2}\tilde{v}\rho^{\ast}(R\bt)_{IJK}$, showing
\begin{align}\label{hbtrbt}
\hbt_{IJK} = \tfrac{1}{2}\tilde{v}\rho^{\ast}(R\bt)_{IJK}.
\end{align}
There is a unique tensor $\widetilde{H}^{IJ}$ such that for any covectors $\mu_{i}$ and $\nu_{i}$ there holds $\widetilde{H}^{IJ}\rho^{\ast}(\mu)_{I}\rho^{\ast}(\nu)_{J} = \widetilde{H^{ij}\mu_{i}\nu_{j}}$ and $\widetilde{H}^{IJ}\be_{J} = 0$. Since $\tilde{R}^{-1}\widetilde{H}^{IQ}\rho^{\ast}(RH)_{QJ} = \delta_{I}\,^{J} - \be_{I}\rad^{J}$ there holds
\begin{align}\label{ginverse}
g^{IJ} = \tilde{v}^{-1}\left(\rad^{I}\rad^{J} + 2\widetilde{R}^{-1}\tilde{H}^{IJ}\right).
\end{align}
It follows that
\begin{align}\label{ghbt}
|\hbt|_{g}^{2} = \tfrac{1}{2}\tilde{v}^{-1}\rho^{\ast}(R^{-1}\nbt).
\end{align}
Using \eqref{hbtrbt}, $2\bt_{ia}\,^{b}\bt_{jb}\,^{a} = \nbt H_{ij}$, \eqref{gmetric}, \eqref{ginverse}, and \eqref{ghbt} it follows straightforwardly that
\begin{align}\label{leinstein}
\begin{split}
4\hsR_{IJ} &= \hbt_{IA}\,^{B}\hbt_{JB}\,^{A} = \tfrac{1}{4}\tilde{v}^{2}\rho^{\ast}(R\bt)_{IA}\,^{B}\rho^{\ast}(R\bt)_{JB}\,^{A} = \tfrac{1}{2}\rho^{\ast}(\nbt H)_{IJ},\\
4\hsR_{g}g_{IJ} &= |\hbt|^{2}_{g}g_{IJ}  = \rho^{\ast}(\nbt H)_{IJ} + 2\rho^{\ast}(R^{-1}\nbt)\be_{I}\be_{J},\\
\hsR_{IJ} &- \tfrac{1}{2}\hsR_{g}g_{IJ} = -\tfrac{1}{4}\rho^{\ast}(R^{-1}\nbt)\be_{I}\be_{J} = -2^{-8/3}\tilde{R}^{-1/3}\tnbt N^{\flat}_{I}N^{\flat}_{J} = \hat{T}_{IJ},
\end{split}
\end{align}
in which $N^{I} = \tilde{u}^{-1}\rad^{I}$ and $N^{\flat}_{I} = N^{Q}g_{IQ} = \tilde{u}\be_{I}$. The last equation of \eqref{leinstein} can be interpreted as an instance of the $2+1$ dimensional general relativistic Einstein equations with a stress energy tensor corresponding to a pressureless perfect fluid (a \textit{dust}), if $N^{I}$ is viewed as the velocity field of the fluid and $-2^{-8/3}\tilde{R}^{-1/3}\tnbt$ is viewed as its mass-energy density. Note that $\hat{T}_{IJ}U^{I}U^{J} \geq 0$ for all vector fields $U^{I}$ on $\hat{M}$.

\subsection{}\label{mongeamperemetricsection}
For $C > 0$ let $\hat{M}_{C} = \{x \in \hat{M}: F(x) > -\log c\}$. Define $\Psi(t)$ by
\begin{align}\label{psit}
\Psi(t) = \int^{C^{1/3}}_{e^{-t/3}}(C - r^{3})^{1/3}\,dr,
\end{align}
for $t > -\log C$, and set $\phi = \Psi(F)$. Let $\phi_{IJ} = \hnabla_{I}d\phi_{J}$. From
\begin{align*}
&\phi_{IJ}  = \dot{\Psi}(F)f_{IJ} + \ddot{\Psi}(F)F_{I}F_{J},&&
&f^{JP}\phi_{IP} = \dot{\Psi}(F)\delta_{I}\,^{P} - \ddot{\Psi}(F)F_{I}\rad^{J},
\end{align*}
it follows that $\phi_{IJ}$ is positive definite on $\hat{M}_{C}$. Noting that  $\dot{\Psi} + 3\ddot{\Psi} = (Ce^{t} - 1)^{-1}\dot{\Psi}$ and using \eqref{detf}, that
\begin{align*}
\det \phi_{IJ} = 27\dot{\Psi}(F)^{2}(\dot{\Psi}(F) + 3 \ddot{\Psi}(F))e^{2F} = 1.
\end{align*}
This shows that $\phi_{IJ}$ is a Riemannian signature Monge-Amp\`ere metric on $\hat{M}_{C}$.
\begin{theorem}\label{mongeamperetheorem}
On a smooth surface $M$, let $(\en, [h])$ be an exact Riemannian signature Einstein AH structure with negative weighted scalar curvature $R$. For each $C > 0$ the metric $\phi_{IJ} = \hnabla_{I}d\phi_{J}$ where $\phi(F) = \Psi(F)$ for $\Psi$ defined by \eqref{psit} and $F$ is the function defined in terms of $R$ as in section \ref{hessiansection} is a Riemannian signature Monge-Amp\`ere metric on $\hat{M}_{C} = \{x \in \hat{M}: F(x) > -\log c\}$.
\end{theorem}
Theorem \ref{mongeamperetheorem} is essentially the same as Proposition $1$ of the unpublished erratum \cite{Loftin-Yau-Zaslow-erratum} to \cite{Loftin-Yau-Zaslow}, and the method of construction, solving the ODE for $\Psi$ that results from requiring $\Psi(F)$ to be Monge-Amp\`ere, simply follows an example in section $2$ of \cite{Loftin-Yau-Zaslow}.

\subsection{}\label{convexsection}
Here the constructions of section \ref{mongesection} are used to show the convexity of the projective structure underlying an Einstein AH structure on a surface of genus $g > 1$. Let the notation be as in that section.

\begin{theorem}\label{convextheorem}
Let $(\en, [h])$ be an exact Riemannian signature Einstein AH structure with negative weighted scalar curvature $R$ on a smooth surface $M$. If a distinguished metric $h \in [h]$ is complete then the Riemannian metric $f_{IJ}$ on $\hat{M}$ is complete and $\en$ is a convex flat real projective structure; in particular this is the case if $M$ is compact.
\end{theorem}

\begin{proof}
Note that $\en$ is projectively flat by Lemma \ref{2deinsteinlemma}. Let $h \in [h]$ be a distinguished metric. If $h_{ij}$ is complete then $g_{ij} = 3P_{ij} = -(3R_{h}/2)h_{ij}$ is complete, and it follows from \eqref{gmetric} that on $\hat{M}$ the metric $f_{IJ}$ has the form $\rho^{\ast}(g)_{IJ} + dt_{I}dt_{J}$ where $t = F/\sqrt{3}$. It is clear from this form that $f_{IJ}$ is complete if $h_{ij}$ is. Theorem $2.1$ of \cite{Shima-Yagi} shows that if a simply connected manifold admits a complete Hessian metric then its affine developing map is a diffeomorphism onto a convex domain $\Omega$ in flat affine space. Applying this to the given structures lifted to the universal cover $\tilde{M}$ of $M$ shows that the affine developing map of the Thomas connection is a diffeomorphism onto a convex domain in $\rea^{3}$. The function $F$ is strictly convex and solves $\det \hnabla_{I}dF_{J} = 27e^{2F}\Psi$. Transferring it to the image of the affine developing map of $\hnabla$ gives on $\Omega$ a function with the same properties, and it follows (see the argument given in the remark following the statement of Theorem $1$ on pages $357-358$ of \cite{Cheng-Yau-realmongeampere}) that $\Omega$ contains no complete affine line. From the equivariance of all the preceeding with respect to the scaling action on $\hat{M}$ it follows that $\Omega$ is a convex cone containing no complete affine line and the projective developing map of $\en$ is a diffeomorphism onto the oriented projectivization of $\Omega$, which by the preceeding argument is a convex domain in the projective sphere the closure of which contains no pair of antipodal points (for if it did, the cone over it, which is $\Omega$, would contain a complete affine line). This shows that $\en$ is a convex projective structure.
\end{proof}

\begin{theorem}\label{conformaldeterminedtheorem}
If a surface $M$ admits a convex flat real projective structure $\en$ then it admits a unique Riemannian signature conformal structure $[h]$ such that $(\en, [h])$ is an exact Einstein AH structure with negative weighted scalar curvature and complete distinguished metric.
\end{theorem}
\begin{proof}
By assumption the pullback of $\hat{M}$ over the universal cover of $M$ is identified by the affine developing map of the Thomas connection of $\en$ with a convex cone $\Omega \subset \rea^{3}$ containing no complete affine line. By Theorem $4.4$ of \cite{Cheng-Yau-realmongeampere} there is a smooth function $F$ on $\Omega$ solving $\det \hnabla_{I}dF_{J} = 27e^{2F}\Psi^{2}$, tending to $+\infty$ at the boundary of $\Omega$, and such that $f_{IJ} = \hnabla_{I}dF_{J}$ is a complete Riemannian metric on the interior of $\Omega$. Passing to the tube domain (in $\com^{3}$) over $\Omega$ and applying the generalized Schwarz lemma for volume forms proved in section $1$ of \cite{Mok-Yau} it can be deduced that $e^{F}$ has positive homogeneity $-3$. Define the density $u$ on $M$ by $\tilde{u} = e^{-F/3}$ and let $h_{ij} = -u^{-1}\B(u)$. Tracing through the identifications in section \ref{mongesection} backwards, it is straightforward to check that $(\en, [h])$ is an Einstein AH structure with parallel negative scalar curvature and distinguished metric $h_{ij}$. The completeness of $h_{ij}$ follows from the splitting \eqref{gmetric} and the completeness of $f_{IJ}$, as in the proof of Theorem \ref{convextheorem}. If $[g]$ is another conformal structure with the same properties as $[h]$, let $G$ be the corresponding function on $\Omega$, which has the same properties as has $F$. Passing to the tube domain over $\Omega$ and applying the Schwarz lemma from \cite{Mok-Yau} shows that $F = G$, and so $g$ and $h$ are homothetic (here the completeness of both $g$ and $h$ is essential).
\end{proof}

Theorem \ref{conformaldeterminedtheorem} shows that an exact Einstein AH structure with negative scalar curvature and complete distinguished metric is already completely determined by its underlying (necessarily flat) projective structure, which is convex. This correspondence is evidently diffeomorphism equivariant, and so combining Theorems \ref{scalarexacttheorem}, \ref{2dmodulitheorem} and \ref{conformaldeterminedtheorem} there results Theorem \ref{summarytheorem}.

\begin{proof}[Proof of Theorem \ref{summarytheorem}]
The only thing that perhaps requires comment is the genus $1$ case. From Theorem \ref{scalarexacttheorem} it follows that an exact Einstein AH structure on a torus is Weyl if and only if it is a flat conformal structure. The analogue in this case of Theorem \ref{2dmodulitheorem} follows straightforwardly from Lemma \ref{genuszerolemma}.
\end{proof}

\subsection{Remark on tractor formalism}
Projective structures are the simplest parabolic geometries, and the powerful general machinery (see \cite{Cap-Slovak-book}) applicable to such geometries should be useful in further understanding Einstein AH structures. In particular, it should be possible to give substance to the analogy between Einstein AH structures and extremal K\"ahler metrics using the cohomological ideas of \cite{Calabi-constantcurvature} (implicitly relating certain conformal and projective BGG sequences), which have yet to be developed within the general parabolic geometry framework. 

It should be mentioned that S. Armstrong has in \cite{Armstrong-einstein} proposed a notion of Einstein structures for general parabolic geometries. Here is not the place for a careful examination of how his notion relates to that of Einstein AH structures, but it is briefly indicated how to pass from the formalism used in this paper to the tractor formalism of \cite{Cap-Gover}. The homogeneities of $\rad$, $\nabla$, $g$, and $\Psi$, are, respectively $0$, $0$, $2$, and $3$. These numbers are explained by passing to the tractor bundle $\tractor \to M$, which is the rank $3$ vector bundle over $M$ the total space of which is the quotient of $T\hat{M}$ by the action of $\reat$ via $r^{-1}TR_{r}$. Sections of $\tractor$ correspond to homogeneity $-1$ vector fields on $\hat{M}$, from which observation it is clear that $\nabla$ induces a connection on $\tractor$ (the \textbf{tractor connection}), and $g$ and $\Psi$ descend, respectively, to give a metric and a volume form on $\tractor$. The relation between the Thomas and tractor connections is explained in a bit more detail in section $3.1$ of \cite{Fox-cproj} or in section $2.3$ of \cite{Cap-Zadnik-cproj}; for background on the tractor formalism see \cite{Cap-Gover}. The vector field $\rad$ does not descend to $\tractor$, but its span does, giving a distinguished line subbundle of $\tractor$. Since a choice of a non-vanishing density on $M$ can be identified with a section of this distinguished line subbundle, it determines a splitting of $\tractor$ as a direct sum of $TM$ and the trivial line bundle. In much of the literature the objects just described appear defined in terms of such a splitting; for example the exposition in sections $2$ and $3$ of \cite{Labourie-flatprojective} is made in this way. Probably an appropriate reformulation of the Einstein AH condition would interpret the section of $S^{2}(\tractor^{\ast})$ corresponding to $g$ as harmonic.

\section{Lagrangian immersions in (para)-K\"ahler space forms}\label{parakahlersection}

\subsection{}
Let $(N, g, J)$ be a (para)-K\"ahler manifold, which means that $J_{i}\,^{j}$ is an endomorphism satisfying $J_{i}\,^{p}J_{p}\,^{j} = \ep \delta_{i}\,^{j}$, where $\ep$ is $-1$ in the K\"ahler case and $+1$ in the para-K\"ahler case; $g_{ij}$ is a metric, respectively Riemannian or split; and the tensor $\Om_{ij}= J_{i}\,^{p}g_{pj}$ is a symplectic form. The definitions for para-K\"ahler manifolds of the Ricci form, the Einstein condition, etc. are formally identical to those in the K\"ahler case. A (para)-K\"ahler manifold has \textbf{constant (para)-holomorphic sectional curvature $4c$} if its curvature has the form
\begin{align}\label{chsc}
R_{ijk}\,^{l} = 2c\left(\delta_{[i}\,^{l}g_{j]k} - \ep J_{[i}\,^{l}\Omega_{j]k} + \ep \Omega_{ij}J_{k}\,^{l} \right).
\end{align}
This is equivalent to the condition that for all $X \neq 0$ there holds $g(R(X, JX)X, JX) = 4c\ep g(X, X)^{2}$. Note that if the para-K\"ahler structure $(g, J)$ has constant para-holomorphic sectional curvature $4c$, then the para-K\"ahler structure $(-g, -J)$, which has the same underlying symplectic structure, has constant para-holomorphic sectional curvature $-4c$. An immersed submanifold of a para-K\"ahler manifold is \textbf{spacelike} if the induced metric is positive definite. While the para-K\"ahler structures $(g, J)$ and $(-g, -J)$ appear similar, their spacelike submanifolds are different.

If $\ste$ is a vector space with dual $\sted$, a flat para-K\"ahler structure $(\bG, \bJ)$ on $\ste \times \sted$ is constituted by the symplectic form $\bOmega((u,\mu), (v, \nu)) = \mu(v) - \nu(u)$ and the para-complex structure $\bJ$ equaling the identity on $\sted$ and minus the identity on $\ste$. The map $\Psi:\ste \times \sted \to \gl(n+1, \rea)^{\ast}$ given by $\lb\Psi(u, \mu), A\ra = \mu(Au)$ is the moment map for the action of the general linear group $GL(n+1, \rea)$ of $\ste$ on $\ste \times \sted$. The level sets of non-zero level of $\lb \Psi, \Id\ra$ are the pseudo-spheres of constant $\bG$ norm. Their images in the quotient of $\{ (u, \mu) \in \ste \times \sted: \mu(u) \neq 0\}$ by the action of the center of $GL(n+1, \rea)$ are the two components $\Sigma_{\pm} = \{([u], [\mu]) \in \projp(\ste) \times \projp(\sted): \pm \mu(u) < 0\}$ of the complement of the incidence correspondence in $\projp(\ste) \times \projp(\sted)$, where $\projp$ denotes oriented projectivization. The flat para-K\"ahler structure on $\ste \times \sted$ descends to $\Sigma_{\pm}$ to give the model para-K\"ahler structures of constant para-holomorphic sectional curvature. This is formally parallel to the construction of the Fubini-Study metric on $\proj^{n}(\com)$ by reduction of the flat K\"ahler structure on complex Euclidean space via the Hopf fibration and the moment map for the action of $U(n)$.

This can be understood in a more general context as follows. A \textbf{para-Hermitian symmetric space} is an affine symmetric space $G/H$ with an almost para-Hermitian structure such that the symmetries act as automorphisms of the almost para-Hermitian structure. The almost para-Hermitian structure of a para-Hermitian symmetric space is necessarily para-K\"ahler, and $G$ acts by para-K\"ahler automorphisms.  This and other basic facts about these spaces are due to S. Kaneyuki and collaborators in a series of papers, from which there results: 

\begin{theorem}[\cite{Hou-Deng-Kaneyuki-Nishiyama, Kaneyuki-compactification, Kaneyuki-Kozai}]\label{kaneyukitheorem}
Let $G$ be a connected, semisimple Lie group and $H \subset G$ a closed subgroup. The following are equivalent
\begin{enumerate}
\item $G/H$ is a homogeneous para-K\"ahler manifold.
\item $H$ is an open subgroup of a Levi subgroup of a parabolic subgroup $P$ of $G$ having abelian nilradical.
\item $G/H$ is a $G$-equivariant covering space of the adjoint orbit of a hyperbolic semisimple element of $\g$.
\end{enumerate}
Up to covering para-Hermitian symmetric spaces of semisimple Lie groups are in bijection with semsimimple graded Lie algebras $\g = \g_{-1} \oplus \g_{0}\oplus \g_{1}$ in such a way that $\g = \mathfrak{lie}(G)$ and $\g_{0} = \mathfrak{lie}(H)$. 
\end{theorem}
In the setting of Theorem \ref{kaneyukitheorem}, $G/H$ is diffeomorphic to the cotangent bundle of $G/P$, and the symplectic form on $G/H$ is the pullback of the Kostant-Kirillov symplectic form pulled back from the coadjoint orbit via some fixed multiple of the Killing form. The para-Hermitian symmetric spaces are Einstein. The proof is formally parallel to the proof that Hermitian symmetric spaces are Einstein (see Proposition $9.7$ of \cite{Kobayashi-Nomizu-2}). The scalar curvature is determined up to a scale factor determined by the choice of invariant symplectic form.

The para-Hermitian symmetric space structure on the adjoint orbit of an element of $\sll(n+1, \rea)$ generating the center of the stabilizer in $SL(n+1, \rea)$ of any element of $\Sigma_{\pm}\subset \projp(\ste) \times \projp(\sted)$ has constant non-zero para-holomorphic sectional curvature. This orbit is identified with the corresponding connected component $\Sigma_{\pm}$, and its para-Hermitian structure agrees up to constant factors with the model para-K\"ahler structure described above. As in Theorem \ref{kaneyukitheorem} the components $\Sigma_{\pm}$ are diffeomorphic to the cotangent bundle $T^{\ast}\projp(\ste)$.

\subsection{}\label{pksection}
Let $i:M \to N$ be a Lagrangian immersion in the $2n$-dimensional (para)-K\"ahler manifold $(N, g, J)$, assumed spacelike in the para-K\"ahler case. Via $J_{i}\,^{j}$ the normal bundle to $i(M)$ is identified with its tangent bundle, and this gives an identification of the second fundamental form $\Pi(X, Y)$ with the completely symmetric tensor $\Pi_{ijk}$ on $M$ defined by $\Pi(X, Y, Z) = \Om(\Pi(X, Y), Z)$ for $X, Y, Z \in \Ga(TM)$. Let $h_{ij} = i^{\ast}(g)_{ij}$ be the induced metric, let $H^{i} = h^{ia}\Pi_{apq}h^{pq}$ be the vector field dual to the mean curvature one-form (which is the one-form identified with the mean curvature vector using the (para)-K\"ahler structure), and let $B_{ijk} = \Pi_{ijk} - \tfrac{3}{n+2}H_{(i}h_{jk)}$ be the completely trace-free part of the second fundamental form. In the K\"ahler case, a lemma of P. Dazord, \cite{Dazord}, shows that $dH^{\flat}_{ij}$ is the pullback via $i$ of the Ricci form, and the same statement is true in the para-K\"ahler case, with a formally identical proof. It follows that if $g$ is Einstein, then $dH^{\flat}_{ij} = 0$; in particular this is true if $g$ has constant (para)-holomorphic sectional curvature. If $(N, g, J)$ is four dimensional and has constant (para)-holomorphic sectional curvature $4c$, the Gau\ss-Codazzi equations yield
\begin{align}\label{ckmc}
\begin{split}
&\sR_{h} - 2c - \ep |B|_{h}^{2} + \tfrac{\ep}{4}|H|_{h}^{2} = 0,\\
&4\div_{h}(B)_{ij} = \tf_{h}(\lie_{H}h)_{ij},\qquad 2D_{[i}B_{j]kl} = h_{k[i}\div_{h}(B)_{j]l} + h_{l[i}\div_{h}(B)_{j]k},
\end{split}
\end{align}
the last of which is vacuous by \eqref{2dci}. Say that the immersion is \textbf{CKMC} (has \textbf{conformal Killing mean curvature}) if $\tf_{h}(\lie_{H}h)_{ij}  = 0$. By \eqref{ckmc} and Lemma \ref{kdifferentialslemma}, the immersion is CKMC if and only if both $B^{(3,0)}$ and $H^{(1,0)}$ are holomorphic. In particular, a spacelike Lagrangian immersion of a surface of genus $g > 1$ is CKMC if and only if it has mean curvature $0$. 

Suppose now that the immersion has mean curvature $0$ so that $H^{i} = 0$, and define $\nabla = D - B_{ijp}h^{kp}$. Then $\nabla_{i}h_{jk} = 2B_{ijk}$, so $\nabla$ generates with $[h]$ an AH structure $(\en, [h])$. Since $h^{pq}\nabla_{i}h_{pq} = 0 = 2h^{pq}\nabla_{p}h_{qi}$, $\nabla$ is the aligned representative of $(\en, [h])$, which is exact, and has $h$ as a distinguished representative. From \eqref{confscal} and \eqref{ckmc} it follows that the curvature $\uR_{h}$ is
\begin{align}\label{fpk}
\uR_{h} = \sR_{h} - |B|^{2}_{h} = 2c + (\ep -1)|B|^{2}_{h}.
\end{align}
In the para-K\"ahler case $\ep = 1$, so $\uR_{h} = 2c$ is a constant. Hence the weighted scalar curvature of $(\en, [h])$ is parallel, and $(\en, [h])$ is Einstein. This proves
\begin{theorem}\label{parakahlertheorem}
On a mean curvature zero spacelike Lagrangian immersion of a surface in a four dimensional para-K\"ahler manifold of constant para-holomorphic sectional curvature $4c$ there is induced an exact Riemannian Einstein AH structure with scalar curvature $2c$.
\end{theorem}
Theorem \ref{parakahlertheorem} is the $n = 2$ special case of Theorem $8.4$ of \cite{Fox-ahs}. Theorem $4.3$ of R. Hildebrand's \cite{Hildebrand-crossratio} is essentially equivalent to Theorem \ref{parakahlertheorem}, although it is stated using different terminology.

\begin{corollary}\label{parakahlercorollary}\quad
\begin{enumerate}
\item\label{fpk1} In a $4$-dimensional para-K\"ahler manifold of constant negative para-holomorphic sectional curvature, there is no mean curvature zero spacelike Lagrangian immersion of a two sphere. 
\item\label{fpk2} In a $4$-dimensional para-K\"ahler manifold of constant positive para-holomorphic sectional curvature, there is no mean curvature zero spacelike Lagrangian immersion of a compact orientable surface of genus greater than one.
\item\label{fpk3} In a flat $4$-dimensional para-K\"ahler manifold, a mean curvature zero spacelike Lagrangian immersion of a compact orientable surface is a totally geodesic Lagrangian immersion of a flat torus. 
\item\label{fpk4} In the flat para-K\"ahler space $\ste \times \sted$, there is no mean curvature zero spacelike Lagrangian immersion of a compact orientable surface.
\end{enumerate}
\end{corollary}
\begin{proof}
By Theorem \ref{scalarexacttheorem} there cannot be an exact Einstein AH structure with negative scalar curvature on the two-sphere, there cannot be an exact Einstein AH structure with positive scalar curvature on a compact orientable surface of genus greater than one, and an exact Einstein AH structure with vanishing scalar curvature on a compact orientable surface is necessarily that generated by a flat metric on a torus. This proves \eqref{fpk1}-\eqref{fpk3}. As the geodesics in the flat para-K\'ahler structure on $\ste \times \sted$ are affine lines, which are contained in no compact subset, there can be no such immersion in $\ste \times \sted$.
\end{proof}

The claims of Corollary \ref{parakahlercorollary} are not the strongest results of this sort possible. For example, conclusion \eqref{fpk4} of Corollary \ref{parakahlercorollary} follows from the much stronger Theorem $4.2$ of J. Jost and Y.~L Xin's \cite{Jost-Xin}, which generalizes a Bernstein type theorem of J\"orgens-Calabi-Pogorelov:
\begin{theorem}[\cite{Jost-Xin}]
If the image of a mean curvature zero spacelike Lagrangian immersion in the flat para-K\"ahler space $\ste \times \sted$ is closed then it is an affine subspace.
\end{theorem}

It would be interesting to know which exact Einstein AH structures on surfaces can be realized as in Corollary \ref{parakahlercorollary} as immersed or embedded mean curvature zero Lagrangian submanifolds of para-K\"ahler space forms. In fact, there is a way to associate to an exact Einstein AH structure on a surface $M$ of genus $g > 1$ a spacelike Lagrangian immersion of its universal cover $\tilde{M}$ in a para-K\"ahler manifold of constant para-holomorphic sectional curvature. This is now sketched. It is planned to report the details elsewhere, although, since the first version of this paper appeared, Hildebrand has describe, in \cite{Hildebrand-crossratio} and \cite{Hildebrand-parakahler}, a closely related correspondence between centro-affine hypersurfaces and Lagrangian submanifolds of para-K\"ahler manifolds; in particular, Theorem $4.1$ and Corollary $4.2$ of \cite{Hildebrand-crossratio} yield conclusions very similar to those of the following paragraph. 

Because it is integrable, the horizontal subbundle of the cotangent bundle determined by a flat affine connection constitutes with the vertical subbundle a pair of transverse foliations which are Lagrangian with respect to the tautological symplectic structure, so determine on the cotangent bundle a pair of para-K\"ahler structures distinguished by the choice of the vertical or the horizontal subbundle as the $+1$ eigensubbundle of the para-complex structure. Parallel transport by the flat affine connection $\hnabla$ on $\ste$ determines an identification $T^{\ast}\ste \simeq \ste \oplus \sted$ under which the horizontal (resp. vertical) subbundle is sent to that corresponding to $\ste \times \{0\}$ (resp $\{0\} \times \sted$). Under this identification, $\bOmega$ corresponds to twice the tautological symplectic form on $T^{\ast}\ste$, the para-complex structure $\bJ$ corresponds to the choice of the the vertical subbundle of $T^{\ast}\ste$ as the $+1$-eigensubbundle, and the graph $\Ga_{\be}$ of a closed one-form $\be$ on $\ste$ is identified with a Lagrangian submanifold of $\ste \times \sted$ which is conical (preserved by positive dilations) if and only if $\be$ has homogeneity $2$. The flat para-K\"ahler structure on the cotangent bundle determined by the choice of the vertical subbundle as the $+1$ eigenbundle of the para-complex structure has the property that the pullback $\be^{\ast}(\bG)$ via the closed one-form $\be$ of the resulting para-K\"ahler metric is $2\hnabla\be$. If this is non-degenerate, then the second fundamental form of $\Ga_{\be}$ is $\hnabla\hnabla\be$, and its mean curvature is the logarithmic covariant derivative of $\det \hnabla \be$. In particular, $\Ga_{\be}$ has mean curvature zero if and only if $\det \hnabla \be$ is parallel; if $\be$ is the differential of a positive homogeneity $2$ function $\tilde{v}$, this is equivalent to $\tilde{v}$ solving the Monge-Amp\`ere equation \eqref{detgpsi}. Mean curvature zero immersed Lagrangian submanifolds of $\Sigma_{\pm}$ correspond to mean curvature zero immersed conical Lagrangian submanifolds of $\ste \times \sted$. This is formally parallel to the correspondence between minimal Lagrangian immersions in complex projective space and minimal Lagrangian cones in complex Euclidean space recounted in section $2$ of \cite{Mcintosh}. The Thomas connection of the convex flat projective structure determined by the lift to $\tilde{M}$ of the given Einstein AH structure on $M$ is identified with the restriction to a proper open convex cone in $\ste$ of the standard flat affine connection $\hnabla$ on $\ste$. On this cone there is, as in section \ref{hessiansection} and the proofs of Theorems \ref{convextheorem} and \ref{conformaldeterminedtheorem}, the positive homogeneity $2$ solution $\tilde{v}$ of the Monge-Amp\`ere equation \eqref{detgpsi}. The graph of the one-form $d\tilde{v}$ is a conical spacelike mean curvature zero Lagrangian submanifold which covers the desired mean curvature zero Lagrangian submanifold of the para-K\"ahler space form. The induced Einstein AH structure coincides with the original one on $\tilde{M}$.

\subsection{}
Given a background metric $\tilde{h}_{ij}$, a cubic holomorphic differential $B^{(3,0)}$, and a holomorphic vector field $H^{(1,0)}$, it makes sense to look for a conformal metric $h_{ij} = e^{\phi}\tilde{h}_{ij}$ such that $(h, B, H)$ are as for the induced tensors on a CKMC Lagrangian immersion in a (para)-K\"ahler manifold $(N, g, J)$. There results the equation
\begin{align}\label{ckmceq}
\lap_{\tilde{h}}\phi - \sR_{\tilde{h}} + 2ce^{\phi} - \tfrac{\ep}{4}e^{2\phi}|H|_{\tilde{h}}^{2} + \ep e^{-2\phi}|B|_{\tilde{h}}^{2} = 0.
\end{align}
The solution of the case of \eqref{ckmceq} in which $\ep = +1$ and $c < 0$ has been described for compact surfaces in section \ref{constructionsection}. In the K\"ahler ($\ep = -1$) case, \eqref{ckmceq} should be interesting for both signs of $c$, corresponding to the complex projective plane and the complex hyperbolic plane. One wonders whether for Lagrangian immersions in a complex hyperbolic $4$-manifold there is a deformation space of solutions like that for mean curvature zero surfaces in three-dimensional hyperbolic space studied by C. Taubes in \cite{Taubes}. After the first version of this paper had been completed there appeared the preprint \cite{Loftin-Mcintosh}, in which question of this sort are treated in detail for the $c = -1$ case of \eqref{ckmceq}. For example, after appropriate changes of notation, Theorem $4.1$ of \cite{Loftin-Mcintosh} is:
\begin{lemma}[\cite{Loftin-Mcintosh}]
If $M$ is a compact orientable surface of genus at least two and $\tilde{h}$ is a metric of constant scalar curvature $-2$ on $M$ then for every cubic holomorphic differential $B^{(3, 0)}$ such that there holds everywhere on $M$ the bound $|B|^{2}_{\tilde{h}} \leq 8/27$ there is a solution $\phi$ to the equation \eqref{ckmceq} with parameters $\ep = -1$ and $c = -1$ satisfying $0 \geq \phi \geq \log 2 - \log 3$.
\end{lemma}
\begin{proof}
Since the short proof is just like the proof of Lemma \ref{wexistencelemma}, it is convenient to give it here. Clearly $0$ is a supersolution of $\ahop(\phi) = \lap_{\tilde{h}}\phi + 2 - 2e^{\phi} - e^{-2\phi}|B|_{\tilde{h}}^{2}$. If $c$ is any constant, then $e^{2c}\ahop(c) \geq -2p(e^{c})$ where $p(r) = r^{3} - r^{2} + \tfrac{1}{2}\max_{M}|B|^{2}_{\tilde{h}}$. This polynomial $p$ is non-negative at $r = 0$ and has at $r = 2/3$ a local minimum at which its value is $\tfrac{1}{2}\max_{M}|B|^{2}_{\tilde{h}} - 4/27$. Hence $p$ has a positive zero if and only if $\max_{M}|B|^{2}_{\tilde{h}} \leq 8/27$, in which case its smallest positive zero $r_{1}$ is no greater than $2/3$. In this case $\ahop(\log r_{1}) \geq 0$, so $\log r_{1}$ is a negative subsolution of $\ahop$. As in the proof of Lemma \ref{wexistencelemma} this suffices to show the existence of a solution to $\ahop(\phi) = 0$ satisfying the indicated bounds. 
\end{proof}
More work has to be done to construct from such a solution a minimal Lagrangian immersion in the complex hyperbolic plane, and this is a part of what is accomplished in \cite{Loftin-Mcintosh}. A similar analysis in the para-K\"ahler case should be interesting as well.

\def\cprime{$'$} \def\cprime{$'$} \def\cprime{$'$}
  \def\polhk#1{\setbox0=\hbox{#1}{\ooalign{\hidewidth
  \lower1.5ex\hbox{`}\hidewidth\crcr\unhbox0}}} \def\cprime{$'$}
  \def\Dbar{\leavevmode\lower.6ex\hbox to 0pt{\hskip-.23ex \accent"16\hss}D}
  \def\cprime{$'$} \def\cprime{$'$} \def\cprime{$'$} \def\cprime{$'$}
  \def\cprime{$'$} \def\cprime{$'$} \def\cprime{$'$} \def\cprime{$'$}
  \def\cprime{$'$} \def\cprime{$'$} \def\cprime{$'$} \def\cprime{$'$}
  \def\dbar{\leavevmode\hbox to 0pt{\hskip.2ex \accent"16\hss}d}
  \def\cprime{$'$} \def\cprime{$'$} \def\cprime{$'$} \def\cprime{$'$}
  \def\cprime{$'$} \def\cprime{$'$} \def\cprime{$'$} \def\cprime{$'$}
  \def\cprime{$'$} \def\cprime{$'$} \def\cprime{$'$} \def\cprime{$'$}
  \def\cprime{$'$} \def\cprime{$'$} \def\cprime{$'$} \def\cprime{$'$}
  \def\cprime{$'$} \def\cprime{$'$} \def\cprime{$'$} \def\cprime{$'$}
  \def\cprime{$'$} \def\cprime{$'$} \def\cprime{$'$} \def\cprime{$'$}
  \def\cprime{$'$} \def\cprime{$'$} \def\cprime{$'$} \def\cprime{$'$}
  \def\cprime{$'$} \def\cprime{$'$} \def\cprime{$'$} \def\cprime{$'$}
  \def\cprime{$'$} \def\cprime{$'$} \def\cprime{$'$} \def\cprime{$'$}
\providecommand{\bysame}{\leavevmode\hbox to3em{\hrulefill}\thinspace}
\providecommand{\MR}{\relax\ifhmode\unskip\space\fi MR }
\providecommand{\MRhref}[2]{%
  \href{http://www.ams.org/mathscinet-getitem?mr=#1}{#2}
}
\providecommand{\href}[2]{#2}

\end{document}